\documentclass{amsart}
\title[Interpolation of numbers of Catalan type]{Interpolation of numbers of Catalan type in a local field of positive characteristic}
\author{Greg W.\ Anderson}
\address{University of Minnesota, Mpls., MN 55455}
\email{gwanders@umn.edu}
\date{March 26, 2006}
\thanks{MSC2000: 11S31, 11S80}
\newcommand{\Gb}{{\mathbb{G}}}

\newcommand{\RR}{{\mathbb{R}}}

\DeclareMathOperator{\Moore}{{\mathrm{Moore}}}
\DeclareMathOperator{\Ore}{{\mathrm{Ore}}}
\newcommand{\DDD}{{\mathfrak{D}}}

\newcommand{\sh}{{\mathrm{sh}}}
\DeclareMathOperator{\Gal}{{\mathrm{Gal}}}
\newcommand{\rigged}{{\mathcal{RVL}}}
\newcommand{\ab}{{\mathrm{ab}}}
\newcommand{\perf}{{\mathrm{perf}}}
\DeclareMathOperator{\Res}{{\mathrm{Res}}}
\DeclareMathOperator{\trace}{{\mathrm{tr}}}
\newcommand{\PPP}{{\mathcal{P}}}

\newcommand{\ord}{{\,\mathrm{ord}\,}}

\newcommand{\OO}{{\mathcal{O}}}
\newtheorem{Theorem}[subsection]{Theorem}
\newtheorem{Subtheorem}[subsubsection]{Theorem}
\newtheorem{Proposition}[subsection]{Proposition}
\newtheorem{Subproposition}[subsubsection]{Proposition}
\newtheorem{Lemma}[subsubsection]{Lemma}
\newtheorem{Corollary}[subsubsection]{Corollary}
\newcommand{\card}{{\#}}
\newcommand{\one}{{\mathbf{1}}}
\newcommand{\norm}[1]{{\Vert #1 \Vert}}
\newcommand{\LLL}{{\mathcal{L}}}
\newcommand{\MMM}{{\mathcal{M}}}
\newcommand{\ZZ}{{\mathbb{Z}}}
\newcommand{\SSS}{{\mathcal{S}}}
\newcommand{\FFF}{{\mathcal{F}}}
\newcommand{\GGG}{{\mathcal{G}}}
\newcommand{\FF}{{\mathbb{F}}}
\newcommand{\ee}{{\mathbf{e}}}
\newcommand{\CC}{{\mathbb{C}}}
\begin{document}
\begin{abstract}
Let $k$ be a locally compact topological field of positive characteristic. Let $L$ be a cocompact discrete additive subgroup of $k$. Let $U$ be an open compact additive subgroup of $k$.
Let $\ell$, $u$ and $a$ be elements of $k$, with $a$ nonzero.
We study the behavior of the product $\prod_{0\neq x\in (\ell+L)\cap a(u+U)}x$
as $a$ varies, using tools from local class field theory and harmonic analysis.  Typically ratios of such products occur as partial products grouped by degree for the infinite products representing special values of Gamma-functions for function fields. Our main result provides local confirmation for a two-variable refinement of the Stark conjecture in the function field case recently proposed by the author.
\end{abstract}
\maketitle
\tableofcontents
\section{Introduction}
\subsection{Preliminary discussion of the main objects of study}
\subsubsection{The product of lattice points in a box}
We are going to study products in
a locally compact topological field $k$ of positive characteristic of the form
\begin{equation}\label{equation:InterestingProducts}
\prod_{0\neq x\in (\ell+L)\cap a(u+U)} x,\;\;\;\;
\;\mbox{for}\;\left\{\begin{array}{l}
\mbox{$L\subset k$: cocompact discrete additive subgroup,}\\
\mbox{$U\subset k$: open compact additive subgroup,}\\
\mbox{$\ell,u\in k$ and $0\neq a\in k$,}
\end{array}\right.
\end{equation}
with the goal of understanding how such numbers depend on the parameter $a$, using tools from harmonic analysis and local class field theory.  We think of $\ell+L$ as a lattice and of $a(u+U)$ as a box which grows as $a$ grows in absolute value.  Our main result (Theorem~\ref{Theorem:Interpolability} below)
gives detailed information about (ratios of) numbers of the form (\ref{equation:InterestingProducts}) and provides local confirmation of a two-variable refinement of the Stark conjecture in the function field case recently proposed in \cite{AndersonStirling}.

\subsubsection{Theta and Catalan symbols} 
The notation used in (\ref{equation:InterestingProducts}) is inconvenient, especially when, as is usually the case in practice, we wish to consider ratios or complicated monomials in such products. We introduce a more streamlined system of notation now.
 Following the probabilists, let $\one_S:X\rightarrow \{0,1\}\subset \CC$ denote the indicator function of a subset $S\subset X$,
 i.~e., 
$\one_S(x)=1$ for $x\in S$ and $\one_S(x)=0$ otherwise. Let $\card S$ denote the cardinality of a set $S$.   Let $A^\times$ denote the group of units of a ring $A$ with unit. 
Let $p$ be the characteristic of $k$.  Let $\ZZ[1/p]$ be the result of adjoining an inverse of $p$ to $\ZZ$. 
Note that for every $x\in k^\times$ and $\nu\in \ZZ[1/p]$, the expression $x^\nu$ well-defines a nonzero element of the perfect closure $\sqrt[p^\infty]{k}$ of $k$. Given $f,g:k\rightarrow \CC$, 
let $f\otimes g:k\times k\rightarrow\CC$ be defined by $(f\otimes g)(x,y)=f(x)g(y)$. We come now to the key definition. We say that
a function $\Phi:k\times k\rightarrow \ZZ[1/p]$
is a {\em rational rigged virtual lattice} if it is a finite
$\ZZ[1/p]$-linear combination of functions of the form $\one_{\ell+ L}\otimes\one_{u+U}$, where $\ell$, $L$, $u$ and $U$
are as in (\ref{equation:InterestingProducts}). (We warn the reader that the definition of rational rigged virtual lattice made here is not quite the same as that made in the main body of the paper, but the difference is inconsequential. A similar 
warning applies to the other definitions made and applied in the introduction.)  Given
a rational rigged virtual lattice $\Phi$ 
and $a\in k^\times$, we define
$$\Theta(a,\Phi)=\sum_{x\in k}
\Phi(x,a^{-1}x)\in \ZZ[1/p],$$
$$\left(\begin{array}{c}
a\\
\Phi
\end{array}\right)_+=\prod_{x\in k^\times}
x^{\Phi(x,a^{-1}x)}\in \sqrt[p^\infty]{k}^\times.$$
We call $\Theta(\cdot,\cdot)$
the {\em theta symbol}. If the function
$\Theta(\cdot,\Phi)$ is supported in a compact subset of $k^\times$, we say that $\Phi$ is {\em proper}. For the most part our attention is going to be focused on proper rational rigged virtual lattices.
 We call $\left(\begin{subarray}{c}
\cdot\\
\cdot
\end{subarray}\right)_+$ the {\em partial Catalan symbol}. 
(Of the Catalan symbol itself we speak a bit later.)
For example, notation as in (\ref{equation:InterestingProducts}), we have
$$\Theta(a,\one_{\ell+L}\otimes \one_{u+U})=\card((\ell+L)\cap a(u+U))$$
and the expression
$$\left(\begin{array}{c}
a\\
\one_{\ell+L}\otimes \one_{u+U}
\end{array}\right)_+$$
represents the product (\ref{equation:InterestingProducts}).
Having introduced the formalism of theta and Catalan symbols we are now able to handle ratios and monomials in products of the form (\ref{equation:InterestingProducts}) in an efficient way.

\subsubsection{Rationale for the ``theta'' terminology and connection with $L$-functions}
The $L$-function evaluators $\Theta_{K/k,S}(s)$ figuring in  Tate's formulation  \cite{Tate} of the Stark conjecture  can in the function field case be related to Mellin transforms of functions of the form $\Theta(\cdot,\Phi)$  by the methods of Tate's thesis. In particular, partial zeta values at $s=0$ can be represented as values of the theta symbol  in certain situations. We do not discuss the connection of the theta symbol with $L$-function evaluators here, but see \cite[\S12]{AndersonStirling} for a discussion of the analogous connection in the global setting. 
Theta symbols and $L$-function evaluators {\em \`{a} la} Tate  also have some relationship with theta divisors.  For example, see \cite[Thm.\ 4.1.1]{AndersonAHarmonic}.
While we do not discuss $L$-functions explicitly in the paper, the desire to make our results easily applicable to the study of $L$-functions has dictated our heavy emphasis on harmonic analysis.

\subsubsection{Rationale for the ``Catalan''
terminology}
We regard the values of the partial Catalan symbol
$\left(\begin{subarray}{c}\cdot\\
\cdot
\end{subarray}\right)_+$ as analogues of the classical Catalan numbers.
The formula
$$\begin{array}{rcl}
\displaystyle\frac{1}{n+1}\left(\begin{array}{c}
2n\\
n
\end{array}\right)&=&
\frac{
\prod_{x\in \ZZ\cap(0,2n]}x}
{
\prod_{x\in \ZZ\cap(0,n]}x\cdot \prod_{x\in \ZZ\cap(0,n+1]}x}\\\\
&=&
\prod_{ x\in \RR^\times} x^{(\one_{\ZZ}\otimes\one_{(0,2n]}-\one_{\ZZ}\otimes\one_{(0,n]}-\one_{\ZZ}\otimes\one_{(0,n+1]})(x,x)}
\end{array}
$$
for the usual Catalan numbers suggests how we see the analogy. Actually the  formula 
$$\left(\begin{array}{c}
2n\\
n
\end{array}\right)=
\prod_{ x\in \RR^\times} x^{(\one_{\ZZ}\otimes\one_{(0,2]}-2\cdot\one_{\ZZ}\otimes\one_{(0,1]})(x,n^{-1}x)}
$$
for the Catalan numbers of type $B$
shows off the analogy in a slightly better way.

\subsubsection{Evaluation of partial Catalan symbols in simple examples over $\FF_q[T]$}
Let $\FF_q[T]$ be the ring
of polynomials in a variable $T$
over the field $\FF_q$ of $q$ elements.
Let $\FF_q(T)$ be the fraction field of $\FF_q[T]$.
Let us temporarily (just in this paragraph and the next) identify the local field $k$ considered in (\ref{equation:InterestingProducts}) with the completion 
$\FF_q((1/T))$ of $\FF_q(T)$ at the place $T=\infty$.
Let $(1/T)\FF_q[[1/T]]$
be the open unit ball in $\FF_q((1/T))$.
We have
\begin{equation}\label{equation:SimpleExampleOne}
\prod_{\begin{subarray}{c}
n\in \FF_q[T]\\
n:\mbox{\tiny monic}\\
\deg n=N
\end{subarray}}n
=
\left(\begin{array}{c}
T^N\\
\one_{\FF_q[T]}\otimes \one_{1+(1/T)\FF_q[[1/T]]}
\end{array}\right)_+
\end{equation}
for all integers $N\geq 0$.  The product on the left is the Carlitz analogue of the factorial of $q^N$. For background concerning the latter,
see \cite[\S4.5]{ThakurBook}. We have
 \begin{equation}\label{equation:SimpleExampleTwo}
 \prod_{\begin{subarray}{c}
n\in \FF_q[T]\\
n:\mbox{\tiny monic}\\
\deg n=N
\end{subarray}}\left(1+\frac{s}{n}\right)
=
\left(\begin{array}{c}
T^N\\
(\one_{s+\FF_q[T]}-\one_{\FF_q[T]})\otimes
\one_{1+(1/T)\FF_q[[1/T]]}
\end{array}\right)_+
\end{equation}
for all $s\in (1/T)\FF_q[[1/T]]$
and integers $N\geq 0$. The product on the left is
(for $s\neq 0$)
a (reciprocal) partial product grouped by degree for the infinite product
$$\frac{1}{s}\prod_{\begin{subarray}{c}
n\in \FF_q[T]\\
n:\mbox{\tiny monic}
\end{subarray}}\left(1+\frac{s}{n}\right)^{-1}$$
representing the value at $s$ of the geometric $\Gamma$-function over $\FF_q[T]$.
For background on $\Gamma$-functions for function fields, see \cite[Chap.\ 9]{GossBook} or \cite[Chap.\ 4]{ThakurBook}.
The general theory of $\Gamma$-functions for function fields provides a vast supply of arithmetically significant infinite products whose 
partial products grouped by degree can be represented as values
of partial Catalan symbols, and thus made accessible to study by our methods.

\subsubsection{Examples of interpolation formulas}
By a calculation with Moore determinants it is possible to rewrite (\ref{equation:SimpleExampleOne})
in the form
\begin{equation}
\label{equation:SimpleExampleOneBis}
T^{q^{N+1}}-T=\left(\begin{array}{c}
T^{N+1}\\
\one_{\FF_q[T]}\otimes ( \one_{1+(1/T)\FF_q[[1/T]]}-q
\one_{(1/T)(1+(1/T)\FF_q[[1/T]])})
\end{array}\right)_+,
\end{equation}
for all integers $N\geq 0$.
There are two things to notice here. Firstly, 
dependence of the right side on $N$ is ``explained'' by a varying exponent $q^{N+1}$ on the left. Secondly, the rational rigged virtual lattice appearing on the right is proper. 
Similarly, by a slightly more involved manipulation of Moore determinants, one verifies that
 \begin{equation}\label{equation:SimpleExampleTwoBis}
\begin{array}{cl}
&\displaystyle 1- \alpha (-T)^{-\frac{q^{N+1}-1}{q-1}}
- \beta (1-T)^{-\frac{q^{N+1}-1}{q-1}}\\\\
=&\displaystyle
\left(\begin{array}{c}
T^N\\
(\one_{s+\FF_q[T]}-\one_{\FF_q[T]})\otimes
\one_{1+(1/T)\FF_q[[1/T]]}
\end{array}\right)_+
\end{array}
\end{equation}
for $0\neq s=\frac{\alpha}{T}+\frac{\beta}{T-1}$ 
with $\alpha,\beta\in \FF_q$ and integers $N\geq 0$. 
Once again, we see that dependence on $N$ of the right side  is explained by a varying exponent $q^{N+1}$ on the left
and that the rational rigged virtual lattice on the right is proper.  Many more interpolation formulas analogous to (\ref{equation:SimpleExampleOneBis}) and (\ref{equation:SimpleExampleTwoBis}) are known
 which could in principle be expressed in terms of values of the partial Catalan symbol on proper rational rigged virtual lattices.
For example, there is a generalization of (\ref{equation:SimpleExampleTwoBis})
for all $0\neq s\in \FF_q(T)\cap(1/T)\FF_q[[1/T]]$ (too complicated to repeat
here---see \cite[\S8.4--5]{ThakurBook}) which plays a key role in the analysis  \cite{ABP} of transcendence properties of special values of the geometric $\Gamma$-function over $\FF_q[T]$. 
These examples suggest the possibility of a general interpolation formula relating partial Catalan symbols
to ``variable exponent expressions''.  We prove here in the local setting that such a thing does indeed exist.  Analogous global phenomena are expected. The author's paper \cite{AndersonStirling} explains the (still conjectural) global picture and makes the link to Stark's conjecture. 

\subsection{The asymptotic interpolation theorem} 
We formulate a weakened version of our main result precisely.
We then  indicate (without many details) the directions in which this weakened result 
is strengthened in order to achieve our main result. We comment briefly  on the global situation.

\subsubsection{Apparatus from local class field theory}
Let $\norm{\cdot}$ be the canonical absolute value
of $k$, i.~e., the modulus for Haar measure on $k$.
Let $k_\ab/k$ be the abelian closure of $k$ in a fixed
algebraic closure of $k$.
Let
$\rho:k^\times\rightarrow \Gal(k_\ab/k)$
be the reciprocity law homomorphism of local class field theory,
normalized as in \cite{TateBackground}. 
Under that normalization we have
$C^{\rho(a)}=C^{\norm{a}}$
for all $C$ in the algebraic closure of the prime field in $k^\ab$ and $a\in k^\times$. 
\subsubsection{Asymptotic interpolability}
Let $\OO$ be the ring of local integers of $k$.
Let $q$ be the cardinality of the residue field of $\OO$.
Let $\FF_q$ be the field of Teichm\"{u}ller representatives
of $\OO$. Let $X$ and $Y$ be independent variables.
Let $\FF_q[[X,Y]][X^{-1},Y^{-1}]$ be the ring obtained by 
adjoining inverses of $X$ and $Y$ to the two-variable
power series ring $\FF_q[[X,Y]]$.  
We call $\xi\in k_\ab$
a {\em basepoint} if there is some finite totally
ramified subextension $K/k$
of $k_\ab/k$ such that $\xi$ uniformizes the ring of local integers of $K$.
Given $F\in \FF_q[[X,Y]][X^{-1},Y^{-1}]$,
a basepoint $\xi$ and a proper rational rigged virtual lattice $\Phi$,
we say that $\Phi$ is {\em asymptotically yoked} to the pair $(F,\xi)$
if 
\begin{equation}\label{equation:AsymptoticInterpolation}
\left(\begin{array}{c}
a\\
\Phi
\end{array}\right)_+=F(\xi,(\rho(a)^{-1} \xi)^{\norm{a}})\;\;
\mbox{for all $a\in k^\times$ such that $\norm{a}\gg 1$.}
\end{equation}
In this situation we also say that $\Phi$ is {\em asymptotically interpolable}.  By a straightforward application of the Weierstrass Division Theorem one verifies that given $\Phi$ and $\xi$ there exists at most one $F$
such that $\Phi$ is asymptotically yoked to $(F,\xi)$.

We will prove the following result.
\begin{Subtheorem}\label{Subtheorem:Immediate}
Every proper rational rigged virtual lattice 
is a finite $\ZZ[1/p]$-linear combination
of asymptotically interpolable rigged virtual lattices.
\end{Subtheorem}
\noindent This result appears in the main body of the paper
as Corollary~\ref{Corollary:Immediate}
to our main result, Theorem~\ref{Theorem:Interpolability}. 

\subsubsection{Problems suggested by the asymptotic interpolability theorem}
Suppose that $\Phi$ is asymptotically yoked to $(F,\xi)$.
We may then consider the following problems:
\begin{enumerate}
\item  Find the prime factorization of $F$ in $\FF_q[[X,Y]][X^{-1},Y^{-1}]$.
\item Find the multiplicities (possibly negative) with which $X$ and $Y$ divide $F$.
\item For all $a\in k^\times$ (not just for $\norm{a}\gg 1$):
\begin{enumerate}
\item find $F(\xi,(\rho(a)^{-1}\xi)^{\norm{a}})$; in case of vanishing,
\item find the ``order of the zero'' of $F$ at $(\xi,(\rho(a)^{-1}\xi)^{\norm{a}})$, and
\item find the ``leading Taylor coefficient'' of $F$ at $(\xi,(\rho(a)^{-1}\xi)^{\norm{a}})$.
\end{enumerate}
\end{enumerate}
These problems are solved completely
by Theorem~\ref{Theorem:Interpolability} and the apparatus introduced to prove 
the theorem.

\subsubsection{Sketch of solution}
It turns out
that one can read the prime factorization of $F$ directly from the theta
symbol $\Theta(\cdot,\Phi)$, and that the $X$- and $Y$-multiplicities
of $F$ can be made explicit in terms of certain integrals involving $\Phi$
and the ramification index of $\xi$ over $k$.
Thus problems 1 and 2 are solved. Moreover, problems 1 and 3b turn out to be essentially the same. To solve problems 3a and 3c, we define a modified version $\GGG_0$
of the Fourier transform  which stabilizes the class of proper rational rigged virtual lattices,
and we define
$$\left(\begin{array}{c}
a\\
\Phi
\end{array}\right)=
\left(\begin{array}{c}
a\\
\Phi
\end{array}\right)_+\cdot\left(\begin{array}{c}
a^{-1}\\
\GGG_0[\Phi]\end{array}\right)_+^{\norm{a}}
$$
for all $a\in k^\times$ and proper rational rigged virtual lattices $\Phi$.
We call $\left(\begin{subarray}{c}
\cdot\\
\cdot
\end{subarray}\right)$ the {\em Catalan symbol}. 
It turns out that whenever we have ``one-sided'' asymptotic interpolation
as in (\ref{equation:AsymptoticInterpolation}), we automatically
also have ``two-sided'' asymptotic interpolation
$$\left(\begin{array}{c}
a\\
\Phi
\end{array}\right)=
F(\xi,(\rho(a)^{-1}\xi)^{\norm{a}})\;\;\;\mbox{for $a\in k^\times$
such that $\max(\norm{a},\norm{a}^{-1})\gg 1$.}$$
Further, it turns out that there holds a more delicate sort of  interpolation
which we call {\em strict interpolation}, valid for all $a\in k^\times$,
equating the Catalan symbol to a leading Taylor coefficient. 
Thus problems 3a and 3c are solved.

\subsubsection{Relationship to the global picture}
In the author's paper \cite{AndersonStirling}
global adelic versions of the theta and Catalan symbols are defined, a conjecture relating these objects to two-variable algebraic functions is proposed, the conjecture is verified in the genus zero case, and it is explained how the conjecture refines the Stark conjecture.
The results obtained in this paper, so we claim, provide detailed local confirmation of the author's conjecture. The claim requires justification which we do not provide here. We will take up the topic in future publications as part of our ongoing effort to prove the conjecture of \cite{AndersonStirling}. This paper is more or less self-contained. We consistently  take here  a local  point of view
complementary to the global point of view of \cite{AndersonStirling}.

\subsection{Organization of the paper}
In the next (second) section we 
fix notation and assumptions in force throughout the paper and we  develop
tools  to handle
products of the form (\ref{equation:InterestingProducts}).
We define  rigged virtual lattices, the theta symbol,
the rational Fourier transform $\GGG_0$,
and the Catalan symbol. We set up a calculus summarized by a short list of simple rules.
We emphasize {\em scaling rules} and {\em functional equations}. In the third
section of the paper we develop what amounts to a theory of expressions of the form
$F(\xi,(\rho(a)^{-1}\xi)^{\norm{a}})$ and in particular make rigorous sense of leading Taylor coefficients for such expressions. We develop a calculus of {\em shadow} theta and
Catalan symbols with scaling rules and functional equations parallel to those satisfied by the theta and Catalan symbols.
In the fourth section
of the paper we define asymptotic interpolability, strict interpolability and interpolability. We work out the formal properties of these notions in some detail. In the fifth section
we explicitly construct examples of strict interpolation.  In the final (sixth)
section of the paper we prove the main result by showing
that from the examples constructed in the fifth section 
all possible examples of strict interpolation can be built up by natural operations.

\section{Theta and Catalan symbols}
\label{section:RiggedStirling}
We develop local analogues of global notions introduced in \cite{AndersonStirling}.  We work in the setting of harmonic analysis on local fields. 
The main result of this section is Theorem~\ref{Theorem:RiggedStirling} 
(see equation (\ref{equation:SimplifiedRiggedStirling}) for a simplified version) which is analogous to the adelic Stirling formula \cite[Thm.\ 7.7]{AndersonStirling}.

\subsection{The set up}\label{subsection:BasicSetUp}
We specify the basic data for our constructions.
We introduce notation and terminology used throughout the paper.

\subsubsection{General notation}
Let $A^\times$ denote the group of units of a ring $A$ with unit.  
Let $\one_S$ denote the function equal to $1$ on a set $S$ and $0$ elsewhere. Let $\card S$ denote the cardinality of a set $S$.

\subsubsection{The local field $k$}
Fix a local field $k$ of characteristic $p>0$ with maximal compact subring $\OO$.
Let $\FF_q=\{x\in k\mid x^q=x\}$, where $q$ is  the (finite) cardinality of the residue field of $\OO$.
Then $\FF_q$ is both a subfield of $k$ and a set of representatives for the residue field of $\OO$.
Moreover, for every uniformizer $\pi\in \OO$,
we have $\OO=\FF_q[[\pi]]$ and $k=\FF_q[[\pi]][\pi^{-1}]$.

\subsubsection{Chevalley differentials}
Let $\Omega$ be the space of locally constant $\FF_q$-linear functionals $k\rightarrow\FF_q$. Elements of $\Omega$ are called {\em Chevalley differentials}.
For every $\delta\in \Omega$,  let $\Res \delta=\delta(1)\in \FF_q$.
For every $x\in k$ and $\delta\in \Omega$, let $x\delta\in \Omega$ be 
defined  by the formula $(x\delta)(y)=\delta(xy)$ for all $y\in k$; in this way $\Omega$  becomes 
a vector space one-dimensional  over $k$.  
It is well-known that there exists a unique
$\FF_q$-linear derivation $d:k\rightarrow\Omega$ such that 
$$\Res(x\,d\pi)=a_{-1}\;\;\;\left(x=\sum_{i=-\infty}^\infty
a_i\pi^i\;\;\left(a_i\in \FF_q,\;a_i=0\;\mbox{for $i\ll 0$}\right)\right)$$
for all uniformizers $\pi\in\OO$ and $x\in k$.
Given $\xi,\eta\in \Omega$ with $\eta\neq 0$
and $x\in k$ such that $\xi=x \eta$,
we write $x=\xi/\eta$.

\subsubsection{The Chevalley differential $\omega$}
Fix $0\neq \omega\in \Omega$
and let $\DDD$ be the fractional $\OO$-ideal defined by the rule
$$\DDD^{-1}=\{\xi\in k\mid \Res(\xi x \omega)=0\;\mbox{for all $x\in \OO$}\}.$$ 
Many constructions below depend on the choice of $\omega$. But for the most part we suppress reference to $\omega$ in the notation.

\subsubsection{Valuations and measures}
Let $\mu$ be Haar measure on $k$.
 Normalize $\mu$ by requiring that
  $$\mu\OO\cdot \mu\DDD^{-1}=1.$$
 Let $\mu^\times$ be Haar measure on $k^\times$.
 Normalize $\mu^\times$ by requiring that $$\mu^\times \OO^\times=1.$$
 For each $a\in k$, put $$\norm{a}=\mu(a\OO)/\mu\OO,\;\;\;\;\ord a=-\log \norm{a}/\log q.$$
The function $\norm{\cdot}$
 is an absolute value of $k$ with respect to which $k$ is complete. The function $\ord$ is an additive valuation of $k$. Note that for every uniformizer $\pi\in \OO$,  we have
 $$\norm{\pi}=q^{-1},\;\;\;\ord \pi=1.$$ Note also that 
 $$\int f(a^{-1}x)d\mu(x)=\norm{a}\int f(x)d\mu(x)$$
 for all $a\in k^\times$ and $\mu$-integrable complex-valued functions $f$.
  For every $\delta\in \Omega$, put $$\norm{\delta}=\norm{\delta/d\pi},\;\;\; \ord \delta=-\log\norm{\delta}/\log q,$$
  where $\pi\in \OO$ is any uniformizer.  Note that $$\mu\OO=\norm{\omega}^{1/2}.$$

\subsubsection{The field $\FF$ and character $\lambda$}
Fix a subfield $\FF\subset \FF_q$.  Fix a nonconstant homomorphism
$\lambda$ from the additive group of $\FF$ to the multiplicative group  of nonzero complex numbers. Many of our constructions depend on $\FF$ or $\lambda$ or both,
but for the most part we suppress reference to $\FF$ and $\lambda$ in the notation.

\subsubsection{The character $\ee$}
We define
\begin{equation}\label{equation:eeomega}
\ee(x)=\lambda(\trace_{\FF_q/\FF}\Res(x\omega))
\end{equation}
for all $x\in k$.  This is the character we will use to define the Fourier transform on $k$.

\subsubsection{Bruhat-Schwartz space}
Given a locally compact totally disconnected Hausdorff space $X$,
e.g., $X=k$ or $X=k\times k$,
let $\SSS(X)$ denote the space of complex-valued compactly supported locally constant
 functions on $X$.

\subsubsection{Further notation}\label{subsubsection:Further}
Let $\bar{k}$
be a fixed algebraic closure of $k$.  Let $k_\perf$ be the closure of $k$ in $\bar{k}$ under extraction of $p^{th}$ roots. Let $k_\ab$ be the abelian closure of $k$
in $\bar{k}$.
For each positive integer $n$, let $\FF_{q^n}$ be the unique subfield of $\bar{k}$ of cardinality $q^n$. 
Put $\FF_{q^\infty}=\bigcup_{n=1}^\infty \FF_{q^n}\subset k_\ab$.
Let the absolute value $\norm{\cdot}$ and additive valuation $\ord$ on $k$
be canonically extended to $\bar{k}$,
and again denoted by $\norm{\cdot}$ and  $\ord$, respectively.
\pagebreak

\subsection{Apparatus from harmonic analysis}
\subsubsection{Fourier transforms}
Given $\varphi\in \SSS(k)$, put
$$\FFF[\varphi](\xi)=\hat{\varphi}(\xi)=
\int \varphi(x)\ee(-x\xi)d\mu(x),$$
thus defining the {\em Fourier transform}
$$\FFF[\varphi]=\hat{\varphi}\in \SSS(k).$$
Given $\varphi\in \SSS(k)$ and $a\in k^\times$, put 
$$\varphi^{(a)}(x)=\varphi(a^{-1}x).$$
We have a squaring rule
$$\FFF^2[\varphi]=\varphi^{(-1)}$$
holding for all $\varphi\in \SSS(k)$. (Our normalization of additive Haar measure $\mu$ was chosen to put the Fourier inversion formula into this simple form.)
Further, we have by an evident transformation of integrals a scaling rule
$$
\FFF[\varphi^{(a)}]=\norm{a}\FFF[\varphi]^{(a^{-1})}
$$
holding for all $\varphi\in \SSS(k)$ and $a\in k^\times$.
\subsubsection{The Poisson summation formula}
Let $L\subset k$ be a cocompact discrete subgroup.
The {\em Poisson summation formula} states that for every $\varphi\in \SSS(k)$ we have
\begin{equation}\label{equation:BasicPoisson}
\sum_{x\in L}\varphi(x)=
\frac{1}{\mu(k/L)}\sum_{\xi\in L^\perp}\hat{\varphi}(\xi)
\end{equation}
where $L^\perp$ is the cocompact discrete subgroup of $k$ defined by the rule
$$\begin{array}{rcl}
L^\perp&=&\{\xi\in k\mid \ee(x\xi)=1\;\mbox{for all}\;x\in L\}\\\\
&=&
\{\xi\in k\mid \Res(x\xi \omega)=0\;\mbox{for all}\;x\in L\},
\end{array}
$$
and $\mu(k/L)$ is the $\mu$-measure of any fundamental domain for $L$.

\subsubsection{Two-variable Fourier transforms}
Given
complex-valued functions $f_1$ and $f_2$ on $k$, set $$(f_1\otimes f_2)(x_1,x_2)=f_1(x_1)f_2(x_2),$$ thus defining a complex-valued function $f_1\otimes f_2$ on $k\times k$. As is well-known, the $\otimes$-operation
identifies $\SSS(k\times k)$ with the tensor square of $\SSS(k)$ over $\CC$.
Given $\varphi\in \SSS(k\times k)$, put
$$\GGG[\varphi](\xi,\eta)=
\int\int \varphi(x,y)\ee(x\xi-y\eta)d\mu(x)d\mu(y),$$
thus defining the {\em two-variable Fourier transform}
$$\GGG[\varphi]\in \SSS(k\times k).$$
Note that
\begin{equation}\label{equation:GFactorization}
\GGG[\varphi_1\otimes \varphi_2]=
\FFF^{-1}[\varphi_1]\otimes \FFF[\varphi_2]=
\hat{\varphi}_1^{(-1)}\otimes \hat{\varphi}_2
\end{equation}
for all $\varphi_1,\varphi_2\in \SSS(k)$. (The asymmetry of the definition of $\GGG$ is dictated by our goal of simplifying formula (\ref{equation:LowerStarPoissonSpecial}) below as much as possible.)
 
\subsubsection{Generalized functions of two variables}
Let $\SSS'(k\times k)$ be the
space of {\em generalized functions} on $k\times k$, i.~e., 
the $\CC$-linear dual of $\SSS(k\times k)$, and let 
$$\langle \cdot,\cdot\rangle:\SSS'(k\times k)
\times \SSS(k\times k)\rightarrow \CC$$
be the canonical pairing.  The theory of generalized functions in our context follows the pattern set classically by the theory of distributions.
We identify $\SSS(k\times k)$ with a subspace of $\SSS'(k\times k)$ by the rule
$$\langle \varphi,\psi\rangle=(\mu\OO)^{-2}
\int\int \varphi(x,y)\psi(x,y)d\mu(x)d\mu(y)$$
for all $\varphi,\psi\in \SSS(k)$. We have inserted the factor $(\mu\OO)^{-2}=\norm{\omega}^{-1}$ to make the identification independent of the choice of $\omega$. We 
extend the Fourier transform $\GGG$ to 
$\SSS'(k\times k)$ by the rule
$$\langle \GGG[\Gamma],\varphi\rangle=
\langle \Gamma,\GGG[\varphi]\rangle$$
for all $\Gamma\in \SSS'(k\times k)$ and $\varphi\in \SSS(k\times k)$.
Given $\varphi\in \SSS(k\times k)$ and $b,c\in k^\times$, put $$\varphi^{(b,c)}(x,y)=\varphi(b^{-1}x,c^{-1}y).$$
We extend the operation $\varphi\mapsto \varphi^{(b,c)}$ to
$\SSS'(k\times k)$ by the rule 
\begin{equation}\label{equation:TwoVariableScalingDef}
\langle \Gamma^{(b,c)},\varphi^{(b,c)}\rangle=
\norm{bc}\langle \Gamma,\varphi\rangle
\end{equation}
for all $\Gamma\in \SSS'(k\times k)$ and $\varphi\in \SSS(k\times k)$. 
We have  squaring and scaling rules
\begin{equation}\label{equation:TwoVariableFourierSquaring}
\GGG^2[\Gamma]=\Gamma^{(-1,-1)},
\end{equation}
\begin{equation}\label{equation:TwoVariableFourierScaling}
\GGG[\Gamma^{(b,c)}]=\norm{bc}\GGG[\Gamma]^{(b^{-1},c^{-1})},
\end{equation} 
respectively, holding
for all $\Gamma\in \SSS'(k\times k)$ and $b,c\in k^\times$, 
in evident analogy with the squaring and scaling rules obeyed by $\FFF$.

\subsection{Rigged virtual lattices}
\subsubsection{Definition}
We say that $\Phi\in \SSS'(k\times k)$ is a {\em primitive rigged virtual lattice}
if there exist
\begin{itemize}
\item a cocompact discrete subgroup $L\subset k$, 
\item numbers $\ell,w\in k$
and 
\item a number $r\in k^\times$
\end{itemize}
 such that
\begin{equation}\label{equation:TypicalRigged}
\langle \Phi,\varphi\rangle=(\mu \OO)^{-1}
\sum_{x\in \ell+L}\int_{w+r\OO} \varphi(x,y)d\mu(y)
\end{equation}
for all $\varphi\in \SSS(k\times k)$.
We have inserted the factor
$(\mu\OO)^{-1}=\norm{\omega}^{-1/2}$ so that the notion of
primitive rigged virtual lattice is defined independently of the choice of the Chevalley differential $\omega$. 
We define a {\em rigged virtual lattice} to be a finite linear combination 
of primitive rigged virtual lattices with complex coefficients,
and we denote the space of such by $\rigged(k)$.
It is convenient to associate to each 
$\Phi\in \rigged(k)$ a complex-valued function $\Phi_*$ on $k\times k$ in the usual (rather than generalized) sense by the rule
\begin{equation}\label{equation:LowerStarRule}
\Phi_*(x,y)=\lim_{\norm{a}\rightarrow 0}
\langle \Phi,\one_{x+a\OO}\otimes \one_{y+a\OO}\rangle/\norm{a}.
\end{equation}
Note that the ``lower star rule'' (\ref{equation:LowerStarRule}) for passing from $ \rigged(k)$ to functions on $k\times k$
is independent of the choice of $\omega$.
To see that $\Phi_*$ is well-defined,
note that if $\Phi\in  \rigged(k)$ takes the form (\ref{equation:TypicalRigged}), then 
 \begin{equation}\label{equation:TypicalRiggedBis}
\Phi_*=\one_{\ell+L}\otimes \one_{w+r\OO},
\end{equation}
as one verifies by a straightforward calculation. 
Via (\ref{equation:TypicalRigged})
and (\ref{equation:TypicalRiggedBis}),  we have
\begin{equation}\label{equation:PhiStarRep}
\langle \Phi,\varphi\rangle=
(\mu\OO)^{-1}\sum_{x\in k}\int \Phi_*(x,y)\varphi(x,y)d\mu(y)
\end{equation}
for all $\Phi\in  \rigged(k)$ and  $\varphi\in \SSS(k)$.
Thus $\Phi_*$  uniquely determines $\Phi$. By  (\ref{equation:TwoVariableScalingDef}) and (\ref{equation:PhiStarRep}), we have the scaling rule
\begin{equation}\label{equation:RiggedScaling}
\Phi^{(b,c)}\in  \rigged(k),\;\;\;(\Phi^{(b,c)})_*(x,y)=\norm{b}\Phi_*(b^{-1}x,c^{-1}y)
\end{equation}
holding for all $\Phi\in  \rigged(k)$ and $b,c\in k^\times$. The rule (\ref{equation:RiggedScaling}) is a crucially important bookkeeping detail in our theory.
The space $ \rigged(k)$ is stable
under the action of the two-variable Fourier transform $\GGG$,
as the next lemma and its proof show.

\begin{Lemma}\label{Lemma:RVFT}
Let $L\subset k$ be a cocompact discrete subgroup
and fix $\ell\in k$. Fix $\rho\in \SSS(k)$.
If $\Phi\in  \rigged(k)$ takes the form 
$$\Phi_*=\one_{\ell+L}\otimes \rho,$$
 then 
 \begin{equation}\label{equation:LowerStarFourierSpecial}
\GGG[\Phi]\in  \rigged(k),\;\;\;\GGG[\Phi]_*(\xi,\eta)=
\frac{\ee(\ell \xi)}{\mu(k/L)}\one_{L^\perp}(\xi)\hat{\rho}(\eta),
\end{equation}
 and moreover
\begin{equation}\label{equation:LowerStarPoissonSpecial}
\sum_{x\in k}\Phi_*(x,x)=\sum_{\xi\in k}
\GGG[\Phi]_*(\xi,\xi).
\end{equation}
\end{Lemma}
\noindent The sums on both sides of (\ref{equation:LowerStarPoissonSpecial})
are well-defined since both have only finitely many nonzero terms.
\proof We claim that
\begin{equation}\label{equation:PoissonClaim}
\langle \Phi,\GGG[\psi]\rangle=
(\mu\OO)^{-1}\sum_{x\in L^\perp}
 \frac{\ee(\ell x)}{\mu(k/L)}\int\psi(x,y)\hat{\rho}(y)d\mu(y)
\end{equation}
for all $\psi\in \SSS(k\times k)$. 
The claim  granted, (\ref{equation:LowerStarFourierSpecial}) follows
by (\ref{equation:PhiStarRep}) and the definitions.
 To prove (\ref{equation:PoissonClaim}), we may assume without loss of generality that $\psi=\varphi_1\otimes \varphi_2$, where $\varphi_1,\varphi_2\in \SSS(k)$. 
Plugging now into (\ref{equation:TypicalRigged}), using (\ref{equation:GFactorization})
to calculate $\GGG[\varphi_1\otimes\varphi_2]$,
the left side of (\ref{equation:PoissonClaim})
takes the form
$$\left(\sum_{x\in \ell+L}\hat{\varphi}_1^{(-1)}(x)\right)\left( (\mu\OO)^{-1}\int\int \varphi_2(x)\rho(y)\ee(-xy)d\mu(x) d\mu(y)\right),$$
whereas the right side takes the form
$$\left(\sum_{\xi\in L^\perp}
\frac{\ee(\ell\xi)\varphi_1(\xi)}{\mu(k/L)}\right)\left( (\mu\OO)^{-1}\int\int \varphi_2(x)\rho(y)\ee(-xy)d\mu(x) d\mu(y)\right).$$
Now from the Poisson summation formula (\ref{equation:BasicPoisson}),
by substituting  $x\mapsto \varphi(x+\ell)$ for $\varphi$,
we deduce the slightly more general formula
\begin{equation}\label{equation:PrePoisson}
\sum_{x\in \ell+L}\varphi(x)=
\sum_{\xi\in L^\perp}\frac{\ee(\ell \xi)\hat{\varphi}(\xi)}{\mu(k/L)}
\end{equation}
holding for all $\varphi\in \SSS(k)$.
Equality of left and right sides of (\ref{equation:PoissonClaim})  reduces to
the special case $\varphi=\FFF^{-1}[\varphi_1]$
of (\ref{equation:PrePoisson}). The claim (\ref{equation:PoissonClaim}) is proved.
Similarly and finally, formula (\ref{equation:LowerStarPoissonSpecial}) reduces to the special case $\varphi=\rho$ of (\ref{equation:PrePoisson}).  The lemma is proved.
 \qed

\subsection{The theta symbol}
\subsubsection{Definition}
For $a\in k^\times$ and $\Phi\in  \rigged(k)$,
 put 
$$\Theta(a,\Phi)=\sum_{x\in k}(\Phi^{(1,a)})_*(x,x)=
\sum_{x\in k}\Phi_*(x,a^{-1}x),$$
where the second equality is justified by scaling rule (\ref{equation:RiggedScaling}).
Only finitely many nonzero terms appear in the sum on the extreme right, as one can verify by 
passing without loss of generality to the special case (\ref{equation:TypicalRigged})
and applying formula (\ref{equation:TypicalRiggedBis}). Indeed, one finds that
\begin{equation}\label{equation:TypicalRiggedTer}
\begin{array}{rl}
&\Phi_*=\one_{\ell+L}\otimes
\one_{w+r\OO}\\\\
\Rightarrow&\displaystyle\Theta(a,\Phi)=\sum_{x\in k}\one_{(\ell+L)\cap a(w+r\OO)}(x)=
\card ((\ell+L)\cap a(w+r\OO)).
\end{array}
\end{equation}
 Therefore $\Theta(a,\Phi)$ is well-defined.
We call $\Theta(\cdot,\cdot)$ the {\em theta symbol}.
The object $\Theta(\cdot,\cdot)$ defined here is a local version of the theta symbol defined  in \cite{AndersonStirling}.

\subsubsection{Basic formal properties}
Fix $\Phi\in  \rigged(k)$ and $a\in k^\times$.
Via  scaling rule (\ref{equation:RiggedScaling}) and the definition of the theta symbol, we have    the scaling rule
\begin{equation}\label{equation:ThetaScaling}
\Theta(a,\Phi^{(b,c)})=\norm{b}\Theta(ac/b,\Phi)
\end{equation}
holding for all $\Phi\in  \rigged(k)$ and $a,b,c\in k^\times$.
We claim the following:
\begin{equation}\label{equation:ThetaFE}
\Theta(a,\Phi)=
 \norm{a}\Theta(a^{-1},\GGG[\Phi]).
 \end{equation}
\begin{equation}\label{equation:ThetaAsymptotics}
\Theta(a,\Phi)=\left\{\begin{array}{cl}
\Phi_*(0,0)&\mbox{if $\norm{a}$ is sufficiently small,}\\
\GGG[\Phi]_*(0,0)\norm{a}&\mbox{if $\norm{a}$ is sufficiently large.}\end{array}\right.
\end{equation}
\begin{equation}\label{equation:ThetaLocalConstancy}
\mbox{The function $\left(a\mapsto \Theta(a,\Phi)\right):k^\times \rightarrow \ZZ[1/p]$ is locally constant.}
\end{equation}
To prove all three claims we may assume without loss of generality that
$\Phi$ is of the form (\ref{equation:TypicalRigged})
to which Lemma~\ref{Lemma:RVFT} applies.
Functional equation (\ref{equation:ThetaFE}) follows from formula (\ref{equation:LowerStarPoissonSpecial}),  scaling rule (\ref{equation:ThetaScaling}) and the definitions. Further, (\ref{equation:ThetaAsymptotics})
for small $\norm{a}$ and (\ref{equation:ThetaLocalConstancy}) are easy to check
using (\ref{equation:TypicalRiggedTer}).
   Finally, (\ref{equation:ThetaAsymptotics})
for large $\norm{a}$ follows by (\ref{equation:ThetaFE}) from (\ref{equation:ThetaAsymptotics})
for small $\norm{a}$. The claims are proved. We call (\ref{equation:ThetaFE}) the {\em functional equation} satisfied by the theta symbol.

\subsubsection{Properness and effectiveness of rigged virtual lattices}
We say that \linebreak $\Phi\in \rigged(k)$ is {\em proper} if $\Theta(\cdot,\Phi)$
is compactly supported, i.~e., if $\Theta(a,\Phi)=0$ 
for $\max(\norm{a},\norm{a}^{-1})\gg 0$. Equivalently, 
$\Phi$ is proper if $\Phi_*(0,0)=0$ and $\GGG[\Phi]_*(0,0)=0$.
The space of  proper rigged virtual lattices is stable under the action
of the two-variable Fourier transform
$\GGG$, stable under formation of finite $\CC$-linear combinations,
and for all $b,c\in k^\times$ stable  under
the operation $\Phi\mapsto \Phi^{(b,c)}$.
We say that $\Phi\in \rigged(k)$ is {\em effective} if $\Theta(a,\Phi)$
is a nonnegative real number for all $a\in k^\times$.
The space of effective rigged virtual lattices is stable under the two-variable Fourier transform
$\GGG$, stable under formation of finite nonnegative real linear combinations,
and for all $b,c\in k^\times$ stable under
the operation $\Phi\mapsto \Phi^{(b,c)}$.

\begin{Theorem}[The Stirling formula for rigged virtual lattices]
\label{Theorem:RiggedStirling}
For all $a\in k^\times$ and  $\Phi\in  \rigged(k)$,
the formula
\begin{equation}\label{equation:RiggedStirling}
\begin{array}{cl}
&\displaystyle \Theta(a,\Phi)\ord \omega+
\sum_{x\in k^\times} \Phi_*(x,a^{-1}x)\ord x+
\sum_{\xi\in k^\times}\norm{a}\GGG[\Phi]_*(\xi,a\xi)\ord \xi\\&\\
=&{\displaystyle
\int\frac{\Theta(at,\Phi)-\left\{\begin{array}{cl}
\norm{at}\GGG[\Phi]_*(0,0)&\mbox{if $\norm{t}>1$}\\
\Theta(a,\Phi)&\mbox{if $\norm{t}=1$}\\
\Phi_*(0,0)&\mbox{if $\norm{t}<1$}
\end{array}\right.}{\norm{1-t}}d\mu^\times(t)}\\&\\
&{\displaystyle -\Phi_*(0,0)\ord
a+\int\left(\Phi_*(0,t)-\Phi_*(0,0)\one_{\OO}(t)\right)d\mu^\times(t)}\\&\\
&\displaystyle
+\norm{a}\bigg(\GGG[\Phi]_*(0,0)\ord a+\int(
\GGG[\Phi]_*(0,t)-\GGG[\Phi]_*(0,0)
\one_{\OO}(t))d\mu^\times(t)\bigg)
\end{array}
\end{equation}
holds.
\end{Theorem}
\noindent The sums on the left side of (\ref{equation:RiggedStirling}) are absolutely convergent since there appear in each only finitely many nonzero terms. The integral involving $\Theta(\cdot,\cdot)$ on the right side of (\ref{equation:RiggedStirling}) is absolutely convergent in view of the formal properties of the theta symbol that we have just proved. The remaining integrals on the right side of (\ref{equation:RiggedStirling}) are absolutely convergent since their integrands are compactly supported in $k^\times$. We remark that (\ref{equation:RiggedStirling}) simplifies considerably if $\Phi$ is proper. The theorem is closely analogous
to \cite[Theorem 7.7]{AndersonStirling}
and a generalization of \cite[Corollary 6.11]{AndersonStirling}.
\proof 
Replacement of the pair $(a,\Phi)$
by the pair $(1,\Phi^{(1,a)})$ leaves both sides of 
(\ref{equation:RiggedStirling}) unchanged, as one verifies with the help of the scaling rules (\ref{equation:TwoVariableFourierScaling}),
(\ref{equation:RiggedScaling}), and (\ref{equation:ThetaScaling}), along with some evident transformations of integrals.
We may therefore assume without loss of generality that $a=1$ for the rest of the proof.
Further, we may assume without loss of generality that $\Phi$ is as in Lemma~\ref{Lemma:RVFT}, namely
$$\Phi=\one_{\ell+L}\otimes \rho\;\mbox{and}\;\GGG[\Phi]_*(\xi,\eta)=
\frac{\ee(\ell \xi)}{\mu(k/L)}\one_{L^\perp}(\xi)\hat{\rho}(\eta),$$
for some $\rho\in \SSS(k)$, in which case functional equation (\ref{equation:ThetaFE}) written out explicitly (with $t\in k^\times$ in place of $a\in k^\times$) takes  the form 
\begin{equation}\label{equation:ThetaFENow}
\sum_{x\in \ell+L}\rho(t^{-1}x)
=\norm{t}\sum_{\xi\in L^\perp}
\frac{\ee(\ell \xi)\hat{\rho}(t\xi)}{\mu(k/L)}.
\end{equation}
The sum on the left equals
$\Theta(t,\Phi)$ and the sum on right equals $\norm{t}\Theta(t^{-1},\GGG[\Phi])$. Now put $H=\one_{\OO}-\frac{1}{2}\one_{\OO^\times}$.
For every $\varphi\in \SSS(k)$ and $x\in k^\times$, put
$$\begin{array}{rcl}
\MMM^{(1)}[\varphi]&=&\displaystyle\int(\varphi(t)-\one_\OO(t)\varphi(0))d\mu^\times (t),\\\\
\LLL^+[\varphi](x)&=&\displaystyle\int\frac{H(t)\varphi(xt^{-1})-\frac{1}{2}\one_{\OO^\times}(t)\varphi(x)}{\norm{1-t}}d\mu^\times(t)\\\\
&=&\displaystyle\int\frac{H(t^{-1})\norm{t}{\varphi}(x  t)-\frac{1}{2}\one_{\OO^\times}(t){\varphi}(x)}{\norm{1-t}}d\mu^\times(t).
\end{array}
$$
By the local Stirling formula \cite[Thm.\ 6.6]{AndersonStirling},
there exists unique $\psi\in \SSS(k)$ such that
$$
\begin{array}{rcl}
\psi(x)&=&
\left\{\begin{array}{rl}
-\rho(x)(\frac{1}{2}\ord\omega+\ord x)+
\LLL^+[\rho](x)&\mbox{if $x\neq 0$,}\\
-\frac{1}{2}\rho(0)\ord\omega+\MMM^{(1)}[\rho]&\mbox{if $x=0$,}
\end{array}\right.\\\\
-\hat{\psi}(\xi)&=&
\left\{\begin{array}{rl}
-\hat{\rho}(\xi)(\frac{1}{2}\ord\omega+\ord \xi)+
\LLL^+[\hat{\rho}](\xi)&\mbox{if $\xi\neq 0$,}\\
-\frac{1}{2}\hat{\rho}(0)\ord\omega+\MMM^{(1)}[\hat{\rho}]&\mbox{if $\xi=0$.}
\end{array}\right.
\end{array}
$$
Using the fact that the left side of (\ref{equation:ThetaFENow})
represents $\Theta(t,\Phi)$, and using \cite[Lemma 6.9]{AndersonStirling} to justify the exchange of summation and integration processes, we have
$$
\begin{array}{cl}
&\displaystyle\sum_{0\neq x\in \ell+L}
\LLL^+[\rho](x)\\\\
=&\displaystyle
\int
\frac{H(t)(\Theta(t,\Phi)-A)-\frac{1}{2}
\one_{\OO^\times}(t)
(\Theta(1,\Phi)-A)}{\norm{1-t}}d\mu^\times(t),
\end{array}
$$
where 
$$A=\one_L(\ell)\rho(0)=\Phi_*(0,0).$$
Reasoning similarly, but this time using the fact that the right side of (\ref{equation:ThetaFENow}) represents $\norm{t}\Theta(t^{-1},\GGG[\Phi])$, we have
$$\sum_{0\neq \xi\in L^\perp}\frac{\ee(\ell\xi)\LLL^+[\hat{\rho}](\xi)}{\mu(k/L)}$$
$$=
\int
\frac{H(t^{-1})\norm{t}(\Theta(t^{-1},\GGG[\Phi])-B)-\frac{1}{2}
\one_{\OO^\times}(t)
(\Theta(1,\GGG[\Phi])-B)}{\norm{1-t}}d\mu^\times(t),
$$
where 
$$B=\frac{\hat{\rho}(0)}{\mu(k/L)}=\GGG[\Phi]_*(0,0).$$
It is now only a matter of bookkeeping to verify that 
relation (\ref{equation:RiggedStirling}) in the case $a=1$ coincides with relation (\ref{equation:ThetaFENow}) in the case $\psi=\rho$ and $t=1$; these calculations are quite similar to those undertaken to prove \cite[Cor.\ 6.11]{AndersonStirling} and so we omit further details. Thus the theorem is proved. 
\qed

\pagebreak

\subsection{The   rational Fourier transform}
\subsubsection{Trivial identities}
Put
\begin{equation}\label{equation:LambdaNought}
\lambda_0:=\left(x\mapsto \left\{\begin{array}{rl}
1&\mbox{if $x=0$,}\\
-1&\mbox{if $x=1$,}\\
0&\mbox{otherwise}
\end{array}\right.\right):\FF_q\rightarrow\{-1,0,1\}.
\end{equation}
 Let $\Sigma'_C$ abbreviate $\sum_{C\in \FF^\times}$.
 We have
 \begin{equation}\label{equation:lambda1}
\card \FF\cdot \lambda_0(x)=\Sigma'_C(1-\lambda(C))\lambda(-Cx),
\end{equation}
\begin{equation}\label{equation:lambda3}
\Sigma'_C\lambda(Cx)=\Sigma'_C\lambda_0(Cx)
\end{equation}
for all $x\in \FF$.

\subsubsection{Definition of $\GGG_0$}\label{subsubsection:RationalFourierDefinition}  Put 
$$\ee_0:=(x\mapsto \lambda_0(\trace_{\FF_q/\FF}\Res x\omega)):k\rightarrow\{-1,0,1\}.$$
Equivalently, via (\ref{equation:eeomega}) and (\ref{equation:lambda1}), we have
\begin{equation}\label{equation:eeToee0Conversion}
\ee_0(x)=(\card \FF)^{-1}\Sigma'_C(1-\lambda(C))\ee(-Cx).
\end{equation}
Given $\varphi\in \SSS(k\times k)$, put
$$\GGG_0[\varphi](\xi,\eta):=
\int\int \varphi(x,y)\ee_0(-x\xi+y\eta)d\mu(x)d\mu(y),$$
thus defining a complex-valued function on $k\times k$ which by (\ref{equation:eeToee0Conversion}) 
belongs to $\SSS(k\times k)$. 
We call $\GGG_0[\varphi]$ the {\em  rational Fourier transform} of $\varphi$. The operator $\GGG_0$ defined here is a local version of the global rational Fourier transform defined in \cite{AndersonStirling}. 
By an evident transformation of integrals we have a scaling rule
\begin{equation}\label{equation:PreRationalFourierScaling}
\GGG_0[\varphi^{(b,c)}]=\norm{bc}\GGG_0[\varphi]^{(b^{-1},c^{-1})}
\end{equation}
holding for all $b,c\in k^\times$. By (\ref{equation:lambda1}), (\ref{equation:lambda3}) 
and the definitions, we have
\begin{equation}\label{equation:PreRelationOfFTs}
\GGG_0[\varphi]=(\card \FF)^{-1}\Sigma'_C(1-\lambda(C))\GGG[\varphi]^{(C^{-1},C^{-1})},
\end{equation}
\begin{equation}\label{equation:PreRelationOfFTsTer}
\Sigma_C'\GGG[\varphi]^{(C,C)}=\Sigma_C'\GGG_0[\varphi]^{(C,C)}.
\end{equation}
We extend $\GGG_0$ to an operator on $\SSS'(k\times k)$ by the rule
$$\langle \GGG_0[\Gamma],\varphi\rangle=\langle\Gamma,\GGG_0[\varphi]\rangle$$
for all $\Gamma\in \SSS'(k\times k)$ and $\varphi\in \SSS(k\times k)$. The operator $\GGG_0$ so extended satisfies the evident generalizations of (\ref{equation:PreRationalFourierScaling}), (\ref{equation:PreRelationOfFTs}) and
(\ref{equation:PreRelationOfFTsTer}).

\begin{Lemma}\label{Lemma:RationalRVFT}
Fix $\ell,w\in k$, $r\in k^\times$ and a cocompact discrete subgroup $L\subset k$.
Given $\Phi\in \rigged(k)$ such that
$$\Phi_*=\one_{\ell+L}\otimes \one_{w+r\OO},$$
we have
\begin{equation}\label{equation:LowerStarFourierSpecialNought}
\begin{array}{l}
\GGG_0[\Phi]\in  \rigged(k),\\\\
\displaystyle\GGG_0[\Phi]_*(x,y)=
\frac{\ee_0(-\ell x+wy)\mu(r\OO)}{\mu(k/L)}
\one_{L^\perp}(x)\one_{r^{-1}\DDD^{-1}}(y).
\end{array}
\end{equation}
\end{Lemma}
\proof By Lemma~\ref{Lemma:RVFT} we have
$$\begin{array}{l}
\GGG[\Phi]\in  \rigged(k),\\\\
\displaystyle\GGG[\Phi]_*(x,y)=
\frac{\ee(\ell x-wy)\mu(r\OO)}{\mu(k/L)}
\one_{L^\perp}(x)\one_{r^{-1}\DDD^{-1}}(y).
\end{array}$$
By identity (\ref{equation:PreRelationOfFTs}) in its generalized form we have
$$\GGG_0[\Phi]=(\card \FF)^{-1}\Sigma'_C(1-\lambda(C))\GGG[\Phi]^{(C^{-1},C^{-1})}\in \rigged(k).$$
We get the claimed formula for $\GGG_0[\Phi]_*$
now by scaling rule (\ref{equation:RiggedScaling}) 
and formula (\ref{equation:eeToee0Conversion}). \qed
 \subsubsection{Basic formal properties of $\GGG_0$}
By the scaling rule (\ref{equation:RiggedScaling}) for the lower star operation and the scaling rule (\ref{equation:PreRationalFourierScaling}) in its generalized form, we have a scaling rule
\begin{equation}\label{equation:RationalFourierScalingRule}
\GGG_0[\Phi^{(b,c)}]_*(x,y)=\norm{c}\GGG_0[\Phi]_*(bx,cy)
\end{equation}
holding for all $\Phi\in \rigged(k)$ and $b,c\in k^\times$.
By the squaring rule (\ref{equation:TwoVariableFourierSquaring}) and scaling rule (\ref{equation:TwoVariableFourierScaling}) for $\GGG$, scaling rule (\ref{equation:RiggedScaling}) for the lower star operation,
and relation (\ref{equation:PreRelationOfFTs}) in its generalized form, we have a squaring rule
\begin{equation}\label{equation:RTSquare}
\card \FF\cdot \GGG_0^2[\Phi]_*(x,y)=\Phi_*(x,y)+\Sigma'_C\Phi_*(Cx,Cy)
\end{equation}
for all $\Phi\in \rigged(k)$, after a brief calculation which we omit.
By the scaling rule (\ref{equation:RiggedScaling}) for the lower star operation and relation (\ref{equation:PreRelationOfFTsTer}) in its generalized form, we have
\begin{equation}\label{equation:RelationOfFTsTer}
\Sigma_C'\GGG[\Phi]_*(Cx,Cy)=\Sigma_C'\GGG_0[\Phi]_*(Cx,Cy).
\end{equation}
By scaling rule (\ref{equation:ThetaScaling})
and relation (\ref{equation:RelationOfFTsTer}),
the functional equation (\ref{equation:ThetaFE}) satisfied by the theta symbol can be rewritten in the form
\begin{equation}\label{equation:NoughtThetaFE}
\Theta(a,\Phi)=\norm{a}\Theta(a^{-1},\GGG_0[\Phi])
\end{equation}
holding for all $a\in k^\times$ and $\Phi\in  \rigged(k)$. In other words, we can simply replace $\GGG$ by $\GGG_0$ in (\ref{equation:ThetaFE}) without invalidating the latter,
and this noted, it follows that we can test $\Phi\in \rigged(k)$ for properness
by checking that $\Phi_*(0,0)=0$ and $\GGG_0[\Phi]_*(0,0)=0$.
It follows in turn that properness and effectivity are
preserved by $\GGG_0$.
It follows similarly by (\ref{equation:RelationOfFTsTer}) that
Theorem~\ref{Theorem:RiggedStirling} remains valid with $\GGG$
replaced by $\GGG_0$.

\subsubsection{Dependence of the operator $\GGG_0$ on $\omega$}
Let us temporarily write $\mu_\omega$,
  $\GGG_{0,\omega}$ instead of $\mu$ and $\GGG_0$,
respectively,
in order
to keep track of dependence on $\omega$. Fix $u\in k^\times$.
We have 
$$\mu_{u\omega}=\norm{u}^{1/2}\mu_\omega,$$
hence
$$\GGG_{0,u\omega}[\varphi]=\norm{u}^{-1}\GGG_{0,\omega}[\varphi^{(u,u)}]
$$
for all $\varphi\in \SSS(k\times k)$, as one verifies by an evident transformation of integrals.
It follows that the analogous relation holds for the canonical extension of $\GGG_\omega$
to $\SSS'(k\times k)$, and in particular
\begin{equation}\label{equation:SeparabilityIndependence}
\GGG_{0,u\omega}[\Phi]_*(x,y)=\GGG_{0,\omega}[\Phi]_*(ux,uy)
\end{equation}
for all $\Phi\in \rigged(k)$, via (\ref{equation:RationalFourierScalingRule}).
Hereafter, until further notice, we revert to the system of notation
in which reference to $\omega$ is largely suppressed.

\subsubsection{Dependence of the operator $\GGG_0$ on $\FF$}
\label{subsubsection:pGNoughtDependence}
In this paragraph we temporarily write $\GGG_{0,\FF}$ instead of $\GGG_0$  in order to keep track of dependence on $\FF$. Consider
an intermediate field $\FF\subset\Gb\subset\FF_q$.
For each $c\in \FF$, let
$$V(c,\Gb,\FF)=\{v\in \Gb\mid \trace_{\Gb/\FF}v=c\}.$$
From Wilson's theorem $\prod_{C\in \FF^\times}=-1$ and its analogue for $\Gb$
we deduce that
\begin{equation}\label{equation:vIdentity}
c^{-1}\prod_{v\in V(c,\Gb,\FF)}v=\prod_{v\in V(1,\Gb,\FF)}v=1=\prod_{0\neq v\in V(0,\Gb,\FF)}v
\end{equation}
holds for all $c\in \FF^\times$.  We have a trivial identity
$$
\lambda_{0}(\trace_{\Gb/\FF}x)=
\sum_{v\in V(1,\Gb,\FF)}\lambda_{0}(v^{-1}x)-\sum_{0\neq v\in V(0,\Gb,\FF)}\lambda_{0}(v^{-1}x)
$$
holding for all $x\in \Gb$ which leads by the definitions
and scaling rule (\ref{equation:RiggedScaling})
 to the identity
\begin{equation}\label{equation:GpChange}
\begin{array}{cl}
&\displaystyle\GGG_{0,\FF}[\Phi]_*(x,y)\\\\
=&\displaystyle
\sum_{v\in V(1,\Gb,\FF)}\GGG_{0,\Gb}[\Phi]_*(v^{-1}x,v^{-1}y)-\sum_{0\neq v\in V(0,\Gb,\FF)}\GGG_{0,\Gb}[\Phi]_*(v^{-1}x,v^{-1}y)
\end{array}
\end{equation}
holding for all $\Phi\in \rigged(k)$.
 Hereafter, until further notice, we revert to the system of notation in which reference to $\FF$ is largely suppressed.

\subsubsection{Rationality and separability of rigged virtual lattices}
We say that $\Phi\in \rigged(k)$ is {\em rational} if $\Phi$ can be expressed as a finite $\ZZ[1/p]$-linear combination of primitive rigged virtual lattices.
By Lemma~\ref{Lemma:RationalRVFT}, the class of rational rigged virtual lattices
is stable under the operation $\GGG_0$. The class of such is stable also for every $b,c\in k^\times$ under the operation $\Phi\mapsto \Phi^{(b,c)}$.
Given rational $\Phi\in \rigged(k)$, we say that $\Phi$ is {\em separable}
if both $\Phi_*$ and $\GGG_0[\Phi]_*$ are $\ZZ$-valued.
For example, the rational primitive rigged virtual lattice
considered in Lemma~\ref{Lemma:RationalRVFT} is separable if
$\mu(r\OO)/\mu(k/L)\in \ZZ$. By (\ref{equation:SeparabilityIndependence}), the notion of
separability is independent of the choice of $\omega$.
By (\ref{equation:GpChange}) we are not permitted to conclude
that the notion of separability is independent of $\FF$,
but at least we can say that if $\FF\subset\Gb\subset\FF_q$,
then ``$\Gb$-separability'' implies ``$\FF$-separability''.

\subsection{The   Catalan symbol}
\label{subsection:CatalanFormalism}
\subsubsection{Definition}
\label{subsubsection:CatalanFormalism}
Recall that $k_\perf$ is
 the closure of $k$ under the extraction of $p^{th}$ roots in a fixed algebraic closure $\bar{k}$. Note that $\Phi_*(x,y)\in \ZZ[1/p]$ for
all rational $\Phi\in  \rigged(k)$ and $x,y\in k$.
Now given $a\in k^\times$ and a rational rigged virtual lattice $\Phi$, put
$$
\begin{array}{lcl}
\displaystyle\left(\begin{array}{c}
a\\
\Phi
\end{array}\right)_+&=&\displaystyle\prod_{x\in k^\times}x^{\Phi_*(x,a^{-1}x)}\in
 k_\perf^\times,\\\\
\displaystyle\left(\begin{array}{c}
a\\
\Phi
\end{array}\right)&=&\displaystyle
\prod_{x\in k^\times} x^{\Phi_*(x,a^{-1}x)}\cdot
\prod_{\xi\in k^\times}\xi^{\GGG_0[\Phi]_*(\xi,a\xi)\norm{a}}\in k_{\perf}^\times,
\end{array}
$$
thus defining 
the {\em  Catalan symbol} $\left(\begin{subarray}{c}
\cdot\\
\cdot
\end{subarray}\right)$
along with a partial version $\left(\begin{subarray}{c}
\cdot\\
\cdot
\end{subarray}\right)_+$ of it. The object $\left(\begin{subarray}{c}
\cdot\\
\cdot\end{subarray}\right)$ defined here is a local version of the global  Catalan symbol defined in  \cite{AndersonStirling}. 
\subsubsection{Basic formal properties}
Fix $a\in k^\times$ and rational $\Phi\in \rigged(k)$.
It is clear that $\left(\begin{subarray}{c}a\\\Phi\end{subarray}\right)_+$  and  
$\left(\begin{subarray}{c}a\\\Phi\end{subarray}\right)$  depend $\ZZ[1/p]$-linearly on $\Phi$. We have a scaling rule
\begin{equation}\label{equation:PartialCatalanScaling}
\left(\begin{array}{c}
a\\
\Phi^{(b,c)}
\end{array}\right)_+=b^{\norm{b}(\Theta(ac/b,\Phi)-\Phi_*(0,0))}
\left(\begin{array}{c}
ac/b\\
\Phi
\end{array}\right)_+^{\norm{b}}
\end{equation}
holding for all $b,c\in k^\times$,
by scaling rule (\ref{equation:RiggedScaling}) for the lower star operation and the definitions. Directly from the definitions we deduce a  factorization
\begin{equation}\label{equation:CatalanPartialFac}
\left(\begin{array}{c}
a\\
\Phi
\end{array}\right)=\left(\begin{array}{c}
a\\
\Phi
\end{array}\right)_+\left(\begin{array}{c}
a^{-1}\\
\GGG_0[\Phi]
\end{array}\right)_+^{\norm{a}}\end{equation}
and an implication
\begin{equation}\label{equation:SeparabilityImplication}
(\mbox{$\Phi$ is separable})\Rightarrow \left(\begin{array}{c}
a\\
\Phi
\end{array}\right)^{\max(1,\norm{a}^{-1})}\in k^\times.
\end{equation}
For the Catalan symbol we have a  functional equation
\begin{equation}\label{equation:CatalanFunctionalEquation}
\left(\begin{array}{c}
a^{-1}\\
\GGG_0[\Phi]
\end{array}\right)^{\norm{a}}=(-1)^{\Theta(a,\Phi)-\Phi_*(0,0)}\left(\begin{array}{c}
a\\
\Phi
\end{array}\right)
\end{equation}
holding by squaring rule (\ref{equation:RTSquare}), 
scaling rule (\ref{equation:PartialCatalanScaling}),
Wilson's theorem  $\prod_{C\in \FF^\times}C=-1$, and (\ref{equation:CatalanPartialFac}).
 Furthermore, we have a scaling rule
\begin{equation}\label{equation:CatalanScaling}
\left(\begin{array}{c}
a\\
\Phi^{(b,c)}
\end{array}\right)=b^{\norm{ac}\GGG_0[\Phi]_*(0,0)-\norm{b}\Phi_*(0,0)}\left(\begin{array}{c}
ac/b\\
\Phi
\end{array}\right)^{\norm{b}}\end{equation}
holding for all $b,c\in k^\times$,
by scaling rule (\ref{equation:RationalFourierScalingRule}),  
 functional equation (\ref{equation:NoughtThetaFE}),
scaling rule (\ref{equation:PartialCatalanScaling}), and (\ref{equation:CatalanPartialFac}).
Note that (\ref{equation:PartialCatalanScaling},\ref{equation:CatalanPartialFac},\ref{equation:CatalanFunctionalEquation},\ref{equation:CatalanScaling})
all simplify considerably if $\Phi$ is proper. 
We claim the following:
\begin{equation}\label{equation:CatalanLocalConstancy}
\mbox{The function $\left(a\mapsto \left(\begin{array}{c}a\\
\Phi
\end{array}\right)\right):k^\times\rightarrow k_\perf^\times$ is locally constant.}
\end{equation}
To prove the claim we may assume without loss of generality that $\Phi$ is of the form (\ref{equation:TypicalRigged}), and it is enough to prove the
analogous statement for the partial  Catalan symbol.
 Local constancy is then easy to check. Indeed, we have
 $$\Phi_*=\one_{\ell+L}\otimes \one_{w+r\OO}\Rightarrow\left(\begin{array}{c}
 a\\\Phi
 \end{array}\right)_+=
 \prod_{0\neq x\in (\ell+L)\cap a(w+r\OO)} x.$$
The claim is proved. We have
\begin{equation}\label{equation:DownAsymptoticPartialCatalan}
\left(\begin{array}{c}
a\\
\Phi
\end{array}\right)_+=1\;\;\mbox{for $\norm{a}$ sufficiently small.}
\end{equation}
This is also proved by reducing to the case 
(\ref{equation:TypicalRigged}). It follows finally that
\begin{equation}\label{equation:DownAsymptoticPartialCatalanBis}
\left(\begin{array}{c}
a\\
\Phi
\end{array}\right)=\left(\begin{array}{c}
a\\
\Phi
\end{array}\right)_+\;\;\mbox{for $\norm{a}\gg 1$}
\end{equation}
via factorization (\ref{equation:CatalanPartialFac}) and functional equation (\ref{equation:CatalanFunctionalEquation}).

\subsubsection{Simplified statement of Theorem~\ref{Theorem:RiggedStirling}}
Fix $a\in k^\times$ and proper rational \linebreak $\Phi\in  \rigged(k)$.
In this case the Stirling formula (\ref{equation:RiggedStirling}) for rigged virtual lattices 
simplifies a lot (especially, terms on the right side can be dropped),
and if we further simplify by writing the left side in terms of the  Catalan symbol (recall that by (\ref{equation:RelationOfFTsTer}) we may replace $\GGG$ by $\GGG_0$
in (\ref{equation:RiggedStirling}) without invalidating the latter)
we arrive at the more comprehensible relation
\begin{equation}\label{equation:SimplifiedRiggedStirling}
\begin{array}{cl}
&\displaystyle\Theta(a,\Phi)\ord \omega+\ord \left(\begin{array}{c}
a\\
\Phi
\end{array}\right)\\\\
=&\displaystyle\int \frac{\Theta(at,\Phi)-\one_{\OO^\times}(t)\Theta(a,\Phi)}{\norm{1-t}}d\mu^\times(t)\\\\
&\displaystyle+\int \Phi_*(0,t)d\mu^\times(t)+\norm{a}\int \GGG_0[\Phi]_*(0,t)
d\mu^\times(t).
\end{array}
\end{equation}
In the sequel we always apply Theorem~\ref{Theorem:RiggedStirling} 
in simplified form (\ref{equation:SimplifiedRiggedStirling}).

\subsubsection{Dependence of the Catalan symbol on $\omega$}
\label{subsubsection:omegaCatalanDependence}
 Fix $a\in k^\times$ and rational \linebreak $\Phi\in  \rigged(k)$.
Let us temporarily write $\GGG_{0,\omega}$ and $\left(\begin{subarray}{c}
a\\
\Phi
\end{subarray}\right)_\omega$ instead of $\GGG_0$
and $\left(\begin{subarray}{c}
a\\
\Phi
\end{subarray}\right)$, respectively,
in order to keep track of dependence on $\omega$. Fix $u\in k^\times$.
We have
\begin{equation}\label{equation:ChangeOfDifferential}
\left(\begin{array}{c}
a\\
\Phi
\end{array}\right)_{u\omega}=
u^{\GGG_{0,\omega}[\Phi]_*(0,0)\norm{a}-\Theta(a,\Phi)}\left(\begin{array}{c}
a\\
\Phi
\end{array}\right)_\omega
\end{equation}
by (\ref{equation:SeparabilityIndependence}), the definition of the Catalan symbol,
and the functional equation
 (\ref{equation:NoughtThetaFE}) for the theta symbol. 
 Suppose now that $\Phi$ is proper and let $\omega_0\in \Omega$
 be a fixed nonzero Chevalley differential.
It follows that the  product
$$\left(\omega/\omega_0\right)^{\Theta(a,\Phi)}\left(\begin{array}{c}
a\\
\Phi
\end{array}\right)_\omega$$
is  independent of $\omega$. As a consistency check,
note that the right side of formula (\ref{equation:SimplifiedRiggedStirling})
does not involve $\omega$ at all.
Hereafter, we revert to the system of
notation in which reference to $\omega$ is largely suppressed.

\subsubsection{Dependence of the Catalan symbol on $\FF$}
\label{subsubsection:pCatalanDependence}
In this paragraph we temporarily write 
$\GGG_{0,\FF}$ and $\left(\begin{subarray}{c}
a\\
\Phi
\end{subarray}\right)_\FF$ instead of $\GGG_{0}$ and $\left(\begin{subarray}{c}
a\\
\Phi
\end{subarray}\right)$, respectively, in order
to keep track of dependence on $\FF$.
By direct substitution into the definition of $(\begin{subarray}{c}
\cdot\\ \cdot
\end{subarray})_\FF$, we find via (\ref{equation:vIdentity}) and
(\ref{equation:GpChange}) that $(\begin{subarray}{c}
\cdot\\ \cdot
\end{subarray})_\FF=(\begin{subarray}{c}
\cdot\\ \cdot
\end{subarray})_\Gb$,  for any intermediate field $\FF\subset \Gb\subset\FF_q$, after a straightforward calculation which we omit. In other words, the Catalan symbol is independent of the choice of subfield $\FF\subset \FF_q$. Hereafter, we revert to the system of notation in which reference to $\FF$ is largely suppressed.

\section{Shadow theta and Catalan symbols}\label{section:InterpolatingGadgets}
We develop a second theory  
of symbols parallel to but independent of
that developed in the
preceding section. This second theory concerns objects built from the raw materials of local class field theory.
The main result of this section is
Theorem~\ref{Theorem:ShadowRiggedStirling}, which is parallel to Theorem~\ref{Theorem:RiggedStirling} in its simplified form (\ref{equation:SimplifiedRiggedStirling}). 
Basic notation introduced in \S\ref{subsection:BasicSetUp}
remains in force. 

\pagebreak

\subsection{Apparatus from local class field theory}

\subsubsection{The reciprocity law homomorphism $\rho$}
Let 
$$\rho:k^\times\rightarrow \Gal(k_\ab/k)$$
be the reciprocity law homomorphism of local class field theory,
``renormalized'' after  the modern fashion so that 
$$
C^{\rho(a)}=C^{\norm{a}}\;\;\;\;(C\in \FF_{q^\infty},\;\;a\in k^\times).
$$
Our normalization is the same as in \cite{TateBackground}, but the opposite of that in
\cite{LubinTate}.

\subsubsection{The Lubin-Tate description of $k_\ab/k$} \label{subsubsection:ExplicitReciprocity}
From \cite{LubinTate} we take the following
simple explicit description of the extension $k_\ab/k$. 
Fix a uniformizer $\pi\in\OO$. Let $X$ be a variable. Put 
$$[1]_\pi(X)=X\in \OO[X].$$
For $n>0$, put 
$$[\pi^n]_\pi(X)=\pi[\pi^{n-1}]_\pi(X)+([\pi^{n-1}]_\pi(X))^q\in \OO[X].$$
For general 
$$a=\sum_{i=0}^\infty a_i\pi^i\in \OO\;\;\;(a_i\in \FF_q),$$
put
$$[a]_\pi(X)=\sum_{i=0}^\infty a_i[\pi^n]_\pi(X)\in \OO[[X]].$$
Note that the power series $[a]_\pi(X)$ is
$\FF_q$-linear.
Put
$$\Xi_\pi=\{0\neq \xi\in \bar{k}\mid [\pi^n]_\pi(\xi)=0\;\mbox{for some positive integer $n$}\}.$$
Since the polynomials $[\pi^n]_\pi(X)$ have distinct roots, every element of $\Xi_\pi$ is separable over $k$.
The set $\Xi_\pi\cup\{0\}$ is a vector space over $\FF_q$ and becomes
an $\OO$-module isomorphic to $k/\OO$ when equipped with the $\OO$-action 
$$((a,\xi)\mapsto [a]_\pi(\xi)):\OO\times(\Xi_\pi\cup\{0\})\rightarrow
\Xi_\pi\cup\{0\}.$$
This $\OO$-action commutes with the action of the Galois group of the separable algebraic closure of $k$ in $\bar{k}$ and hence
$\Xi_\pi\subset k_\ab$.
According to the explicit reciprocity law of Lubin-Tate \cite{LubinTate},
we have
\begin{equation}\label{equation:ExplicitReciprocityLaw}
\rho(\pi)\xi=\xi,\;\;\;\rho(a)\xi=[a]_\pi(\xi)
\end{equation}
for all $a\in \OO^\times$ and $\xi\in \Xi_\pi$.
(Lest the reader be jarred by the seeming inconsistency of (\ref{equation:ExplicitReciprocityLaw})
with \cite[Formula 2, p.\ 380]{LubinTate},
we remind him/her that we renormalized the reciprocity law.) 
Now for every positive integer $n$ the polynomial
$$[\pi^n]_\pi(X)/[\pi^{n-1}]_\pi(X)$$
satisfies the Eisenstein criterion for irreducibility over $\OO$. It follows via
(\ref{equation:ExplicitReciprocityLaw})  that 
\begin{equation}\label{equation:Eisenstein}
\ord\xi=1/[k(\xi):k]=\mu^\times\{a\in \OO^\times \mid
\rho(a)\xi=\xi\}
\end{equation} 
for all $\xi\in \Xi_\pi$.

\pagebreak

\subsection{Path calculus}
\subsubsection{Paths} Let $X$ be a variable.
Put
$$\PPP=\left(\bigcup_{n=1}^\infty
X\FF_{q^n}[[X]]\right)\setminus\{0\}.$$
The set $\PPP$ forms a monoid with unit
under power series composition. We equip the monoid $\PPP$ with an
action of $\Gal(k_\ab/k)$ by declaring that
$$\sigma V=
\sum_{i=1}^\infty (\sigma a_i)X^i\;\;\;
\left(V=\sum_{i=1}^\infty a_iX^i\;\;(a_i\in \FF_{q^\infty})\right)$$
for all $\sigma\in \Gal(k_\ab/k)$ and $V\in \PPP$. We call elements of $\PPP$ {\em paths}.

\subsubsection{Basepoints}
Let $\PPP^*$ be the set of $\xi\in k_\ab$ satisfying the condition
$$\ord\xi=\mu^\times \{a\in \OO^\times \mid
\rho(a)\xi=\xi\}=1/\card\{\rho(a)\xi\mid a\in \OO^\times\}.$$
Clearly $\Gal(k_\ab/k)$ stabilizes the set $\PPP^*$.
Note that $\Xi_\pi\subset \PPP^*$
for every uniformizer $\pi\in \OO$ by (\ref{equation:Eisenstein}).
If $\xi\in \PPP^*$, then for every element $x$ of the field generated over $k$ by $\FF_{q^\infty}$ and $\xi$, there exists a positive integer $n$ such that $x$ can be expanded as a Laurent series in $\xi$ with coefficients in $\FF_{q^n}$.   
Given a pair $\xi,\eta\in \PPP^*$ such that $\xi$ belongs to the field generated over $k$ by $\FF_{q^\infty}$ and $\eta$, we write $\xi\leq \eta$.
Note that for each $\xi\in \PPP^*$ and uniformizer $\pi\in \OO$ there exists some
$\eta\in \Xi_\pi$ such that $\xi\leq \eta$.
Note that  any two elements of $\PPP^*$ have a common upper bound. 
We call elements of $\PPP^*$ {\em basepoints}.

\subsubsection{The special paths $U_{\xi,\eta}$}\label{subsubsection:PathRules}
Given a pair $\xi,\eta\in \PPP^*$ such that $\xi\leq \eta$, we note that there exists unique $$U_{\xi,\eta}=U_{\xi,\eta}(X)\in\PPP$$
such that $$U_{\xi,\eta}(\eta)=\xi.$$ 
In particular, 
$$U_{\xi,\xi}(X)=X.$$
We find it helpful to think of $U_{\xi,\eta}$
as a path leading from $\eta$ to $\xi$. 
We have
$$U_{\xi,\eta}(X)/X^{\frac{\ord\xi}{\ord\eta}}\in
\bigcup_{n=1}^\infty \FF_{q^n}[[X]]^\times.
$$
Since $k_\ab/k(\xi)$ is separable, we have
$$U'_{\xi,\eta}(X)=\frac{d}{dX}U_{\xi,\eta}(X)\neq 0.$$
We have $$U_{\sigma\xi,\sigma\eta}=
\sigma U_{\xi,\eta}$$
for all $\sigma\in \Gal(k_\ab/k)$ and
$$U_{\xi,\eta}(X)^{\norm{a}}=(\rho(a)U_{\xi,\eta})(X^{\norm{a}})$$
for all $a\in k^\times$ such that $\norm{a}\geq 1$. 
We have $$U_{\xi,\eta}\circ U_{\eta,\zeta}=U_{\xi,\zeta}$$
whenever we have basepoints $\xi\leq \eta\leq \zeta$. 
Combining the rules noted above, we have
$$
U_{\xi,\rho(a)\eta}(U_{\eta,\rho(b)\zeta}(X^{\norm{b}})^{\norm{a}})=
U_{\xi,\rho(ab)\zeta}(X^{\norm{ab}})
$$
for all $a,b\in k^\times$ and $\xi,\eta,\zeta\in \PPP^*$ such that $\xi\leq \eta\leq \zeta$.
Finally, note that
$$U_{\xi,\rho(a)\xi}(X^{\norm{a}})=U_{\xi,\rho(b)\xi}(X^{\norm{b}})
\Leftrightarrow (\norm{a}=\norm{b}\;\mbox{and}\;\rho(a)\xi=\rho(b)\xi)$$
for all $a,b\in k^\times$ such that $\min(\norm{a},\norm{b})\geq 1$
and $\xi\in \PPP^*$.  We take these rules for granted in the calculations to follow.

\begin{Subproposition}\label{Subproposition:PathFactorization} Fix $b,c\in k^\times$ such that
$\norm{b}\geq \norm{c}$. Fix $\xi,\eta\in \PPP^*$ such that $\xi\leq \eta$. Fix a set $S\subset \OO^\times$ of representatives for the quotient
$\frac{\{a\in\OO^\times\mid \rho(b)\xi=\xi\}}{\{a\in \OO^\times \mid \rho(a)\eta=\eta\}}
$.
Then we have
\begin{equation}\label{equation:PreWeierstrassNuff}
\begin{array}{cl}
&\displaystyle U_{\xi,\rho(b)\eta}(X^{\norm{b}})-
U_{\xi,\rho(c)\eta}(Y^{\norm{c}})\\\\
=&\displaystyle V(X,Y)\prod_{s\in S}
(U_{\eta,\rho(sb/c)\eta}(X^{\norm{b/c}})-Y)^{\norm{c}}
\end{array}
\end{equation}
for some $V(X,Y)\in \bigcup_{n=1}^\infty \FF_{q^n}[[X,Y]]^\times$.
 \end{Subproposition}
 \proof Since
$$
U_{\xi,\rho(b)\eta}(X^{\norm{b}})-
U_{\xi,\rho(c)\eta}(Y^{\norm{c}})=
(U_{\rho(c)^{-1}\xi,\rho(b/c)\eta}(X^{\norm{b/c}})-
U_{\rho(c)^{-1}\xi,\eta}(Y))^{\norm{c}},
$$
we may assume without  loss of generality that $c=1$.
Fix a positive integer $n$ such that 
 $$\{U_{\xi,\rho(b)\eta}\}\cup \{U_{\xi,\eta}\}
 \cup\{U_{\eta,\rho(sb)\eta}\mid s\in S\}\subset \FF_{q^n}[[X]].$$
By the Weierstrass Division Theorem we have a unique factorization
$$U_{\xi,\rho(b)\eta}(X^{\norm{b}})-U_{\xi,\eta}(Y)=V(X,Y)W(X,Y)$$
where  
$$V(X,Y)\in \FF_{q^n}[[X,Y]]^\times,$$
$$\mbox{$W(X,Y)\in \FF_{q^n}[[X]][Y]$ is monic in $Y$ of degree $e=\frac{\ord \xi}{\ord \eta}$, and}$$
$$W(0,Y)=Y^e.$$  It suffices to prove that
\begin{equation}\label{equation:WeierstrassNuff}
(-1)^eW(X,Y)= \prod_{s\in S} (U_{\eta,\rho(sb)\eta}(X^{\norm{b}})-Y).
\end{equation}
Since for each $s\in S$ we have
$$\begin{array}{rcl}
U_{\xi,\eta}(U_{\eta,\rho(sb)\eta}(X^{\norm{b}}))&=&
U_{\xi,\rho(sb)\eta}(X^{\norm{b}})\\\\
&=&(\rho(s)U_{\rho(s)^{-1}\xi,\rho(b)\eta})(X^{\norm{b}})\\\\
&=&U_{\xi,\rho(b)\eta}(X^{\norm{b}}),
\end{array}
$$
each factor on the right side of (\ref{equation:WeierstrassNuff}) divides 
the left side of (\ref{equation:WeierstrassNuff}) in the ring $\FF_{q^n}[[X]][Y]$.
Since for $s$ ranging over $S$ the power series
$U_{\eta,\rho(sb)\eta}(X^{\norm{b}})$ are distinct,
the left side of (\ref{equation:WeierstrassNuff})
divides the right side of (\ref{equation:WeierstrassNuff})
in the ring $\FF_{q^n}[[X]][Y]$.
Finally, since $e$ equals the cardinality of $S$ by definition of $\PPP^*$, equality indeed holds in (\ref{equation:WeierstrassNuff}). \qed

\pagebreak
\subsection{Two-variable extension of path calculus}
\subsubsection{The rings $\Lambda$
and $\Lambda_0$}
Let $X$ and $Y$ be independent variables. Put
$$\Lambda_0=\bigcup_{n=1}^\infty
\FF_{q^n}[[X,Y]],\;\;\;
\Lambda=\Lambda_0[X^{-1},Y^{-1}].$$
Given nonzero $F,G\in \Lambda$, we write 
$F\sim G$ if there exists $H\in \Lambda_0^\times$ such that $F=GH$.
Given any $\alpha,\beta\in \bar{k}$ such that $0<|\alpha|<1$ and $0<|\beta|<1$, note that each power series $F=F(X,Y)\in \Lambda$ specializes to a value $F(\alpha,\beta)\in \bar{k}$,
and  that
if $\alpha,\beta\in k_\ab$, then $F(\alpha,\beta)\in k_\ab$.  

\begin{Subproposition}\label{Subproposition:LambdaStart}
$\Lambda$
is a principal ideal domain.
\end{Subproposition}
\proof 
Given an integer $d\geq 0$ and 
$W=W(X,Y)\in \Lambda$, let us say that $W$ is {\em distinguished of degree $d$}
if $W(X,Y)\in \Lambda_0$, $W(X,Y)$ is a monic polynomial in $Y$ of degree $d$,
$W(0,Y)=Y^d$, and $W(X,0)\neq 0$.  Now suppose we are given nonzero $F\in \Lambda$. There are unique integers $\mu$ and $\nu$ and unique $F_0\in \Lambda_0$
divisible in $\Lambda_0$ by neither $X$ nor $Y$ such that $F=X^\mu Y^\nu F_0$. 
Choose $n$ such that $F_0\in \FF_{q^n}[[X,Y]]$.
By the Weierstrass Division Theorem there exists unique distinguished $F_1\in \FF_{q^n}[[X,Y]]$ and unique  $F_2\in \FF_{q^n}[[X,Y]]^\times$ 
such that $F_0=F_1F_2$, and moreover $F_1$ and $F_2$ are independent of the choice of $n$.  We declare the degree in $Y$  of $F_1$ to be the {\em $\Lambda$-degree} of $F$.  If $F$ is a polynomial in $Y$, then the $\Lambda$-degree of $F$ cannot exceed the degree in $Y$ of $F$. 
Now let there be given a nonzero ideal $I\subset \Lambda$.
Let $0\neq F\in I$ be of minimal possible $\Lambda$-degree. 
Choose $G\in I$ arbitrarily. We claim that $F$ divides $G$.
To prove the claim, after multiplying both $F$ and $G$ by suitably chosen invertible elements of $\Lambda$, we may assume without loss of generality that $F$
is distinguished and that $G$ is a polynomial in $Y$. 
After dividing $G$ by $F$ according to the
usual polynomial division algorithm, we may assume that $G$ is of degree
in $Y$ strictly less than that of $F$, in which case $G$ must vanish identically---otherwise we have a contradiction to minimality of the $\Lambda$-degree of $F$.
 The claim is proved,
and with it the proposition. \qed

\subsubsection{Valuations of $\Lambda$}
Given nonzero $F\in \Lambda$, we  define the $X$- and {\em $Y$-valuations} of  $F$ to be the unique integers $m$ and $m$, respectively,
such that $X^{-m}Y^{-n}F$ belongs to $\Lambda_0$
but is divisible in $\Lambda_0$ by neither $X$ nor $Y$. 
Given $P\in \Lambda_0$ remaining irreducible in $\Lambda$ (e.g., any power series $P$ of the form $X-f(Y)$ or $f(X)-Y$ with $f\in \PPP$) and
nonzero $F\in \Lambda$, we define the {\em $P$-valuation}
of $F$ to be the multiplicity with which $P$ divides $F$ in $\Lambda$, which is
well-defined since $\Lambda$ is a principal ideal domain.

\subsubsection{Natural operations on $\Lambda$}\label{subsubsection:LambdaOps}
We equip $\Lambda$ with an action of  $\Gal(k_\ab/k)$ stabilizing $\Lambda_0$ by declaring that
$$\sigma F=\sum_i\sum_j (\sigma a_{ij})X^iY^i\;\;\;
\left(F=\sum_i\sum_j a_{ij}X^iY^j\;\;(a_{ij}\in \FF_{q^\infty})\right)$$
for all $\sigma\in \Gal(k_\ab/k)$ and $F\in \Lambda$.
Given $F=F(X,Y)\in \Lambda$, put 
$$F^\dagger(X,Y)=F(Y,X).$$
Given also $\xi,\eta\in \PPP^*$ such that $\xi\leq \eta$
and $b,c\in k^\times$ such that 
$\min(\norm{b},\norm{c})\geq 1$, we define
$$\begin{array}{rcl}
F^{[b,c;\xi,\eta]}(X,Y)&=&F(U_{\xi,\rho(b)\eta}(X^{\norm{b}}),U_{\xi,\rho(c)\eta}(Y^{\norm{c}}))\\
&=&F(U_{\rho(b)^{-1}\xi,\eta}(X)^{\norm{b}},U_{\rho(c)^{-1}\xi,\eta}(X)^{\norm{c}}),
\end{array}$$
which belongs again to $\Lambda$.
To abbreviate notation we write
$$F^{[b,c;\xi]}=F^{[b,c;\xi,\xi]},\;\;\;F^{[b;\xi,\eta]}=F^{[b,b;\xi,\eta]},\;\;\;F^{[b;\xi]}=F^{[b,b;\xi]},$$
i.~e., we drop doubled letters when possible. 
From the rules noted in \S\ref{subsubsection:PathRules} we deduce the following rules obeyed by
the ``square bracket'' and ``dagger'' operations.
We have
$$
(F^{[b,c;\xi,\eta]})^\dagger=(F^\dagger)^{[c,b;\xi,\eta]}.
$$
We have
$$F^{[b,c;\xi,\eta]}(\eta,\eta)=F((\rho(b)^{-1}\xi)^{\norm{b}},(\rho(c)^{-1}\xi)^{\norm{c}}).$$
Given also $\sigma\in \Gal(k_\ab/k)$, we have
$$
\sigma(F^{[b,c;\xi,\eta]})=(\sigma F)^{[b,c;\sigma\xi,\sigma\eta]}.
$$
Given also $a\in k^\times$, we have
$$
\begin{array}{rcl}
(F^{[b,c;\xi,\eta]}(X,Y))^{\norm{a}}&=&(\rho(a)F)^{[b,c;\rho(a)\xi,\rho(a)\eta]}(X^{\norm{a}},Y^{\norm{a}})\\
&=&(\rho(a)F)^{[ab,ac;\rho(a)\xi,\eta]}(X,Y).\end{array}
$$
Given also $d,e\in k^\times$ such that $\norm{d},\norm{e}\geq 1$ and $\zeta\in \PPP^*$ such that $\eta\leq \zeta$, we have
 $$
 (F^{[b,c;\xi,\eta]})^{[d,e;\eta,\zeta]}=F^{[bd,ce;\xi,\zeta]}.
$$
We take these rules for granted in the calculations to follow.

\subsubsection{The power series $Z_{t,\xi}$}
Given $\xi\in\PPP^*$ and $t\in k^\times$, put
$$
Z_{t,\xi}=\left\{\begin{array}{ll}
X^{\norm{t}}-U_{\rho(t)\xi,\xi}(Y)&\mbox{if $\norm{t}\geq 1$,}\\
Y^{\norm{t^{-1}}}-U_{\rho(t^{-1})\xi,\xi}(X)&\mbox{if $\norm{t}<1$,}
\end{array}\right.
$$
which is a power series belonging to the ring $\Lambda_0$. Note that
$$Z_{1,\xi}=X-Y.$$
Note that
$$\sigma Z_{t,\xi}=Z_{t,\sigma\xi}
$$
for all $\sigma\in  \Gal(k_\ab/k)$.
We claim that
$$
Z_{t,\xi}\sim 
\left\{\begin{array}{ll}
(X-Y)^{[t,1;\xi]}&\mbox{if $\norm{t}\geq 1$,}\\
(X-Y)^{[1,t^{-1};\xi]}&\mbox{if $\norm{t}\leq 1$.}
\end{array}\right.
$$
It is enough by the Weierstrass Division Theorem to show that the power series on the right divides the power series on the left.
In the case $\norm{t}>1$, we have such divisibility because $U_{\rho(t)\xi,\xi}(U_{\xi,\rho(t)\xi}(X^{\norm{t}}))=X^{\norm{t}}$. The remaining cases of the claim are handled similarly. Thus the claim is proved.
The symmetry
$$
Z_{t,\xi}^\dagger\sim Z_{t^{-1},\xi}
$$
follows. 
Denoting by $(Z_{t,\xi})$ the ideal of the principal ideal domain $\Lambda$ generated by $Z_{t,\xi}$, it is clear that $(Z_{t,\xi})$ is prime, and moreover  that 
$$
(Z_{t,\xi})=(Z_{t',\xi})\Leftrightarrow Z_{t,\xi}=Z_{t',\xi}\Leftrightarrow (\norm{t}=\norm{t'}\;\mbox{and}\;
\rho(t)\xi=\rho(t')\xi).
$$
We take these properties of $Z_{t,\xi}$ for granted in the calculations to follow.

\begin{Subproposition}\label{Subproposition:ZFac}
Fix $b,c\in k^\times$ such that $\min(\norm{b},\norm{c})\geq 1$.
Fix $\xi,\eta\in \PPP^*$ such that $\xi\leq \eta$.
Fix a set $S$ of representatives for the quotient 
$\frac{\{a\in \OO^\times \mid \rho(a)\xi=\xi\}}{\{a\in \OO^\times \mid \rho(a)\eta=\eta\}}$.
We have
\begin{equation}\label{equation:KeyZFac}
Z_{t,\xi}^{[b,c;\xi,\eta]}\sim \prod_{s\in S}
Z_{stb/c,\eta}^{\min(\norm{bt},\norm{c})/\min(\norm{t},1)}
\end{equation}
for every $t\in k^\times$.  We also have
\begin{equation}\label{equation:XYFac}
X^{[b,c;\xi,\eta]}\sim X^{\norm{b}\frac{\ord\xi}{\ord\eta}},\;\;\;
Y^{[b,c;\xi,\eta]}\sim Y^{\norm{c}\frac{\ord\xi}{\ord\eta}}.
\end{equation}
For all nonzero $F\in \Lambda$ we have
\begin{equation}\label{equation:ValuationLift}
(\mbox{$(X-Y)$-valuation of $F^{[b,c;\xi,\eta]}$})=
\max(\norm{b},\norm{c})\cdot (\mbox{$Z_{c/b,\xi}$-valuation of $F$}).
\end{equation}
\end{Subproposition}
\proof
Proposition~\ref{Subproposition:PathFactorization} in the present system of notation takes the form
$$\norm{b}\geq \norm{c}\Rightarrow
(X-Y)^{[b,c;\xi,\eta]}\sim
\prod_{s\in S}((X-Y)^{[sb/c,1;\eta]})^{\norm{c}}.$$
Symmetrically we have
$$\norm{b}\leq \norm{c}
\Rightarrow (X-Y)^{[b,c;\xi,\eta]}\sim
\prod_{s\in S}((X-Y)^{[1,s^{-1}c/b;\eta]})^{\norm{b}}.$$
Formula (\ref{equation:KeyZFac}) now follows after a slightly tedious case analysis the further details of which we omit.
Statement (\ref{equation:XYFac}) is clear. We record it here for convenient reference
since it is usually applied in conjunction with (\ref{equation:KeyZFac}).
We turn finally to the proof of (\ref{equation:ValuationLift}), which consists of a study of divisibility properties in the principal ideal domain $\Lambda$.
Let $\ell$ be the $Z_{c/b,\xi}$-valuation of $F$ and write
$F=Z_{c/b,\xi}^\ell W$, where $Z_{c/b,\xi}$ and $W$ are relatively prime.
It follows that $Z_{c/b,\xi}^{[b,c;\xi,\eta]}$ and
$W^{[b,c;\xi,\eta]}$ are also relatively prime.
 By (\ref{equation:KeyZFac}) we can write
$Z_{c/b,\xi}^{[b,c;\xi,\eta]}=(X-Y)^{\max(\norm{b},\norm{c})}V$, where $X-Y$ and $V$ are relatively prime.
We conclude that
$$F^{[b,c;\xi,\eta]}=(X-Y)^{\ell\cdot \max(\norm{b},\norm{c})}W^{[b,c;\xi,\eta]}V^\ell$$
where
$X-Y$ and $W^{[b,c;\xi,\eta]}V^\ell$ are relatively prime, which proves (\ref{equation:ValuationLift}). \qed

\begin{Subproposition}\label{Subproposition:RingHomomorphism}
Fix $b,c\in k^\times$ such that $\min(\norm{b},\norm{c})\geq 1$.
Fix $\xi,\eta\in \PPP^*$ such that $\xi\leq \eta$.
Let $\Lambda'$ be a copy of $\Lambda$ viewed as a $\Lambda$-algebra via the homomorphism $F\mapsto F^{[b,c;\xi,\eta]}$.
(i) $\Lambda'$ is a free $\Lambda$-module of rank $\norm{bc}\left(\frac{\ord\xi}{\ord\eta}\right)^2$.
(ii) If $\norm{b}=\norm{c}=1$, then $\Lambda'/\Lambda$ is a Galois extension with Galois group 
isomorphic to the product of two copies of the quotient $\frac{\{u\in \OO^\times\mid \rho(u)\eta=\eta\}}{
\{u\in \OO^\times \mid \rho(u)\xi=\xi\}}$. For all $F\in \Lambda$ we have
(iii) $F=0\Leftrightarrow F^{[b,c;\xi,\eta]}=0$ and (iv) $F\in \Lambda^\times\Leftrightarrow F^{[b,c;\xi,\eta]}\in \Lambda^\times$.
 \end{Subproposition}
 \proof 
  (i) Let $\Lambda_0'$ be a copy of $\Lambda_0$
  viewed as $\Lambda_0$-algebra via the ring homomorphism  $F\mapsto F^{[b,c;\xi,\eta]}$. In view of formula (\ref{equation:XYFac}),
it is enough to prove that $\Lambda'_0$ is a free \linebreak $\Lambda_0$-module of rank $\norm{bc}\left(\frac{\ord\xi}{\ord\eta}\right)^2$. Let $V$ and $W$ be independent variables independent also of $X$ and $Y$. Let $N$ be a positive integer such that $U_{\xi,\sigma \eta}(X)\in \FF_{q^N}[[X,Y]]$ for all $\sigma\in \Gal(k_\ab/k)$.
For any positive integer $n$ divisible by $N$ the quotient
$$M_n=\FF_{q^n}[[X,Y,V,W]]/(X-U_{\xi,\rho(b)\eta}(V^{\norm{b}}),
Y-U_{\xi,\rho(c)\eta}(W^{\norm{c}}))$$
is a free $\FF_{q^n}[[V,W]]$-module with basis $\{1\}$ and a free $\FF_{q^n}[[X,Y]]$-module with basis
$$\left\{V^iW^j\left| i=0,\dots\norm{b}\frac{\ord\xi}{\ord\eta}-1,
j=0,\dots,\norm{c}\frac{\ord\xi}{\ord\eta}-1\right.\right\}$$ by the Weierstrass Division Theorem.  It follows that $M=\bigcup_{N\mid n} M_n$ is a free \linebreak $\Lambda_0$-module of rank $\norm{bc}\left(\frac{\ord\xi}{\ord\eta}\right)^2$. Since $M$ is a $\Lambda_0$-module isomorphic to $\Lambda'_0$, we are done.

(ii) 
The discriminant of the finite flat extension $\Lambda'/\Lambda$ of principal ideal domains factors as a product of a function of $X$ times a function of $Y$ and hence is a unit of $\Lambda$.
Thus $\Lambda'/\Lambda$ is finite \'{e}tale.
Call the group in question $G$. 
The group $G$ acts faithfully on the ring extension $\Lambda'/\Lambda$ by the ``square bracket'' operation,
and since $G$ has cardinality
equal to the rank of the extension $\Lambda'/\Lambda$,
in fact $G$ must be the Galois group.

(iii) By (i) we may identify $\Lambda$ with a subring of $\Lambda'$, whence the result.

(iv)  By (i) and the Cohen-Seidenberg theorem,
every maximal ideal of $\Lambda$ lies below some
maximal ideal of $\Lambda'$, whence the result.
\qed

\subsection{Interpolating gadgets and Coleman units}
\subsubsection{Interpolating gadgets}
Given a pair $(F,\xi)$ where $0\neq F\in \Lambda$
and $\xi\in \PPP^*$, we call $(F,\xi)$ an {\em interpolating
gadget} if $F^{[a;\xi]}=F^{\norm{a}}$ for all $a\in k^\times$
such that $\norm{a}\geq 1$.

\begin{Subproposition}\label{Subproposition:InterpolatingGadget}
Fix a basepoint $\xi\in \PPP^*$. Fix nonzero $F\in \Lambda$.
The following conditions are equivalent:
\begin{enumerate}
\item $(F,\xi)$ is an interpolating gadget.
\item $(F^{[b,c;\xi,\eta]},\eta)$ is an interpolating gadget
for every $\eta\in \PPP^*$  and $b,c\in k^\times$ such that $\min(\norm{b},\norm{c})\geq 1$.
\item $(F^{[b,c;\xi,\eta]},\eta)$ is an interpolating gadget
for some $\eta\in \PPP^*$  and $b,c\in k^\times$ such that $\min(\norm{b},\norm{c})\geq 1$.
\item $(\sigma F)^{[1;\sigma\xi,\xi]}=F$
for all $\sigma\in \Gal(k_\ab/k)$.
\item $F(\xi,(\rho(a)^{-1}\xi)^{\norm{a}})\in k$ for all $a\in k^\times$
such that $\norm{a}\gg 1$.
\end{enumerate}
\end{Subproposition}
\proof The scheme of proof is
(1$\Rightarrow$2$\Rightarrow$3$\Rightarrow$1$\Rightarrow$4$
\Rightarrow$1$\Rightarrow$5$\Rightarrow$1). To prove the implication (1$\Rightarrow$2), just note that
$$(F^{[b,c;\xi,\eta]})^{[a;\eta]}=
(F^{[a;\xi]})^{[b,c;\xi,\eta]}=(F^{[b,c;\xi,\eta]})^{\norm{a}}$$
for all $a\in k^\times$ such that $\norm{a}\geq 1$.
The implication (2$\Rightarrow$3) is trivial.
To prove the implication (3$\Rightarrow$1), just note that for any $a\in k^\times$ such that $\norm{a}\geq 1$ we have
$$(F^{[a;\xi]}-F^{\norm{a}})^{[b,c;\xi,\eta]}=
(F^{[b,c;\xi,\eta]})^{[a;\eta]}-(F^{[b,c;\xi,\eta]})^{\norm{a}}=0$$
and hence $F^{[a;\xi]}-F^{\norm{a}}=0$
by
Proposition~\ref{Subproposition:RingHomomorphism}.
To prove the equivalence (1$\Leftrightarrow$4), fix $F\in \Lambda$,
$\sigma\in \Gal(k_\ab/k)$ and $a\in k^\times$ such that
$$\rho(a)^{-1}\xi=\sigma \xi,\;\;\rho(a)^{-1}F=\sigma F.$$
Then we have
$$F^{[a;\xi]}=((\rho(a)^{-1}F)^{[1;\rho(a)^{-1}\xi,\xi]})^{\norm{a}}=
((\sigma F)^{[1;\sigma \xi,\xi]})^{\norm{a}}.$$
Under hypothesis 4, the last term equals $F^{\norm{a}}$.
Under hypothesis 1 the first term equals $F^{\norm{a}}$.
Thus the equivalence (1$\Leftrightarrow$4) is proved.
We turn next to the proof of the implication (1$\Rightarrow$5). Fix $\sigma\in \Gal(k^\ab/k)$ and $a\in k^\times$ such that $\norm{a}\geq 1$. Let $G=F^{[1,a;\xi]}$, in which case
$F(\xi,(\rho(a)^{-1}\xi)^{\norm{a}})=G(\xi,\xi)$. Note that 
since (1$\Rightarrow$2) holds, $(G,\xi)$
is an interpolating gadget, and that since (1$\Rightarrow$4) holds, $(G,\xi)$
has property 4. We now have
$$\begin{array}{rcl}
\sigma(F(\xi,(\rho(a)^{-1}\xi)^{\norm{a}}))&=&\sigma(G(\xi,\xi))=(\sigma G)(\sigma\xi,\sigma\xi)\\
&=&(\sigma G)^{[1;\sigma\xi,\xi]}(\xi,\xi)=G(\xi,\xi)=F(\xi,(\rho(a)^{-1}\xi)^{\norm{a}}),
\end{array}
$$
and hence $F(\xi,(\rho(a)^{-1}\xi)^{\norm{a}})\in k$.
Thus the implication (1$\Rightarrow$5) is proved.
We turn finally to the proof of the implication (5$\Rightarrow$1). Fix $a,b\in k^\times$ with $\min(\norm{a},\norm{b})\geq 1$ and choose $\norm{a}$ large enough so that
$F(\xi,(\rho(a)^{-1}\xi)^{\norm{a}})\in k$. As before, put $G=F^{[1,a;\xi]}$, so that $G(\xi,\xi)=F(\xi,(\rho(a)^{-1}\xi)^{\norm{a}})$. We  have
$$\begin{array}{rcl}
F^{[b;\xi]}(\xi,(\rho(a)^{-1}\xi)^{\norm{b}})&=&F^{[b,ab;\xi]}(\xi,\xi)=G^{[b;\xi]}(\xi,\xi)\\
&=&
G((\rho(b)^{-1}\xi)^{\norm{b}},(\rho(b)^{-1}\xi)^{\norm{b}})
\\
&=&(\rho(b)^{-1}G)(\rho(b)^{-1}\xi,\rho(b)^{-1}\xi)^{\norm{b}}\\
&=&\rho(b)^{-1}G(\xi,\xi)^{\norm{b}}=G(\xi,\xi)^{\norm{b}}\\
&=&F(\xi,(\rho(a)^{-1}\xi)^{\norm{a}})^{\norm{b}}.
\end{array}
$$
It follows that  $F^{[b;\xi]}-F^{\norm{b}}$ belongs to infinitely many maximal ideals
of $\Lambda$ and hence vanishes identically. Thus $(F,\xi)$ is an interpolating gadget.
Thus the implication (5$\Rightarrow$1) is proved. 
\qed
\begin{Subproposition}\label{Subproposition:TotallyRamified}
Given a uniformizer $\pi\in \OO$,
a basepoint $\xi\in \Xi_\pi\subset \PPP^*$, and an interpolating gadget $(F,\xi)$, we have
$F\in \FF_q[[X,Y]][X^{-1},Y^{-1}]$.
\end{Subproposition}
\proof We have $\rho(\pi)F=(\rho(\pi)F)^{[1;\rho(\pi)\xi,\xi]}=F$
by  the explicit reciprocity law (\ref{equation:ExplicitReciprocityLaw})
and the fourth characterization of interpolating gadgets stated in the preceding proposition. \qed

\subsubsection{Coleman units}
Given $F\in \Lambda$ and $\xi\in \PPP^*$, we say that $F\in \Lambda$ is a {\em Coleman unit} at $\xi$
if $F$ factors as a product of elements of the set $\{Z_{t,\xi}\mid t\in k^\times\}$ times a unit of $\Lambda$.
We say that an interpolating gadget $(F,\xi)$
is of {\em Coleman type} if $F$ is a Coleman unit at $\xi$. The rationale for the terminology comes ultimately from Coleman's marvelous paper \cite{Coleman}. We have previously discussed 
our interest in (and admiration for) Coleman's paper in our papers \cite{AndersonStickelberger}
and \cite{AndersonStirling}. We refer the reader to the latter for more discussion.

\begin{Subproposition}\label{Subproposition:ColemanUnits}
Fix $\xi\in \PPP^*$ and $F\in \Lambda$.
The following statements are equivalent:
\begin{enumerate}
\item $F$ is a Coleman unit at $\xi$.
\item For all $\eta\in \PPP^*$ such that $\xi\leq \eta$
and $b,c\in k^\times$ such that $\min(\norm{b},\norm{c})\geq 1$,
the power series
 $F^{[b,c;\xi,\eta]}$ is a Coleman unit at $\eta$.
\item For some $\eta\in \PPP^*$ such that 
$\xi\leq \eta$ and $b,c\in k^\times$ such that $\min(\norm{b},\norm{c})\geq 1$, the power series $F^{[b,c;\xi,\eta]}$ is  a Coleman unit at $\eta$.
\end{enumerate}
\end{Subproposition}
\proof (1$\Rightarrow$2) This follows immediately from Proposition~\ref{Subproposition:ZFac}.  (2$\Rightarrow$3) Trivial.
(3$\Rightarrow$1) In any case, we have a factorization $F=GH$ where $H$ is a Coleman unit at $\xi$ and $G$ is relatively prime to $Z_{t,\xi}$ for all $t\in k^\times$. It follows that $G^{[b,c;\xi,\eta]}$ is relatively
prime to $Z_{t,\xi}^{[b,c;\xi,\eta]}$ for all $t\in k^\times$.
By Proposition~\ref{Subproposition:ZFac}
it follows that $G^{[b,c;\xi,\eta]}$ is relatively
prime to $Z_{t,\eta}$
for all $t\in k^\times$. But $F^{[b,c;\xi,\eta]}$ is by hypothesis
a Coleman unit and $G^{[b,c;\xi,\eta]}$ divides $F^{[b,c;\xi,\eta]}$.
Therefore we have $G^{[b,c;\xi,\eta]}\in\Lambda^\times$. By Proposition~\ref{Subproposition:RingHomomorphism}
it follows that $G\in \Lambda^\times$.
Therefore $F$ is a Coleman unit.
\qed
\begin{Subproposition}\label{Subproposition:GaloisDescent}
Fix an interpolating gadget $(G,\eta)$.
Fix $\xi\in \PPP^*$ such that $\xi\leq \eta$.
Assume that $G^{[u,v;\eta]}=G$ for all
$u,v\in\OO^\times$ such that $\rho(u)\xi=\xi=\rho(v)\xi$. Then there exists
a unique interpolating gadget $(F,\xi)$ with basepoint $\xi$ such that
$F^{[1;\xi,\eta]}=G$. Moreover, if $(G,\eta)$
is of Coleman type, then so is $(F,\xi)$.
\end{Subproposition}
\proof Existence and uniqueness of $F\in \Lambda$ such that
$F^{[1;\xi,\eta]}=G$ we get from part (ii)
of Proposition~\ref{Subproposition:RingHomomorphism}. One verifies that $(F,\xi)$
is an interpolating gadget by Proposition~\ref{Subproposition:InterpolatingGadget}.
Assuming that $(G,\eta)$ is of Coleman type,
one verifies that $(F,\xi)$ is of Coleman type by Proposition~\ref{Subproposition:ColemanUnits}.
\qed

\subsection{Definition and basic formal properties of the shadow symbols}

\begin{Subproposition}\label{Subproposition:RaisonDetreBis}
Fix an interpolating gadget $(F,\xi)$. Fix $\eta\in \PPP^*$ such that $\xi\leq \eta$. Fix $b,c\in k^\times$ such that $\min(\norm{b},\norm{c})\geq 1$. Put
$$\ell= (\mbox{$(X-Y)$-valuation of $F^{[b,c;\xi,\eta]}$}).$$
(i) Then $\ell/\norm{b}$ is  independent of $\eta$ and
depends only the ratio $c/b$.
(ii) Furthermore, the expression
\begin{equation}\label{equation:DiagonalExpression}
\left(\left(\frac{d\pi}{\omega}\right)^\ell\cdot
\left(\frac{1}{U_{\pi,\eta}'(\eta)}\right)^{\ell}\cdot \left(\left.\frac{F^{[b,c;\xi,\eta]}(X,Y)}{(X-Y)^\ell}\right|_{X=Y=\eta}\right)\right)^{1/\norm{b}}
\end{equation}
 defines a nonzero element of $k_\perf$ which is independent of $\eta$ and $\pi$, and which depends only on the ratio $c/b$.
\end{Subproposition}
\proof 
For each nonzero $H\in \Lambda$, in the unique possible way, write
$$H\equiv \epsilon(H)\cdot (X-Y)^{\delta(H)}\bmod{(X-Y)^{\delta(H)+1}}
$$
where
$$\delta(H)\in \ZZ\cap[0,\infty),\;\;\;\epsilon(H)=\epsilon(H)(Y)\in \bigcup_{n=1}^\infty\FF_{q^n}[[Y]][Y^{-1}]^\times.$$
 In other words, $\delta(H)$ is the order of vanishing of $H$ along the ideal $(X-Y)$,
and $\epsilon(H)$ is the leading coefficient of the Taylor expansion of $H$ in powers of $X-Y$
with coefficients in $\bigcup_{n=1}^\infty\FF_{q^n}[[Y]][Y^{-1}]$.
The $\delta$ and $\epsilon$ notation defined here will only be used in this  proof, not elsewhere in the paper.  It is convenient to set
$$G=F^{[b,c;\xi]}.$$
Note that $(G,\xi)$ is again an interpolating gadget
by Proposition~\ref{Subproposition:InterpolatingGadget}.

(i) By the definitions of $\ell$ and $G$,
and relation (\ref{equation:ValuationLift}) of Proposition~\ref{Subproposition:ZFac},
 we have
 $$
\ell=\delta(F^{[b,c;\xi,\eta]})=\delta(G^{[1;\xi,\eta]})=\delta(G).
$$
Thus $\ell$ is independent of $\eta$.
 By definition of interpolating gadget we have
$$\delta(F^{[ab,ac;\xi]})=\delta(G^{[a;\xi]})=
\delta(G^{\norm{a}})=\norm{a}\delta(G)$$
for all $a\in k^\times$ such that $\norm{a}\geq 1$. 
Thus $\ell/\norm{b}$ depends only on the ratio $c/b$.

(ii) Let $r$ denote the number
defined by (\ref{equation:DiagonalExpression}). 
By the Chain Rule in the form
$$
U'_{\pi,\eta}(X)=U'_{\pi,\xi}(U_{\xi,\eta}(X))U'_{\xi,\eta}(X),
$$
we have
\begin{equation}\label{equation:DiagonalExpressionBis}
\left(U'_{\pi,\xi}(\xi)\frac{\omega}{d\pi}\right)^{\delta(G)}r^{\norm{b}}=
\epsilon(G^{[1;\xi,\eta]})(\eta)/U'_{\xi,\eta}(\eta)^{\delta(G)}.
\end{equation}
The expression $U'_{\pi,\xi}(\xi)\frac{\omega}{d\pi}$ defines an element of $k_\ab$ independent of $\pi$  by a further application of the Chain Rule. Thus $r$ is independent of $\pi$.
By L'Hopital's Rule in the form
$$\delta((X-Y)^{[1;\xi,\eta]})=1,\;\;\;\epsilon((X-Y)^{[1;\xi,\eta]})=U'_{\xi,\eta}(Y)$$
we have
\begin{equation}\label{equation:ChainRuleEpsilon}
\epsilon(G^{[1;\xi,\eta]})
=\epsilon(G)(U_{\xi,\eta}(Y))\cdot U'_{\xi,\eta}(Y)^{\delta(G)},
\end{equation}
and hence
the right side of (\ref{equation:DiagonalExpressionBis}) is independent of $\eta$. Thus $r$ is independent of $\eta$. We therefore have simply
\begin{equation}\label{equation:DiagonalExpressionTer}
\left(U'_{\pi,\xi}(\xi)\frac{\omega}{d\pi}\right)^{\delta(G)}r^{\norm{b}}=
\epsilon(G)(\xi).
\end{equation}
By definition of interpolating gadget  we have
$$\epsilon(F^{[ab,ac;\xi]})=
\epsilon(G^{[a;\xi]})=\epsilon(G^{\norm{a}})=
\epsilon(G)^{\norm{a}}$$
for all $a\in k^\times$
such that $\norm{a}\geq 1$.
Thus $r$ depends only on the ratio
$c/b$.

It remains only to show that $r^{\norm{b}}\in k^\times$. In any case we have $r^{\norm{b}}\in k^\times_\ab$ by (\ref{equation:DiagonalExpressionTer}).  Fix $\sigma\in \Gal(k_\ab/k)$ arbitrarily.  It
will be enough to show that $\sigma r^{\norm{b}}=r^{\norm{b}}$.
We have
$$\epsilon(G)(\xi)=
\epsilon((\sigma G)^{[1;\sigma \xi,\xi]})(\xi)=
\epsilon(\sigma G)(\sigma \xi)\cdot U'_{\sigma \xi,\xi}(\xi)^{\delta(G)}
=\sigma (\epsilon(G)(\xi))\cdot U'_{\sigma \xi,\xi}(\xi)^{\delta(G)}$$
by the fourth characterization of an interpolating gadget given in Proposition~\ref{Subproposition:InterpolatingGadget} followed by an application of  (\ref{equation:ChainRuleEpsilon}).
We have 
$$\sigma(U'_{\pi,\xi}(\xi))=U'_{\pi,\sigma\xi}(\sigma\xi)=
U_{\pi,\sigma\xi}'(U_{\sigma\xi,\xi}(\xi))=U_{\pi,\xi}'(\xi)/U_{\sigma\xi,\xi}'(\xi)$$
by  the Chain Rule.
Finally, by (\ref{equation:DiagonalExpressionTer}), we indeed have $\sigma r^{\norm{b}}=r^{\norm{b}}$. \qed
\subsubsection{Definitions}
For every interpolating gadget $(F,\xi)$
and $a\in k^\times$ there exist by the preceding proposition unique 
$$\Theta_\sh(a,F,\xi)\in \ZZ[1/p],\;\;\;\left(\begin{array}{c}
a\\
F,\xi
\end{array}\right)_{\sh}\in k_\perf^\times$$
such that the following statement holds: for all $v,w\in k^\times$, $\eta\in \PPP^*$, uniformizers $\pi\in \OO$, 
and integers $\ell$ such that
$$a=w/v,\;\;\;\min(\norm{v},\norm{w})\geq 1,\;\;\;\xi\leq \eta,$$
$$\ell=(\textup{\mbox{$(X-Y)$-valuation of $F^{[v,w;\xi,\eta]}$}}),
$$
we have
\begin{equation}\label{equation:ShadowThetaDef}
\norm{v}\Theta_\sh(a,F,\xi)=\ell,
\end{equation}
\begin{equation}\label{equation:ShadowCatalanDef}
\left(\begin{array}{c}
a\\
F,\xi
\end{array}\right)_\sh^{\norm{v}}=
\left(\frac{d\pi}{\omega}\right)^{\ell}\cdot
\left(\frac{1}{U'_{\pi,\eta}(\eta)}\right)^{\ell}\cdot \left.
\frac{F^{[v,w;\xi,\eta]}(X,Y)}{(X-Y)^{\ell}}\right|_{X=Y=\eta}.
\end{equation}
We call $\Theta_\sh(\cdot,\cdot,\cdot)$ the {\em shadow theta symbol} and $\left(\begin{subarray}{c}\cdot\\
\cdot,\cdot
\end{subarray}\right)_\sh$ the {\em shadow Catalan symbol}. 
The somewhat complicated definitions of these symbols
have been contrived  to trivialize the verification of their formal properties.
\subsubsection{Formal properties of the shadow theta symbol}
Let $(F,\xi)$ and $(G,\xi)$ be interpolating gadgets
with a common basepoint and fix $a\in k^\times$ arbitrarily.
We have
 $$\Theta_\sh(a,FG,\xi)=
 \Theta_\sh(a,F,\xi)+\Theta_\sh(a,G,\xi).$$
 We have a functional equation
 \begin{equation}\label{equation:ShadowThetaFE}
\Theta_{\sh}(a,F,\xi)=\norm{a}\Theta_{\sh}(a^{-1},F^\dagger,\xi).
\end{equation}
This should be compared with functional
equation (\ref{equation:NoughtThetaFE}) for the theta symbol.
Given also $b,c\in k^\times$ such that $\min(\norm{b},\norm{c})\geq 1$, we have a scaling rule
\begin{equation}
\label{equation:ShadowThetaScaling}
\Theta_{\sh}(a,F^{[b,c;\xi,\eta]},\eta)=\norm{b}\Theta_{\sh}(ac/b,F,\xi).\end{equation}
 This should be compared with scaling rule (\ref{equation:ThetaScaling})
for the theta symbol.
We have
\begin{equation}\label{equation:ScholiumBis}
(\mbox{$Z_{t,\xi}$-valuation of $F$})=\Theta_\sh(t,F,\xi)/\max(\norm{t},1)
\end{equation}
for all $t\in k^\times$.
In principle (\ref{equation:ScholiumBis}) could have been taken as the 
definition  of the shadow theta symbol.
The function 
$$(a\mapsto\Theta_\sh(a,F,\xi)):k^\times\rightarrow \ZZ[1/p]$$ 
is compactly supported and factors through the quotient  
$$k^\times/\{a\in \OO^\times \mid \rho(a)\xi=\xi\}.$$
This should be compared with properties (\ref{equation:ThetaAsymptotics}) and (\ref{equation:ThetaLocalConstancy})
of the theta symbol. The deductions of these formal properties 
of the shadow theta symbol from the definitions are straightforward. We omit the details.

\subsubsection{Formal properties of the shadow Catalan symbol}
Again, let $(F,\xi)$ and $(G,\xi)$ be interpolating gadgets
with a common basepoint and fix $a\in k^\times$ arbitrarily.
We have
 $$\left(\begin{array}{c}a\\
 FG,\xi\end{array}\right)_\sh=\left(\begin{array}{c}a\\
 F,\xi\end{array}\right)_\sh\left(\begin{array}{c}a\\
 G,\xi\end{array}\right)_\sh.$$
We have a functional equation 
\begin{equation}\label{equation:ShadowCatalanFunctionalEquation}
\left(\begin{array}{c}
a^{-1}\\
F^\dagger,\xi
\end{array}\right)_\sh^{\norm{a}}=(-1)^{\Theta_\sh(a,F,\xi)}\left(\begin{array}{c}
a\\
F,\xi
\end{array}\right)_\sh.
\end{equation}
This should be compared with functional equation (\ref{equation:CatalanFunctionalEquation}) 
for the Catalan symbol.
Given also $b,c\in k^\times$ such that $\min(\norm{b},\norm{c})\geq 1$ and $\eta\in\PPP^*$
such that $\xi\leq \eta$ we have a scaling rule
\begin{equation}\label{equation:ShadowCatalanScaling}
\left(\begin{array}{c}
a\\
F^{[b,c;\xi,\eta]},\eta
\end{array}\right)_\sh
=\left(\begin{array}{c}
ac/b\\
F,\xi
\end{array}\right)_\sh^{\norm{b}}.\end{equation}
This should be compared with scaling rule
(\ref{equation:PartialCatalanScaling}) for the Catalan symbol. 
We have
\begin{equation}\label{equation:ShadowCatalanDefBis}
\left(\begin{array}{c}
a\\
F,\xi
\end{array}\right)_\sh=F(\xi,(\rho(a)^{-1}\xi)^{\norm{a}})\;\;\mbox{for $\max(\norm{a},\norm{a}^{-1})\gg 0$.}
\end{equation}
(More precisely, the equality above holds for $a\in k^\times$
if and only if $\Theta_\sh(a,F,\xi)=0$.)
 The function
$$\left(a\mapsto \left(\begin{array}{c}
a\\
F,\xi
\end{array}\right)_\sh\right):k^\times \rightarrow
k_\perf^\times$$
factors through the quotient
$k^\times/\{a\in \OO^\times \mid \rho(a)\xi=\xi\}$.
This should be compared with property (\ref{equation:CatalanLocalConstancy})
of the Catalan symbol. The deductions of these formal properties 
of the shadow Catalan symbol from the definitions are straightforward. We omit the details.

\begin{Theorem}[Shadow Stirling formula]\label{Theorem:ShadowRiggedStirling}
$\;$

Fix an interpolating gadget $(F,\xi)$.\\

(i) For all $a\in k^\times$ we have
\begin{equation}\label{equation:ShadowRiggedStirling}
\begin{array}{cl}
&\displaystyle\Theta_\sh(a,F,\xi)\ord \omega+\ord \left(\begin{array}{c}
a\\
F,\xi
\end{array}\right)_\sh\\\\
\geq&\displaystyle\int \frac{\Theta_\sh(at,F,\xi)-\one_{\OO^\times}(t)\Theta_\sh(a,F,\xi)}{\norm{1-t}}d\mu^\times(t)\\\\
&+
(\textup{\mbox{$X$-valuation of $F$}})\cdot\ord \xi\\\\
& + \norm{a}\cdot (\textup{\mbox{$Y$-valuation of $F$}})\cdot\ord \xi.
\end{array}
\end{equation}\\

(ii) The following statements are equivalent:\\

\begin{itemize}
\item $(F,\xi)$ is of Coleman type.
\item Equality holds in (\ref{equation:ShadowRiggedStirling}) for all $a\in k^\times$.
\item Equality holds in (\ref{equation:ShadowRiggedStirling}) for at least one $a\in k^\times$.\\
\end{itemize}

(iii) Let $\delta(a)$ denote the left side of (\ref{equation:ShadowRiggedStirling}) minus the right side of  (\ref{equation:ShadowRiggedStirling}).
Then we have 
$$\limsup_{\norm{a}\rightarrow\infty} \delta(a)<\infty,\;\;\;\limsup_{\norm{a}\rightarrow 0} \delta(a)/\norm{a}<\infty.$$
\end{Theorem}
\noindent This theorem is parallel to the simplified
statement (\ref{equation:SimplifiedRiggedStirling}) of the Stirling formula for rigged virtual lattices. 
\proof We turn first to the proof of part (i) of the theorem, which is most of the work. Choose $b,c\in k^\times$ such that
$\min(\norm{b},\norm{c})\geq 1$ and $a=c/b$. 
 Also choose a uniformizer $\pi\in \OO$ and $\eta\in \Xi_\pi$ such that $\xi\leq \eta$.
Put 
$$\ell =(\mbox{$(X-Y)$-valuation of $F^{[b,c;\xi,\eta]}$}),$$
$$\tilde{F}=\tilde{F}(X,Y)=F^{[b,c;\xi,\eta]}(X,Y)/(X-Y)^\ell,$$
$$G(X,Y)=((X-Y)^{[1;\pi,\eta]})^\ell,\;\;\;\tilde{G}(X,Y)=G(X,Y)/(X-Y)^\ell,$$
noting that
$$\tilde{G}(\eta,\eta)=U'_{\pi,\eta}(\eta)^\ell,$$
and also that
$$\ell=\norm{b}\Theta_\sh(a,F,\xi)$$
by definition of the shadow theta symbol.
By definition of the shadow Catalan symbol,
$$\ord\tilde{F}(\eta,\eta)-\ord\tilde{G}(\eta,\eta)$$
 equals the lefthand side of  (\ref{equation:ShadowRiggedStirling})
 multiplied by $\norm{b}$.

 We turn to the investigation of the righthand side of  (\ref{equation:ShadowRiggedStirling}).
 Choose a set of representatives $T_\xi\subset k^\times$ for the quotient 
$$k^\times/\{a\in \OO^\times\mid\rho(a)\xi=\xi\}.$$
By  (\ref{equation:ScholiumBis}) we can write
\begin{equation}\label{equation:PrelimFFac}
\begin{array}{rcl}
F&\sim &W \cdot X^{(\mbox{\scriptsize $X$-valuation of $F$})} \cdot Y^{(\mbox{\scriptsize $Y$-valuation of $F$})}\\\\
&&\displaystyle \cdot \prod_{t\in T_\xi} Z_{t,\xi}^{\Theta_\sh(t,F,\xi)/\max(\norm{t},1)}
\end{array}
\end{equation}
where $W\in \Lambda_0$ is a power series divisible neither by $X$, nor by $Y$, nor by $Z_{t,\xi}$ for any $t\in T_\xi$. Now let $n$ be the smallest positive integer such that
$[\pi^n]_\pi(\eta)=0$. 
By the explicit reciprocity law (\ref{equation:ExplicitReciprocityLaw}) of Lubin-Tate,
$$\{a\in \OO^\times \mid \rho(a)\eta=\eta\}=1+\pi^n\OO.$$
Choose a set $T_\eta$ of representatives for the quotient 
$k^\times/(1+\pi^n\OO)$, making this choice so that $1\in T_\eta$.
We have
\begin{equation}\label{equation:InfamousW}
\begin{array}{rcl}
\tilde{F}&\sim& W^{[b,c;\xi,\eta]}\cdot X^{(\mbox{\scriptsize $X$-valuation of $F$})\norm{b}\frac{\ord\xi}{\ord \eta}}\\\\
&& \cdot Y^{(\mbox{\scriptsize $Y$-valuation of $F$})\norm{c}\frac{\ord\xi}{\ord \eta}}\\\\
&&\displaystyle\cdot \prod_{1\neq t\in T_\eta}Z_{t,\eta}^{\norm{b}\Theta_\sh(ta,F,\xi)/\max(\norm{t},1)}
\end{array}
\end{equation}
by  
(\ref{equation:KeyZFac}), (\ref{equation:XYFac}), and (\ref{equation:PrelimFFac}) after a calculation which we omit. Note that $W^{[b,c;\xi,\eta]}$ is divisible by neither $X$, nor $Y$, nor by $Z_{t,\eta}$ for any $t\in k^\times$. We have also
$$\tilde{G}\sim\prod_{1\neq t\in T_\eta\cap \OO^\times}Z_{t,\eta}^{\norm{b}\Theta_\sh(a,F,\xi)}$$
by (\ref{equation:KeyZFac}) and the definition of the shadow theta symbol. For $1\neq t\in T_\eta$ we have via the explicit reciprocity law (\ref{equation:ExplicitReciprocityLaw}) and the definition of $Z_{t,\xi}$ that
\begin{equation}\label{equation:ShadowClincher}
\frac{\ord Z_{t,\eta}(\eta,\eta)}{\max(\norm{t},1)}=
\frac{\mu^\times(1+\pi^n\OO)}{\norm{1-t}}.
\end{equation}
A straightforward calculation now shows that
$$\ord \tilde{F}(\eta,\eta)-\ord \tilde{G}(\eta,\eta)-\ord W((\rho(b)^{-1}\xi)^{\norm{b}},(\rho(c)^{-1}\xi)^{\norm{c}})$$
equals the right side of (\ref{equation:ShadowRiggedStirling})
multiplied by $\norm{b}$. Since 
$$\ord W((\rho(b)^{-1}\xi)^{\norm{b}},(\rho(c)^{-1}\xi)^{\norm{c}})\geq 0,$$
the proof of part (i) of the theorem is complete.

We turn to the proof of part (ii) of the theorem.
Now $(F,\xi)$ is of Coleman type if and only if $W\sim 1$
only if 
equality holds in (\ref{equation:ShadowRiggedStirling}) for all  $a\in k^\times$. Thus the first of the three given statements implies the second. Of course the second statement trivially implies the third.
Clearly, $(F,\xi)$ is not of Coleman type if and only if
$W$ is a nonunit of $\Lambda_0$ if and only if inequality
holds in (\ref{equation:ShadowRiggedStirling}) for all $a\in k^\times$.
Thus the negation of the first statement implies the negation of the third statement.
The proof of part (ii) of the theorem is complete.

We turn finally to the proof of part (iii) of the theorem.
By the Weierstrass Division Theorem applied
in $\FF_{q^m}[[X,Y]]$ for some $m$, we may assume that $W$ is a monic polynomial in $X$ of degree $d$ and that $W(X,0)=X^d$. 
Then we have
$$\norm{b}^{-1}\ord W((\rho(b)^{-1}\xi)^{\norm{b}},(\rho(c)^{-1}\xi)^{\norm{c}})=d\ord\xi$$
for $\norm{c/b}=\norm{a}\gg 1$, which proves 
$\limsup_{\norm{a}\rightarrow\infty}\delta(a)<\infty$.
The finiteness of
$\limsup_{\norm{a}\rightarrow 0}\delta(a)/\norm{a}$
is proved similarly---one instead reduces to the case in which $W$ is a monic polynomial in $Y$ of degree $e$ such that $W(0,Y)=Y^e$.  The proof of part (iii) of the theorem is complete. The proof of the theorem is complete.
\qed

\section{Interpolability and related notions}
We now bring the two theories of symbols into alignment. We define the notion of interpolability
along with several variants. 
 We work out key consequences of these definitions.

\subsection{Definitions}
\subsubsection{Asymptotic interpolability}
Given  a proper rational 
rigged virtual lattice $\Phi$ and a pair $(F,\xi)$
with $0\neq F\in \Lambda$ and $\xi\in \PPP^*$,
if
$$\left(\begin{array}{c}
a\\
\Phi
\end{array}\right)_+=F(\xi,(\rho(a)^{-1}\xi)^{\norm{a}})$$
for all $a\in k^\times$ such that $\norm{a}\gg 1$,
we say that $\Phi$ is {\em asymptotically yoked} to  $(F,\xi)$.  
Note that $(F,\xi)$ is in this situation automatically an interpolating gadget by Proposition~\ref{Subproposition:InterpolatingGadget}.
We say that a rigged virtual lattice $\Phi$ is {\em asymptotically interpolable}
if $\Phi$ is rational, proper and asymptotically yoked to some interpolating gadget.
\subsubsection{Remark}
The notion of asymptotic interpolability
defined above seems to differ from that
defined in the introduction. 
Let us now reconcile the definitions by proving asymptotic interpolability in the two senses to be equivalent. In any case, it is clear that 
given any proper rational virtual rigged lattice $\Phi$,
asymptotic interpolability of $\Phi_*$ in the sense defined in the introduction
implies asymptotic interpolability of $\Phi$ in the sense defined immediately above. Conversely,
suppose that $\Phi$ is asymptotically yoked to $(F,\xi)$ in the sense defined immediately above. The problem here is that we might not have $F\in \FF_q[[X,Y]][X^{-1},Y^{-1}]$.
Fix a uniformizer $\pi\in \OO$ and then choose $\eta\in \Xi_\pi$ such that $\xi\leq \eta$.
Put $G=F^{[1;\xi,\eta]}$.  By Proposition~\ref{Subproposition:InterpolatingGadget} both $(F,\xi)$ and $(G,\eta)$ are interpolating gadgets.
By Proposition~\ref{Subproposition:TotallyRamified},
we have $G\in \FF_q[[X,Y]][X^{-1},Y^{-1}]$.
One verifies easily that
$(G,\eta)$ is asymptotically yoked to $\Phi_*$
in the sense defined in the introduction.
The proof of equivalence is complete.

\subsubsection{Strict interpolability}
\noindent Given a proper rational rigged virtual lattice $\Phi$
and  an interpolating gadget $(F,\xi)$ such that
\begin{equation}\label{equation:Yoke}
\Theta(a,\Phi)=\Theta_{\sh}(a,F,\xi),\;\;\;
\left(\begin{array}{c}
a\\
\Phi
\end{array}\right)=
\left(\begin{array}{c}
a\\
F,\xi
\end{array}\right)_\sh
\end{equation} for all $a\in k^\times$,
we say that $\Phi$ and $(F,\xi)$ are {\em yoked}.
We say that a rigged virtual lattice $\Phi$ is {\em strictly interpolable}
if $\Phi$ is proper, rational, and yoked to some interpolating gadget. A strictly interpolable rigged virtual lattice is necessarily effective,
and moreover necessarily asymptotically interpolable
by (\ref{equation:DownAsymptoticPartialCatalanBis}) 
and (\ref{equation:ShadowCatalanDefBis}).

\subsubsection{Dependence of the yoking relation on $\omega$ and $\FF$}\label{subsubsection:YokeIndependence}
The relation ``$\Phi$ and $(F,\xi)$ are yoked'' is invariant under change of the Chevalley differential $\omega$  as one sees by comparing
the definition (\ref{equation:ShadowCatalanDef}) of the shadow Catalan symbol
to formula (\ref{equation:ChangeOfDifferential}) describing the effect of change of $\omega$ on the Catalan symbol.
In other words, $\omega$ ``cancels'' when one forms the ratio $\left(\begin{subarray}{c}
a\\
\Phi
\end{subarray}\right)/\left(\begin{subarray}{c}
a\\
F,\xi
\end{subarray}\right)_\sh$. The yoking relation also remains invariant under
change of subfield $\FF\subset \FF_q$ since neither the Catalan symbol nor its shadow depend on $\FF$. 
Thus the yoking relation is  intrinsic. Similarly, but much more trivially,
the asymptotic yoking relation is intrinsic since neither $\omega$
nor $\FF$ have anything to do with that definition.

\subsubsection{Interpolability}
Given $\Phi\in \rigged(k)$, we say that
$\Phi$ is
{\em interpolable} if $$\Phi=\sum_{i=1}^N \alpha_i \Phi_i^{(b_i,c_i)}$$
for some integer $N\geq 0$,
numbers $\alpha_i\in \ZZ[1/p]$, numbers $b_i,c_i\in k^\times$, and strictly interpolable rigged virtual lattices $\Phi_i$.  Given also for each $i$
the basepoint $\xi_i$ of some interpolating gadget yoked to $\Phi_i$,
and some common upper bound $\xi\in \PPP^*$
for the $\xi_i$, we say that $\Phi$ is  {\em of conductor $\leq \xi$}.
 An interpolable rigged virtual lattice is automatically proper and rational.
By definition the space of interpolable rigged virtual lattices is a $\ZZ[1/p]$-module, and also 
for every $b,c\in k^\times$ stable under the
operation $\Phi\mapsto \Phi^{(b,c)}$.

\begin{Proposition}\label{Proposition:InterpolatingUniqueness}
Let $(F,\xi)$ and $(G,\eta)$ be interpolating gadgets asymptotically
yoked to the same proper rational rigged virtual lattice.
  Let $\zeta\in \PPP^*$ be a common upper bound for $\xi$ and $\eta$. 
Then we have $F^{[1;\xi,\zeta]}=G^{[1;\eta,\zeta]}$.
\end{Proposition}
\proof Put
$\tilde{F}=F^{[1;\xi,\zeta]}$ and $\tilde{G}=G^{[1;\eta,\zeta]}$.
For every $a\in k^\times$ such that $\norm{a}\gg 1$ we have
$$
\begin{array}{rcl}
\tilde{F}(\zeta,(\rho(a)^{-1}\zeta)^{\norm{a}})&=&
F(\xi,(\rho(a)^{-1}\xi)^{\norm{a}})\\
&=&
G(\eta,(\rho(a)^{-1}\eta)^{\norm{a}})=
\tilde{G}(\zeta,(\rho(a)^{-1}\zeta)^{\norm{a}}).
\end{array}$$
Thus the difference  $\tilde{F}-\tilde{G}$ belongs to infinitely many maximal ideals of the principal ideal domain $\Lambda$ and hence vanishes identically.
\qed

\begin{Proposition}\label{Proposition:DaggerInterpretation}
Let $\Phi$ be a strictly interpolable rigged virtual lattice yoked to an interpolating gadget $(F,\xi)$.
Then $\GGG_0[\Phi]$ is strictly interpolable and yoked to $(F^\dagger,\xi)$. 
\end{Proposition}
\proof Fix $a\in k^\times$ arbitrarily.
We have
$$\norm{a}\Theta(a^{-1},\GGG_0[\Phi])=\Theta(a,\Phi)=\Theta_\sh(a,F,\xi)=
\norm{a}\Theta_\sh(a^{-1},F^\dagger,\xi)$$
where the first and third steps are justified by the functional equations
(\ref{equation:NoughtThetaFE}) and (\ref{equation:ShadowThetaFE})
satisfied by the theta symbol and its shadow, respectively.
We have
$$\begin{array}{rcl}
\displaystyle\left(\begin{array}{c}
a^{-1}\\
\GGG_0[\Phi]
\end{array}\right)^{\norm{a}}&=&\displaystyle
(-1)^{\Theta(a,\Phi)}\left(\begin{array}{c}
a\\
\Phi
\end{array}\right)\\&=&(-1)^{\Theta_\sh(a,F,\xi)}\left(\begin{array}{c}
a\\
F,\xi
\end{array}\right)_\sh=\left(\begin{array}{c}
a^{-1}\\
F^\dagger,\xi
\end{array}\right)^{\norm{a}},
\end{array}$$
where the first and third steps are justified by the functional equations (\ref{equation:CatalanFunctionalEquation})
and (\ref{equation:ShadowCatalanFunctionalEquation}) satisfied by the Catalan symbol and its shadow,
respectively. \qed

\begin{Proposition}\label{Proposition:SquareBracketInterpretation}
Let $\Phi$ be a strictly interpolable rigged virtual lattice yoked to an interpolating gadget $(F,\xi)$.
Let $b,c\in k^\times$ be given such that $\min(\norm{b},\norm{c})\geq 1$.
Let $\eta\in \PPP^*$ be given such that $\xi\leq \eta$. Then $\Phi^{(b,c)}$ is strictly interpolable and yoked to $(F^{[b,c;\xi,\eta]},\eta)$.
\end{Proposition}
\proof Fix $a\in k^\times$ arbitrarily.
We have
$$\Theta(a,\Phi^{(b,c)})=
\norm{b}\Theta(ac/b,\Phi)=
\norm{b}\Theta(ac/b,F,\xi)_\sh=
\Theta(a,F^{[b,c;\xi,\eta]}),$$
where the first and third steps are justified by the scaling rules 
(\ref{equation:ThetaScaling}) and (\ref{equation:ShadowThetaScaling})
for the theta symbol and its shadow, respectively.
We have
$$\left(\begin{array}{c}
a\\
\Phi^{(b,c)}
\end{array}\right)=
\left(\begin{array}{c}
ac/b\\
\Phi
\end{array}\right)_\sh^\norm{b}=
\left(\begin{array}{c}
ac/b\\
F,\xi
\end{array}\right)_\sh^{\norm{b}}
=
\left(\begin{array}{c}
a\\
F^{[b,c;\xi,\eta]},\eta
\end{array}\right)_\sh$$
where the first and third steps are justified by the scaling rules 
(\ref{equation:CatalanScaling}) and (\ref{equation:ShadowCatalanScaling})
for the Catalan symbol and its shadow, respectively. 
\qed

\begin{Proposition}\label{Proposition:CheapExamples}
Let $\Phi$ be a rational rigged virtual lattice such that $\Phi_*(0,0)$ and $\GGG_0[\Phi]_*(0,0)$ are integers. Fix a uniformizer $\pi\in \OO$. Fix    $b\in k^\times$.
Let $B(X)\in \FF_q[[X]][X^{-1}]$ be the unique Laurent series such that $B(\pi)=b$.
 Put
$$F(X,Y)=B(X)^{-\Phi_*(0,0)}B(Y)^{\GGG_0[\Phi]_*(0,0)}
\in\Lambda^\times.$$
Then: (i) $(F,\pi)$ is an interpolating gadget.
(ii) $\Phi^{(b,b)}/\norm{b}-\Phi$ is strictly interpolable
and yoked to $(F,\pi)$.
\end{Proposition}
\proof  (i) For all $G\in \FF_q[[X,Y]][X^{-1},Y^{-1}]$
and $d,e\in k^\times$ such that  \linebreak $\min(\norm{d},\norm{e})\geq 1$ we have
$G^{[d,e;\pi]}(X,Y)=G(X^{\norm{d}},Y^{\norm{e}})$, and hence (consider the case $d=e$ in the preceding) the
pair $(G,\pi)$ is an interpolating gadget.
In particular,  $(F,\pi)$ is an interpolating gadget. 
(ii) Fix $a\in k^\times$ arbitrarily.
By scaling rule (\ref{equation:ThetaScaling}) we have 
$$\Theta(a,\Phi^{(b,b)}/\norm{b}-\Phi)=0.$$
By scaling rule (\ref{equation:CatalanScaling}) we have
$$\left(\begin{array}{c}
a\\
\Phi^{(b,b)}/\norm{b}-\Phi
\end{array}\right)=b^{\norm{a}\GGG_0[\Phi]_*(0,0)-\Phi_*(0,0)}.$$
Since $F$ is a unit of $\Lambda$,
$$\Theta_\sh(a,F,\xi)=0$$
and
$$
\left(\begin{array}{c}
a\\
F,\pi
\end{array}\right)_\sh=F^{[1,a;\xi]}\vert_{X=Y=\pi}=F(\pi,\pi^{\norm{a}})=
b^{\norm{a}\GGG_0[\Phi]_*(0,0)-\Phi_*(0,0)}.$$
Thus we have strict interpolation as claimed.
\qed

 \begin{Proposition}\label{Proposition:Underscore}
Let $\Phi$ be a strictly interpolable rigged virtual lattice
yoked to an interpolating gadget $(F,\xi)$. Then $(F,\xi)$ is of Coleman type,
and moreover
\begin{equation}\label{equation:SinhaProblem}
\begin{array}{rcl}
(\textup{\mbox{$X$-valuation of $F$}})\cdot \ord \xi&=&\displaystyle
\int \Phi_*(0,t)d\mu^\times(t),\\\\
(\textup{\mbox{$Y$-valuation of $F$}})\cdot \ord\xi&=&
\displaystyle
\int \GGG_0[\Phi]_*(0,t)d\mu^\times(t).
\end{array}
\end{equation}
\end{Proposition}
\proof By our hypotheses combined with Theorem~\ref{Theorem:RiggedStirling} in its simplified form (\ref{equation:SimplifiedRiggedStirling}) and part (i) of Theorem~\ref{Theorem:ShadowRiggedStirling}, we have
$$\begin{array}{rcl}
\delta(a)&=&\displaystyle\left(\int \Phi_*(0,t)d\mu^\times(t)-(\textup{\mbox{$X$-valuation of $F$}})\cdot \ord \xi\right)\\\\
&&\displaystyle+\norm{a}\left(\int \GGG_0[\Phi]_*(0,t)d\mu^\times(t)-
(\textup{\mbox{$Y$-valuation of $F$}})\cdot \ord\xi\right)\geq 0.
\end{array}$$
By part (iii) of Theorem~\ref{Theorem:ShadowRiggedStirling} this is possible only if $\delta(a)$
vanishes identically, whence (\ref{equation:SinhaProblem}). By part (ii) of Theorem~\ref{Theorem:ShadowRiggedStirling}, since $\delta(a)$ vanishes identically, $(F,\xi)$ must be of Coleman type. 
\qed
\begin{Proposition}\label{Proposition:Descent}
Let $\Phi$ be a rigged virtual lattice. Let $r$ be a power of the characteristic $p$ of $k$. Assume that 
$r\Phi$ is strictly interpolable and that 
\begin{equation}\label{equation:PreSeparability}
\left(\begin{array}{c}
a\\
\Phi
\end{array}\right)^{\max(1,\norm{a}^{-1})}\in k^\times\;\;\mbox{for $\max(\norm{a},\norm{a}^{-1})\gg 1$.}
\end{equation}
Then $\Phi$ is strictly interpolable.
\end{Proposition}
\noindent Recall that  if $\Phi$ is separable, then (\ref{equation:PreSeparability})
is  satisfied,  by (\ref{equation:SeparabilityImplication}).
\proof Let $r\Phi$ be yoked to $(F,\xi)$. It is enough to show that  $F^{1/r}\in \Lambda$, because once that is shown, it is easy to verify that $(F^{1/r},\xi)$ is an interpolating gadget  yoked to $\Phi$. By formula (\ref{equation:ShadowCatalanDefBis}) and hypothesis we have
$$
F(\xi,(\rho(a)^{-1}\xi)^{\norm{a}})^{\max(1,\norm{a}^{-1})}\in (k^\times)^r\;\;
\mbox{for $\max(\norm{a},\norm{a}^{-1})\gg 1$.}
$$
  Choose $a_0\in k^\times$
 such that 
 $$\norm{a_0}\geq r,\;\;\;\;\rho(a_0)^{-1}\xi=\xi,\;\;\;
 \rho(a_0)^{-1}F=F.$$
 Then we have
 \begin{equation}\label{equation:DescentHyp}
F(\xi^{\norm{a_0}^i},\xi),F(\xi,\xi^{\norm{a_0}^i})\in (k^\times)^r\;\;\mbox{for $i\gg 0$.}
\end{equation}
   Write 
 $$F(X,Y)=\sum_{j=0}^{r-1} F_j(X^r,Y)X^j \;\;\;(F_j(X,Y)\in \Lambda)$$
 in the unique possible way.  Let $K$ be the field generated over $k$ by $\FF_{q^\infty}$ and $\xi$.  For some integer $i_0\geq 0$ and all integers $i\geq i_0$ we have
$$\sum_{j=0}^{r-1}F_j(\xi^r,\xi^{\norm{a_0}^i})\xi^j \in (k^\times)^r\subset K^r$$
by (\ref{equation:DescentHyp}).
But the powers $\{\xi^j\}_{j=0}^{r-1}$ form a basis of $K$ as a vector space
over $K^r$, and $F_j(\xi^r,\xi^{\norm{a_0}^i})\in K^r$ for all $i\geq 0$ and $j=0,\dots,r-1$.
We therefore have $F_j(\xi^r,\xi^{\norm{a_0}^i})=0$
for all $i\geq i_0$ and $j=1,\dots,r$. But then, for $j=1,\dots,r$,
since $F_j(X,Y)$ is  contained in infinitely many distinct maximal ideals of the principal ideal domain $\Lambda$, necessarily $F_j(X,Y)$ vanishes identically.
Thus we have 
$$F(X^{1/r},Y)=F_0(X,Y)\in \Lambda.$$ 
Symmetrically, we have 
$$F(X,Y^{1/r})\in \Lambda.$$
Finally, we have $F(X^{1/r},Y^{1/r})\in \Lambda$
and hence $F^{1/r}\in \Lambda$.
\qed
\begin{Proposition}\label{Proposition:InterpolabilitySimplification}
Fix a basepoint $\xi\in \PPP^*$. 
Let $\Phi$ be an interpolable rigged virtual lattice
of conductor $\leq \xi$. Then there exist strictly
interpolable rigged virtual gadgets $\Phi_{\pm}$
each of which is yoked to an interpolating gadget with basepoint $\xi$
such that $r\Phi=\Phi_+-\Phi_-$ for some power $r$ of the characteristic $p$ of $k$.
\end{Proposition}
\proof By hypothesis
$$\Phi=\sum_{i=1}^N \alpha_i\Phi_i^{(b_i,c_i)}$$
where $\alpha_i\in \ZZ[1/p]$, $b_i,c_i\in k^\times$, $\Phi_i$ is strictly interpolable, and $\Phi_i$ is yoked to some interpolating gadget $(F_i,\xi_i)$
such that $\xi_i\leq \xi$.
For each index $i$ choose $d_i\in k^\times$ such that
$\min(\norm{d_ib_i},\norm{d_ic_i})\geq 1$.
Then we have
$$\Phi=\sum_{i=1}^N\alpha_i\norm{d_i}^{-1}
\Phi^{(d_ib_i,d_ic_i)}+\sum_{i=1}^N \alpha_i(\Phi_i^{(b_i,c_i)}-\norm{d_i}^{-1}\Phi_i^{(d_ib_i,d_ic_i)}).
$$
By Propositions~\ref{Proposition:SquareBracketInterpretation} and \ref{Proposition:CheapExamples} we have now exhibited $\Phi$ as a finite $\ZZ[1/p]$-linear combination of strictly interpolable rigged virtual lattices each yoked to an interpolating gadget with basepoint bounded by $\xi$. After relabeling, we may simply assume that 
$$r\Phi=\sum_{i=1}^N \alpha_i\Phi_i$$
where $\alpha_i\in \ZZ$, $r$ is a power of characteristic $p$ of $k$, $\Phi_i$ is strictly interpolable, and $\Phi_i$ is yoked to an interpolating gadget $(F_i,\xi_i)$ such that $\xi_i\leq \xi$.
By Proposition~\ref{Proposition:SquareBracketInterpretation} we may assume that
$\xi_i=\xi$ for all $i$.  Finally, after grouping terms according to the signs of their coefficients,
we may assume that $N=2$, $\alpha_1=1$ and $\alpha_2=-1$,
in which case $\Phi_+=\Phi_1$ and $\Phi_-=\Phi_2$ have the desired properties.
 \qed
\begin{Proposition}\label{Proposition:StrictnessCriterion}
Let a rigged virtual lattice $\Phi$ and a basepoint $\xi\in \PPP^*$ be given.
If $\Phi$ is interpolable of conductor $\leq \xi$,
effective and separable,
then $\Phi$ is strictly interpolable
and yoked to an interpolating gadget with basepoint $\xi$.
\end{Proposition}
\proof By Proposition~\ref{Proposition:InterpolabilitySimplification} we may assume that
for some power $r$ of the characteristic $p$ of $k$
and strictly interpolable rigged virtual lattices
$\Phi_{\pm}$ yoked to interpolating gadgets $(F_{\pm},\xi)$ we have $r\Phi=\Phi_+-\Phi_-$.
Put $F=F_+/F_-$. Now $F$ {\em a priori} belongs to the fraction field of the ring $\Lambda$, but by Proposition~\ref{Proposition:Underscore} the power series
$F_{\pm}$ are  Coleman units with respect to $\xi$, and hence, since $\Phi$ is effective,
one can by means of formula (\ref{equation:ScholiumBis}) verify that $F_-$ divides $F_+$ and hence $F\in \Lambda$.
It is not difficult to verify that $(F,\xi)$ is an interpolating gadget yoked to $r\Phi$; these details we omit. Finally, $(F^{1/r},\xi)$ is an interpolating gadget yoked to $\Phi$ by Proposition~\ref{Proposition:Descent} and our hypothesis that $\Phi$ is separable. \qed

\section{Concrete examples of strict interpolation}
\begin{Theorem}\label{Theorem:InterpolationExamples}
Fix a uniformizer $\pi\in \OO$. 
Fix a sequence $\{\xi_i\}_{i=0}^\infty$ in $k_\ab$ satisfying the relations
$$\xi_0=\pi,\;\;\xi_{i-1}=
\left\{\begin{array}{rl}
-\xi_1^{q-1}&\mbox{if $i=1$}\\
\pi\xi_i+\xi_i^q&\mbox{if $i>1$}
\end{array}\right.
$$
for $i>0$.
Put 
$$L=\bigoplus_{i=1}^\infty \FF_q\cdot \pi^{-i}\subset k.$$
Fix $t\in k^\times$ and an integer $M\geq 0$ such that
$\norm{\pi^{-M}}\geq \norm{t}$. Let $\Psi_t$ be the unique primitive rigged virtual lattice such that
$$(\Psi_t)_*=\one_{t+L}\otimes \one_{t+\OO}.$$
The $\Psi_t-\Psi_1$ is strictly interpolable and can be yoked to an interpolating gadget with basepoint $\xi_M$.
\end{Theorem}
\noindent  The proof of the theorem, which is by explicit calculation and construction, takes up the rest of this section.  The theorem is going to play a vital role in the proof of our main result.
We point out that the sequence $\{\xi_i\}_{i=1}^\infty$ is a $\pi$-division tower for 
the Lubin-Tate formal group discussed in \S\ref{subsubsection:ExplicitReciprocity}.    For simplicity we assume that
\begin{equation}\label{equation:SpecialAssumption}
\FF=\FF_q,\;\;\;\omega=d\pi.
\end{equation}
As explained in \S\ref{subsubsection:YokeIndependence}, these assumptions entail no loss of generality.
The calculations undertaken here are similar to those
in \cite[\S10]{AndersonStirling}, but on the whole simpler since local rather than adelic. 

\subsection{Review of some determinant identities}

\subsubsection{The Moore identity}
Given $x_1,\dots,x_n\in k$,
we define the {\em Moore determinant}
$$\Moore(x_1,\dots,x_n)=
\left|\begin{array}{cccc}
x_1^{q^{n-1}}&x_1^{q^{n-2}}&\dots&x_1\\
\vdots&\vdots&&\vdots\\
x_n^{q^{n-1}}&x_n^{q^{n-2}}&\dots&x_n\\
\end{array}\right|=\det_{i,j=1}^n x_i^{q^{n-j}}\in k.$$
Note that $\Moore(x_1,\dots,x_n)$ depends $\FF_q$-linearly
on each of its arguments, and vanishes unless its arguments are $\FF_q$-linearly independent.
Assuming that $x_1,\dots,x_n$ are $\FF_q$-linearly independent, we have \begin{equation}\label{equation:MooreIdentity}
\Moore(x_1,\dots,x_n)=
\prod_{i=1}^n\prod_{v\in V_i}(x_i+v),
\end{equation}
where $V_i$ is the $\FF_q$-span
of $x_{i+1},\dots,x_n$. This is the well-known {\em Moore identity}.
\subsubsection{The Ore identity}
Given $\zeta\in\Omega$ and $x_1,\dots,x_n\in k$, put
$$\Ore(\zeta, x_1,\dots,x_{n})=\left|\begin{array}{cccc}
\Res(x_1\zeta)&x_1^{q^{n-2}}&\dots&x_1\\
\vdots&\vdots&&\vdots\\
\Res(x_n\zeta)&x_n^{q^{n-2}}&\dots&x_n\\
\end{array}\right|\in k.$$
We call $\Ore(\zeta, x_1,\dots,x_{n})$ the {\em Ore determinant} of $\zeta,x_1,\dots,x_n$. 
It is clear that $\Ore(\zeta,x_1,\dots,x_n)$
depends $\FF_q$-linearly on each of its arguments.
Suppose now that $x_1,\dots,x_n$ are $\FF_q$-linearly independent
and that  $\Res(x_i\zeta)\neq 0$ for some $i$. 
Let $V$ be the $\FF_q$-span of $x_1,\dots,x_n$.
From the Moore identity one easily deduces that
\begin{equation}\label{equation:OreIdentity}
\frac{\Ore(\zeta,x_1,\dots,x_n)}{\Moore(x_1,\dots,x_n)}
=\prod_{\begin{subarray}{c}
v\in V\\
\Res(v\zeta)=1
\end{subarray}}v^{-1}.
\end{equation}
The latter we call the {\em Ore identity}, because it makes available to us 
certain useful features of Ore-Elkies-Poonen duality \cite[\S4.14]{GossBook}. 
The definition of the Ore determinant here is just a slight modification 
of the definition made in \cite[\S2.8]{AndersonStirling}.
\subsubsection{Special cases of the Jacobi-Trudi identity}
Let $x_1,\dots,x_N$ be distinct elements of $k$.
Let $Z$ be a variable.
We have
\begin{equation}\label{equation:JacobiTrudi}
\begin{array}{cl}
&\Res_{Z=0}
\left(Z^{-n}
\prod_{i=1}^N(1-x_iZ)^{-1}\frac{dZ}{Z}\right)\\\\
=&(\mbox{$n^{th}$ complete symmetric function of $x_1,\dots,x_N$})\\\\
=&\left|\begin{array}{cccc}
x_1^{n+N-1}&\dots&x_N^{n+N-1}\\
x_1^{N-2}&\dots&x_N^{N-2}\\
\vdots &&\vdots\\
1&\dots&1
\end{array}\right|\bigg/\left|\begin{array}{cccc}
x_1^{N-1}&\dots&x_N^{N-1}\\
\vdots&&\vdots\\
1&\dots&1
\end{array}\right|,
\end{array}
\end{equation}
\begin{equation}\label{equation:JacobiTrudiBis}
\begin{array}{cl}
&\Res_{Z=0}
\left(Z^{-n}
\prod_{i=1}^N(1+x_iZ)\frac{dZ}{Z}\right)\\\\
=&(\mbox{$n^{th}$ elementary symmetric function of $x_1,\dots,x_N$})\\\\
=&\left|\begin{array}{cccc}
x_1^N&\dots&x_N^{N}\\
\vdots&&\vdots\\
\widehat{x_1^{N-n}}&\dots&\widehat{x_N^{N-n}}\\
\vdots&&\vdots\\
1&\dots&1
\end{array}\right|\bigg/\left|\begin{array}{cccc}
x_1^{N-1}&\dots&x_N^{N-1}\\
\vdots&&\vdots\\
1&\dots&1
\end{array}\right|,
\end{array}
\end{equation}
where the ``hat'' indicates omission. In the first identity $n$ may take any nonnegative value, whereas in the second we restrict to $0\leq n\leq N$. The expressions on the left sides are residues of meromorphic differentials on the $Z$-line over $k$; we will make more use of such expressions later in our calculations. The denominators on the right sides are Vandermonde determinants. 
These identities are special cases of the {\em Jacobi-Trudi identity}. Concerning the latter, see \cite[I,3]{Macdonald} for background and details.

\subsection{Calculation of symbols}
Fix $a\in k^\times$. We will calculate the symbols
$$
\Theta(a,\Psi_t-\Psi_1),\;\;\;\left(\begin{array}{c}
a\\
\Psi_t-\Psi_1
\end{array}\right),$$
expressing the Catalan symbol in terms of Moore and Ore determinants.
Along the way we verify that $\Psi_t-\Psi_1$ is proper.
\subsubsection{Calculation of the theta symbol}
Note that
$$k=L\oplus \OO.$$
For every $x\in k$ we write
$$x=\lfloor x\rfloor +\langle x\rangle\;\;\;(\lfloor x\rfloor \in L,\;\;
\langle x\rangle\in\OO)$$
in the unique possible way. We have
\begin{equation}\label{equation:ThetaSetCalc}
\begin{array}{cl}
&(t+L)\cap a(t+\OO)\neq \emptyset \\
\Leftrightarrow&
(a-1)t\in L+a\OO\\
\Leftrightarrow&\langle (a-1)t\rangle\in a\OO\\
\Leftrightarrow&\lfloor at \rfloor+\langle t\rangle \in (t+L)\cap a(t+\OO)\\
\Rightarrow &\lfloor at \rfloor+\langle t\rangle\neq 0.
\end{array}
\end{equation}
By scaling rule (\ref{equation:RiggedScaling}) we have
$$(\Psi_t^{(1,a)})_*=\one_{t+L}\otimes \one_{a(t+\OO)},$$
hence via (\ref{equation:ThetaSetCalc}) we have
\begin{equation}\label{equation:DiagonalPsiT}
\begin{array}{cl}
&(\Psi_t^{(1,a)}-\Psi_1^{(1,a)})_*(x,x)\\\\
=&\left\{\begin{array}{rl}
\left(\one_{\lfloor at\rfloor+\langle t\rangle+L\cap a\OO}-
\one_{1+L\cap a\OO}\right)(x)&\mbox{if $\norm{a}\geq 1$,}\\
\one_{\{\lfloor at\rfloor+\langle t\rangle\}}(x)&\mbox{if
$\norm{a}<1$ and $\langle (a-1)t\rangle\in a \OO$,}\\
0&\mbox{if $\norm{a}<1$ and $\langle (a-1)t\rangle\not\in a \OO$,}\\
\end{array}\right.
\end{array}
\end{equation}
and hence via scaling rule (\ref{equation:ThetaScaling}) for the theta symbol and the definitions
we have

\begin{equation}\label{equation:BaseThetaCalcBis}
\begin{array}{cl}
&\Theta(a,\Psi_t-\Psi_1)=\Theta(1,\Psi_t^{(1,a)}-\Psi_1^{(1,a)})\\\\
=&\displaystyle\left\{\begin{array}{cl}
1&\mbox{if $\norm{a}<1$ and $\langle (a-1)t\rangle\in a \OO$,}\\
0&\mbox{otherwise.}
\end{array}\right.
\end{array}
\end{equation} 
Since $\Theta(a,\Psi_t-\Psi_1)=0$ for $\max(\norm{a},\norm{a}^{-1})\gg 0$,
we have verified that $\Psi_t-\Psi_1$ is proper.

\subsubsection{Calculation of the rational Fourier transform}
Under special assumption (\ref{equation:SpecialAssumption}) we have
$$L^\perp=L,\;\;\;\DDD=\OO,\;\;\;\mu(k/L)=\mu(\OO),$$
and hence 
\begin{equation}\label{equation:SampleRationalFourier}
\GGG_0[\Psi_t]_*(x,y)=
\lambda_0(\Res(t(y-x)d\pi))\one_L(x)\one_{\OO}(y)
\end{equation}
by Lemma~\ref{Lemma:RationalRVFT}, where (recall that)
 we have
$$\lambda_0(C)=\left\{\begin{array}{rl}
1&\mbox{if $C=0$,}\\
-1&\mbox{if $C=1$,}\\
0&\mbox{if $C\not\in \{0,1\}$}
\end{array}\right.
$$
for all $C\in \FF_q$.  

\subsubsection{Analysis of the rational Fourier transform}
From (\ref{equation:SampleRationalFourier}) via scaling rule
(\ref{equation:RationalFourierScalingRule}), we deduce that 
$$
\begin{array}{rcl}
\GGG_0[\Psi_t^{(1,a)}]_*(x,y)&=&\norm{a}\GGG_0[\Psi_t]_*(x,ay)\\\\
&=&
\norm{a}\lambda_0(\Res(t(ay-x)d\pi))\one_L(x)\one_{a^{-1}\OO}(y)
\end{array}
$$
and hence
$$
\GGG_0[\Psi_t^{(1,a)}]_*(x,x)=\norm{a}(\one_{L(t,a,0)}(x)-\one_{L(t,a,1)}(x)),
$$
where
$$
L(t,a,C)=\left\{\ell\in L\cap a^{-1} \OO\left|
\Res\left((a-1)t\ell\,d\pi\right)=C\right.\right\}
$$
for every $C\in \FF_q$.
Also put 
$$L(t,a,*)=\bigcup_{C\in \FF_q^\times}L(t,a,C).$$
Note that
\begin{equation}\label{equation:EllNoteThat}
L(t,a,*)=\emptyset\Leftrightarrow \langle (a-1)t\rangle\in a\OO,\end{equation}
\begin{equation}\label{equation:Trick1}
\one_{L(t,a,0)}-\one_{L(1,a,0)}=\one_{L(1,a,*)}-\one_{L(t,a,*)},
\end{equation}
\begin{equation}\label{equation:Trick2}
L(1,a,1)=
\left\{\begin{array}{rl}
-\pi^{-1}+\pi^{-1}L\cap a^{-1}\OO&\mbox{if $\norm{a}<1$,}\\
\emptyset&\mbox{if $\norm{a}\geq 1$.}
\end{array}\right.
\end{equation}
It follows that
 \begin{equation}\label{equation:GenFourierCalc}
\begin{array}{rl}
&\GGG_0[\Psi_t^{(1,a)}-\Psi_1^{(1,a)}]_*(x,x)\\\\
=&
\left\{\begin{array}{ll}
0\\
\mbox{if $\norm{a}\geq 1$,}\\\\
\norm{a}(\one_{-\pi^{-1}+\pi^{-1}L\cap a^{-1}\OO}+\sum_{C\in \FF_q^\times}\one_{C\pi^{-1}+\pi^{-1}L\cap a^{-1}\OO})(x)\\
\mbox{if $\norm{a}<1$ and $\langle (a-1)t\rangle\in a \OO$,}\\\\
\norm{a}(\one_{L(1,a,1)}+\one_{L(1,a,*)}-
\one_{L(t,a,1)}-\one_{L(t,a,*)}
)(x)
\\\mbox{if $\norm{a}<1$ and $\langle (a-1)t\rangle\not\in a \OO$}.
\end{array}\right.
\end{array}
\end{equation}

\subsubsection{Calculation of the Catalan symbol}
We now obtain a formula for the Catalan symbol $\left(\begin{subarray}{c}
a\\\Psi_t-\Psi_1\end{subarray}\right)$ by combining 
formulas (\ref{equation:DiagonalPsiT}--\ref{equation:GenFourierCalc}), Wilson's theorem $\prod_{C\in \FF_q^\times}C=-1$,
the definitions, the Moore identity (\ref{equation:MooreIdentity}),
the Ore identity (\ref{equation:OreIdentity}),
and  the scaling rule (\ref{equation:CatalanScaling}) for the Catalan symbol.
Write 
$$\norm{a}=\norm{\pi^{-N}}=\norm{\pi^{\nu+1}}.$$ 
We have
\begin{equation}\label{equation:BaseSymbolCalcBisSimp}
\begin{array}{cl}
&\displaystyle \left(\begin{array}{c}
a\\
\Psi_t-\Psi_1
\end{array}\right)=\left(\begin{array}{c}
1\\
\Psi_t^{(1,a)}-\Psi_1^{(1,a)}\end{array}\right)\\\\
=&
\left\{\begin{array}{l}
\displaystyle 
\frac{\Moore(\lfloor at\rfloor +\langle t\rangle,\pi^{-1},\dots,\pi^{-N})}
{\Moore(1,\pi^{-1},\dots,\pi^{-N})}\\\\
\mbox{if $\norm{a}\geq 1$,}\\\\\\
(\lfloor at \rfloor +\langle t\rangle)
\displaystyle\left(\frac{\Moore(\pi^{-1},\dots,\pi^{-\nu-1})}
{\Moore(\pi^{-2},\dots,\pi^{-\nu-1})}\right)^{q^{-\nu}}\\\\
\mbox{if $\norm{a}<1$ and $\Theta(a,\Psi_t-\Psi_1)=1$,}\\\\\\
\displaystyle\left(\frac{\Ore(t(1-a)d\pi,\pi^{-1},\dots,\pi^{-\nu-1})}
{\Ore(d\pi,\pi^{-1},\dots,\pi^{-\nu-1})}\right)^{q^{-\nu}}\\\\
\mbox{if $\norm{a}<1$
and $\Theta(a,\Psi_t-\Psi_1)=0$.}
\end{array}\right.
\end{array}
\end{equation}

\subsection{Construction of a candidate}
We construct the interpolating gadget 
ultimately to be yoked to $\Psi_t-\Psi_1$.

\subsubsection{Special elements of $\Lambda$}
Put
$$\xi=\xi_M$$
to simplify writing below. Put
$$\left.\begin{array}{rcl}
X_i&=&U_{\xi_i,\xi}(X)\\
Y_i&=&U_{\xi_i,\xi}(Y)
\end{array}\right\}\in \FF_q[[X,Y]]\subset \Lambda\;\;\;\mbox{for $i=0,\dots,M$.}$$
We have
\begin{equation}\label{equation:XYTowerRelations}
\begin{array}{rcl}
X_{i-1}&=&\left\{\begin{array}{cl}
X_0X_i+X_i^q&\mbox{if $i>1$,}\\
-X_1^{q-1}&\mbox{if $i=1$,}
\end{array}\right.\\\\
Y_{i-1}&=&\left\{\begin{array}{cl}
Y_0Y_i+Y_i^q&\mbox{if $i>1$,}\\
-Y_1^{q-1}&\mbox{if $i=1$,}
\end{array}\right.
\end{array}\;\;\mbox{for $i=1,\dots,M$.}
\end{equation}
By the explicit reciprocity law (\ref{equation:ExplicitReciprocityLaw}) and the definitions, we have
\begin{equation}\label{equation:ExplicitReciprocityXY}
\begin{array}{rcl}
\displaystyle X_i^{[b^{-1}\pi^{-\beta},c^{-1}\pi^{-\gamma};\xi]}
&=&\displaystyle\left\{\begin{array}{rl}
X_0^{q^\beta}&\mbox{if $i=0$,}\\
\sum_{j=0}^{i-1}b_jX_{i-j}^{q^\beta}&\mbox{if $i>0$,}
\end{array}\right.\\\\
\displaystyle Y_i^{[b^{-1}\pi^{-\beta},c^{-1}\pi^{-\gamma};\xi]}
&=&\left\{\begin{array}{rl}
Y_0^{q^\gamma}&\mbox{if $i=0$,}\\
\sum_{j=0}^{i-1}c_jY_{i-j}^{q^\gamma}&\mbox{if $i>0$,}
\end{array}\right.
\end{array}
\end{equation}
for all
$$b=\sum_{i=0}^\infty b_i\pi^i,\;\;
c=\sum_{i=0}^\infty c_i\pi^i\;\;\;(b_i,c_i\in \FF_q,\;\;b_0\neq 0,\;c_0\neq 0),$$
 integers $\beta,\gamma\geq 0$,
and $i=0,\dots,M$.

\subsubsection{Generating functions}
Put
$$\left.\begin{array}{rcl}
\displaystyle f=f(X,Y,Z)=\sum_{n=0}^\infty f_n(X,Y)Z^n&=&\displaystyle\frac{\displaystyle \prod_{i=0}^\infty
(1-Y_0^{q^i}Z)}{\displaystyle \prod_{i=0}^\infty
(1-X_0^{q^i}Z)}\\\\
\displaystyle g=g(X,Y,Z)=\sum_{m=0}^\infty g_m(X,Y)Z^m&=&\displaystyle
\frac{\displaystyle \sum_{i=0}^{M-1}X_{i+1}Z^i}{\displaystyle
\sum_{i=0}^{M-1}Y_{i+1}Z^i}
\end{array}\right\}\in \Lambda[[Z]]^\times.$$
We remark that
$$f_n(X,Y)\in (X_0,Y_0)^n\subset \FF_q[[X,Y]],\;\;\;g_m(X,Y)\in \FF_q[[X,Y]][X^{-1},Y^{-1}]$$
for all $m,n\geq 0$.
\subsubsection{Definition of $F_t$}
Write
$$t=\sum_{i=-M}^\infty t_i\pi^i\;\;\;(t_i\in \FF_q),$$
and then put
$$F_t=F_t(X,Y)=\sum_{i=0}^\infty t_if_i(X,Y)-\sum_{i=1}^M t_{-i}g_{i-1}(X,Y)\in \Lambda.$$
The sum $\sum_{i=0}^\infty t_if_i(X,Y)$ converges $(X_0,Y_0)$-adically in $\FF_q[[X,Y]]$
and hence $F_t$ is well-defined. We remark
that in the case $M=0$, one simply ignores here
and below all functions and formulas connected with $g$.

\subsubsection{An explicit reciprocity law obeyed by $F_t$}
Given $a\in \OO$, we write $(a\vert Z)\in \FF_q[[Z]]$ for the result of expanding
$a$ as a power series in $\pi$ with coefficients in $\FF_q$, and then replacing $\pi$ by $Z$.
Let us agree to extend the ``square bracket operations'' from $\Lambda$ to $\Lambda[[Z]]$
coefficient-by-coefficient.   From  (\ref{equation:ExplicitReciprocityXY}) we deduce that 
$$\begin{array}{rcl}
\displaystyle f^{[b\pi^{-\beta},c\pi^{-\gamma};\xi]}(X,Y,Z)&=&
\displaystyle f(X^{q^\beta},Y^{q^\gamma},Z),\\\\
\displaystyle
g^{[b\pi^{-\beta},c\pi^{-\gamma};\xi]}(X,Y,Z)&\equiv_M
&\displaystyle\left(\frac{c}{b}\vert Z\right)g(X^{q^\beta},Y^{q^\gamma},Z)
\end{array}
$$
where $\equiv_M$ denotes congruence modulo the ideal $(Z^M)\subset \Lambda[[Z]]$,
and hence
\begin{equation}\label{equation:UpdatedReciprocity}
F_t^{[b\pi^{-\beta},c\pi^{-\gamma};\xi]}(X,Y)=
F_{\lfloor ct/b\rfloor+\langle t\rangle}(X^{q^\beta},Y^{q^\gamma})
\end{equation}
for all $b,c\in \OO^\times$ and  integers $\beta,\gamma\geq 0$. In particular, it is clear that the pair $(F_t^{[1,\pi^{-1};\xi]},\xi)$ is an interpolating gadget.  We will yoke $(F_t^{[1,\pi^{-1};\xi]},\xi)$ to $\Psi_t-\Psi_1$.
\subsubsection{Reduction of the proof} Note that
\begin{equation}\label{equation:PreliminaryRed}
(\Psi_t-\Psi_1)^{(1,u)}=\Psi_{\lfloor ut\rfloor+\langle t\rangle}-\Psi_1
\end{equation}
for all $u\in \OO^\times$ via scaling rule (\ref{equation:RiggedScaling}).
On account of (\ref{equation:PreliminaryRed})
and its analogue (\ref{equation:UpdatedReciprocity}),
along with the scaling rules obeyed by the four types of symbols, we have only to prove that
 \begin{equation}\label{equation:InterpolationNuff}
\begin{array}{rcl}
\displaystyle\Theta_\sh(\pi^{-N-1},F_t,\xi)&=&
\displaystyle\Theta(\pi^{-N},\Psi_t-\Psi_1),\\\\
\displaystyle\left(\begin{array}{c}
\pi^{-N-1}\\
F_t,\xi
\end{array}\right)_\sh&=&\displaystyle\left(\begin{array}{c}
\pi^{-N}\\
\Psi_t-\Psi_1
\end{array}\right)
\end{array}
\end{equation}
for all integers $N$ in order to yoke $(F^{[1,\pi^{-1};\xi]}_t,\xi)$ to $\Psi_t-\Psi_1$, thereby completing the proof of the theorem.

\subsection{Completion of the proof}
Throughout the calculations remaining,
let $N$, $\nu$ and $n$ be nonnegative integers
and let $m$ be an integer in the range $0\leq m<M$.

\subsubsection{Specialization of $f$ and $g$}
We write
$$\frac{\partial F}{\partial X_0 } =\frac{1}{U'_{\pi,\xi}(X)}\frac{\partial F}{\partial X}$$
for all $F\in \Lambda[[Z]]$. 
From (\ref{equation:XYTowerRelations}) we deduce the functional 
and differential equations
$$
\begin{array}{rcl}
f(X,Y,Z)&=&\left(1-Y_0Z\right)f(X,Y^q,Z),\\
f(X^q,Y,Z)&=&(1-X_0 Z)f(X,Y,Z),\\
0&=&\displaystyle-Zf(X,Y,Z)+\left(1-X_0 Z\right)\frac{\partial f}{\partial X_0 }(X,Y,Z),\\\\
g(X,Y,Z)&\equiv_M&\left(Z-Y_0 \right)g(X,Y^q,Z),\\
g(X^q,Y,Z)&\equiv_M &\left(Z-X_0 \right)g(X,Y,Z),\\
0&\equiv_M&\displaystyle-g(X,Y,Z)+\left(Z-X_0 \right)\frac{\partial g}{\partial X_0 }(X,Y,Z).\end{array}
$$
Here, as previously, $\equiv_M$ denotes congruence
modulo the ideal $(Z^M)\subset \Lambda[[Z]]$.
We deduce by iteration and specialization the following list of identities in the power series ring $k[[Z]]$:
\begin{equation}\label{equation:fgFEbis}
\begin{array}{rll}
f(\xi,\xi^{q^{N+1}},Z)&=&\prod_{i=0}^N\frac{1}{1-\pi^{q^i}Z}\\
f(\xi^{q^\nu},\xi,Z)&=&\prod_{i=0}^{\nu-1}(1-\pi^{q^i}Z)\\
\frac{\partial f}{\partial X_0 }(\xi^{q^\nu},\xi,Z)&=&
\frac{Z}{1-\pi^{q^\nu}Z}\prod_{i=0}^{\nu-1}(1-\pi^{q^i}Z)\\\\
g(\xi,\xi^{q^{N+1}},Z)&\equiv_M&\prod_{i=0}^{N}\frac{1}{Z-\pi^{q^i}}\\
g(\xi^{q^\nu},\xi,Z)&\equiv_M& \prod_{i=0}^{\nu-1}(Z-\pi^{q^i})\\
\frac{\partial g}{\partial X_0 }(\xi^{q^\nu},\xi,Z)&\equiv_M&
\frac{1}{Z-\pi^{q^\nu}}\prod_{i=0}^{\nu-1}(Z-\pi^{q^i})
\end{array}
\end{equation}
Here and below $\equiv_M$ denotes congruence
modulo the ideal $(Z^M)\subset k[[Z]]$.
\subsubsection{Application of Jacobi-Trudi identity}
One deduces from first, second, fourth and fifth identities on the list
(\ref{equation:fgFEbis}) via the special cases (\ref{equation:JacobiTrudi})
and (\ref{equation:JacobiTrudiBis}) of the Jacobi-Trudi identity
 that
$$\begin{array}{rcl}
f_n(\xi,\xi^{q^{N+1}})&=&
\frac{\Moore(\pi^n,\pi^{-1},\dots,\pi^{-N})}
{\Moore(1,\pi^{-1},\dots,\pi^{-N})},\\
f_n(\xi^{q^\nu},\xi)&=&
\frac{\Ore(\pi^n d\pi,\pi^{-1},\dots,\pi^{-\nu-1})}
{\Ore(d\pi,\pi^{-1},\dots,\pi^{-\nu-1})},\\
g_m(\xi,\xi^{q^{N+1}})&=&
-\frac{\Moore(\pi^{-(m+1)-N},\pi^{-1},\dots,\pi^{-N})}
{\Moore(1,\pi^{-1},\dots,\pi^{-N})},\\
g_m(\xi^{q^\nu},\xi)&=&\frac{\Ore(\pi^{-(m+1)+\nu+1} d\pi,\pi^{-1},\dots,\pi^{-\nu-1})}
{\Ore(d\pi,\pi^{-1},\dots,\pi^{-\nu-1})},\\
\end{array}
$$
and hence
\begin{equation}\label{equation:ClincherOne}
F_t(\xi,\xi^{q^{N+1}})=
\frac{\Moore(\lfloor \pi^{-N}t\rfloor+\langle t\rangle,\pi^{-1},\dots,\pi^{-N})}{
\Moore(1,\pi^{-1},\dots,\pi^{-N})},
\end{equation}
\begin{equation}\label{equation:ClincherTwo}
F_t(\xi^{q^\nu},\xi)=
\frac{\Ore((1-\pi^{\nu+1})t\,d\pi,\pi^{-1},\dots,\pi^{-\nu-1})}
{\Ore(d\pi,\pi^{-1},\dots,\pi^{-\nu-1})},
\end{equation}
by definition of $F_t$. From the latter formula it follows that
\begin{equation}\label{equation:ClincherThree}
F_t(\xi^{q^\nu},\xi)=0\Leftrightarrow \langle (1- \pi^{\nu+1})t\rangle
\in \pi^{\nu+1}\OO.
\end{equation}

\subsubsection{Application of some residue formulas}
We have
$$
\begin{array}{rcl}
\frac{\partial f_n}{\partial X_0 }(\xi^{q^\nu},\xi)&=&-
\Res_{Z=\infty}\left(\frac{Z^{n}}{Z-\pi^{q^\nu}}(\prod_{i=0}^{\nu-1}
(1-\pi^{q^i}/Z))\frac{dZ}{Z}\right),\\\\
\frac{\partial g_m}{\partial X_0 }(\xi^{q^\nu},\xi)&=&
\Res_{Z=0}\left(\frac{Z^{-(m+1)+\nu+1}}{Z-\pi^{q^\nu}}
(\prod_{i=0}^{\nu-1}(1-\pi^{q^i}/Z))\frac{dZ}{Z}\right),\\\\
\pi^{n q^{\nu}}\frac{\Moore(\pi^{-1},\dots,\pi^{-\nu-1})}
{\Moore(\pi^{-2},\dots,\pi^{-\nu-1})}&
=&
\Res_{Z=\pi^{q^{\nu}}}\left(\frac{Z^{n}}{Z-\pi^{q^{\nu}}}
\left(\prod_{i=0}^{\nu-1}(1-\pi^{q^i}/Z)\right)\frac{dZ}{Z}
\right).
\end{array}
$$
The first two formulas follow directly from
the third and sixth identities on the list (\ref{equation:fgFEbis}),
while
the last is an application of the Vandermonde determinant identity.
By a straightforward if quite tedious calculation exploiting ``sum-of-residues-equals-zero'' one verifies that
\begin{equation}\label{equation:ClincherFour}
\begin{array}{cl}
&\langle (1- \pi^{\nu+1})t\rangle
\in \pi^{\nu+1}\OO\\\\
\Rightarrow &\displaystyle\frac{\partial F_t}{\partial X_0}(\xi^{q^\nu},\xi)=
(\lfloor \pi^{\nu+1}t\rfloor+\langle t\rangle)^{q^\nu}\frac{\Moore(\pi^{-1},\dots,\pi^{-\nu-1})}
{\Moore(\pi^{-2},\dots,\pi^{-\nu-1})}.
\end{array}
\end{equation}

\subsubsection{Last details}
By definition $\Theta_\sh(\pi^{-N-1},F_t,\xi)$
is the  $(X-Y)$-valuation of $F_t(X,Y^{q^{N+1}})$, hence
$\Theta_\sh(\pi^{-N-1},F_t,\xi)$ vanishes by (\ref{equation:ClincherOne}),
and hence  \linebreak $\Theta_\sh(\pi^{-N-1},F_t,\xi)=\Theta(\pi^{-N},\Psi_t-\Psi_1)$
by (\ref{equation:BaseThetaCalcBis}). By definition $q^\nu\Theta_\sh(\pi^\nu,F_t,\xi)$
is the $(X-Y)$-valuation of $F_t(X^{q^\nu},Y)$, hence $\Theta_\sh(\pi^\nu,F_t,\xi)$
is the $(X-Y^{q^\nu})$-valuation of $F_t$,
 hence
$\Theta_\sh(\pi^\nu,F_t,\xi)=1,0$ according as $\langle (1- \pi^{\nu+1})t\rangle$
does or does not belong to $\pi^{\nu+1}\OO$, respectively, by (\ref{equation:ClincherThree})
and (\ref{equation:ClincherFour}), and hence
$\Theta_\sh(\pi^\nu,F_t,\xi)=\Theta(\pi^{\nu+1},\Psi_t-\Psi_1)$
by (\ref{equation:BaseThetaCalcBis}).
Recall that since we made our calculations of Catalan symbols under the special assumption
(\ref{equation:SpecialAssumption}), and we are obliged to make our calculations of shadow Catalan symbols under the same special assumption.
By definition $\left(\begin{subarray}{c}
\pi^{-N-1}\\
F_t,\xi
\end{subarray}\right)=F(\xi,\xi^{q^{N+1}})$ and hence
$\left(\begin{subarray}{c}
\pi^{-N-1}\\
F_t,\xi
\end{subarray}\right)=\left(\begin{subarray}{c}
\pi^{-N}\\
\Psi_t-\Psi_1
\end{subarray}\right)$ by (\ref{equation:BaseSymbolCalcBisSimp})
and (\ref{equation:ClincherOne}). By definition $\left(\begin{subarray}{c}
\pi^\nu\\
F_t,\xi
\end{subarray}\right)^{q^\nu}=\frac{\partial F_t}{\partial X_0}(
\xi^{q^\nu},\xi),F(\xi^{q^\nu},\xi)$
according as $\Theta_\sh(\pi^\nu,F_t,\xi)=\Theta(\pi^{\nu+1},\Psi_t-\Psi_1)=1,0$, 
respectively, and hence $\left(\begin{subarray}{c}
\pi^\nu\\
F_t,\xi
\end{subarray}\right)=\left(\begin{subarray}{c}
\pi^{\nu+1}\\
\Psi_t-\Psi_1
\end{subarray}\right)$ by (\ref{equation:BaseThetaCalcBis}), (\ref{equation:BaseSymbolCalcBisSimp}),
(\ref{equation:ClincherThree}) and (\ref{equation:ClincherFour}).
Thus (\ref{equation:InterpolationNuff})
holds, and the proof of Theorem~\ref{Theorem:InterpolationExamples} is complete.
\qed

\section{The interpolation theorem}

 \begin{Theorem}\label{Theorem:Interpolability}
Every proper rational rigged virtual lattice
is interpolable.
\end{Theorem}
\noindent This is our main result.  The proof of Theorem~\ref{Theorem:Interpolability} takes up the rest of this section.  Fortunately the proof of Theorem~\ref{Theorem:Interpolability}
requires none of the details of the explicit constructions used to prove Theorem~\ref{Theorem:InterpolationExamples}---we will need only the bare statement of the latter.
Before commencing the proof we prove two corollaries
and make a remark.

\begin{Corollary}\label{Corollary:Immediate}
Every proper rational rigged virtual lattice 
is a finite $\ZZ[1/p]$-linear combination
of asymptotically interpolable rigged virtual lattices.
\end{Corollary}
\noindent This is a reiteration of Theorem~\ref{Subtheorem:Immediate}. \proof By Proposition~\ref{Proposition:InterpolabilitySimplification} every interpolable rigged virtual lattice is a $\ZZ[1/p]$-linear combination of strictly interpolable rigged virtual lattices, and as we have already noted strict interpolability implies asymptotic interpolability.  So the corollary follows directly from the theorem.   \qed

\begin{Corollary}\label{Corollary:AnalyticContinuation}
Let $\Phi$ be a proper rational rigged virtual lattice.
Assume that $\Phi$ is asymptotically interpolable and asymptotically
yoked to an interpolating gadget $(F,\xi)$. Then $\Phi$ is strictly interpolable
and yoked to $(F,\xi)$. 
\end{Corollary}
\proof Choose a basepoint $\eta\in \PPP^*$ such that $\xi\leq \eta$
and
$\Phi$ is of conductor $\leq \eta$.
By the theorem and Proposition~\ref{Proposition:InterpolabilitySimplification}, for some power $r$ of the characteristic $p$ of $k$, we can write
$r\Phi=\Phi_+-\Phi_-$, where  $\Phi_{\pm}$ are strictly interpolable rigged virtual lattices yoked to  interpolating gadgets $(G_{\pm},\eta)$, respectively. Put $G=F^{[1;\xi,\eta]}$.
Note that $(G,\eta)$ is asymptotically yoked to $\Phi$.
It follows that $(G^rG_-,\eta)$
and $(G_+,\eta)$ are both asymptotically yoked
to $\Phi_+$, and hence equal by Proposition~\ref{Proposition:InterpolatingUniqueness}.
It follows in turn that $(G,\eta)$ is yoked to $\Phi$.
Finally, by the scaling rules for the shadow symbols,
$(F,\xi)$ is yoked to $\Phi$, as claimed.
\qed

\subsection{Remark} Fix a uniformizer $\pi\in \OO$.
Let $W\subset k$ be a cocompact discrete \linebreak $\FF_q$-subspace
such that $\mu(k/W)=\mu(\pi \OO)$.
Consider the rigged virtual lattice $\Phi$
of the form
$$\Phi_*=\one_W\otimes \one_{1+\pi\OO}-
\one_{\FF_q[\pi^{-1}]}\otimes \one_{1+\pi\OO}.$$
It is easy to see that $\Phi$ is rational,
proper, separable and effective. By Theorem~\ref{Theorem:Interpolability},
the rigged virtual lattice $\Phi$ is interpolable,
and by Proposition~\ref{Proposition:StrictnessCriterion} it follows that $\Phi$ is strictly interpolable. By means of Proposition~\ref{Subproposition:GaloisDescent}
and the scaling rules for the Catalan symbol and its shadow it can be verified that $\Phi$ is yoked to an interpolating gadget of the form $(\tau_W,\sqrt[q-1]{-\pi})$ where $\tau_W=\tau_W(X,Y)\in \FF_q[[X,Y]][X^{-1},Y^{-1}]$.
The power series $\tau_W$ is essentially the same as the one figuring in \cite[Theorem 2.4.1]{AndersonAHarmonic}, and so admits interpretation as a {\em $\tau$-function}.
The method used to prove Theorem~\ref{Theorem:Interpolability} is actually a refinement of that used to prove
\cite[Theorem 2.4.1]{AndersonAHarmonic}. 

\subsection{Setting for the proof of the main result}
Fix a uniformizer $\pi\in \OO$.
Fix a sequence $\{\xi_i\}_{i=0}^\infty$ in $k_\ab$ satisfying the relations
$$\xi_0=\pi,\;\;\xi_{i-1}=
\left\{\begin{array}{rl}
-\xi_1^{q-1}&\mbox{if $i=1$}\\
\pi\xi_i+\xi_i^q&\mbox{if $i>1$}
\end{array}\right.
$$
for $i>0$.  For $t\in k$, let $\Psi_t$ be the unique primitive rigged virtual lattice such that
$$(\Psi_t)_*=\one_{t+\pi^{-1}\FF_q[\pi^{-1}]}\otimes \one_{t+\OO}.$$
 This is more or less
the same setting as that in which we proved Theorem~\ref{Theorem:InterpolationExamples}. However, 
in the present setting:  we allow $t=0$ in the definition of $\Psi_t$; we reserve the symbol ``$L$'' for denoting a
general cocompact discrete subgroup of $k$; and we do not make the special assumption (\ref{equation:SpecialAssumption}).

\subsection{{\em Ad hoc} terminology}

\subsubsection{$\pi$-regularity and level}
Given a rigged virtual lattice $\Phi$ and an integer $M\geq 0$,
we say that $\Phi$ is {\em $\pi$-regular of level $\leq M$}
if it is possible to decompose
$\Phi$ as a finite $\ZZ[1/p]$-linear combination
of primitive rigged virtual lattices of the form $\Psi_t^{(\pi^i,\pi^j)}$
for  $t\in k$ such that $\norm{t}\leq \norm{\pi^{-M}}$
and integers $i$ and $j$.
If $\Phi$ is $\pi$-regular of some level we say that
$\Phi$ is {\em $\pi$-regular}.
Note that if $\Phi$ is a primitive rigged virtual lattice and
$\Phi_*=\one_{\ell+L}\otimes \one_{w+r\OO}$
where $L\subset k$ is a cocompact discrete subgroup, 
$\ell,w\in k$ and $r\in k^\times$, a sufficient condition  for $\Phi$ to be $\pi$-regular
of level $\leq M$ is that
$\norm{\pi^{-M}}\geq \norm{w/r}$
and $L\supset\pi^{-N}\FF_q[\pi^{-1}]$
for some $N$.

\subsubsection{Softness}
We say that a rigged virtual lattice $\Phi$ is {\em soft} if $\Phi$ is rational, 
$\Theta(a,\Phi)=0$ for  all $a\in k^\times$ and $\Phi_*(0,t)=0=\GGG_0[\Phi]_*(0,t)$ for all $t\in k$. Note that softness implies properness and effectiveness.
If $\Phi$ is 
soft, strictly interpolable and yoked to an interpolating gadget of the form $(F,\xi_M)$
for some integer $M\geq 0$, then  $F\in \Lambda_0^\times$
by Proposition~\ref{Proposition:Underscore} and moreover $F\in \FF_q[[X,Y]]^\times$
by Proposition~\ref{Subproposition:TotallyRamified}.

\begin{Proposition}
Let $\Phi$ be a proper rigged virtual lattice.
If $\Phi$ can be written as a finite $\ZZ[1/p]$-linear combination of primitive rigged virtual lattices of the form $\Psi_e^{(\pi^i,1)}$ 
where $e\in \{0,1\}$ and $i\in \ZZ$, then $\Phi$ is interpolable of conductor $\leq \pi$.
\end{Proposition}
\proof After replacing $\Phi$
by $\Phi^{(\pi^\beta,0)}$ for some integer $\beta\geq 0$, we may assume that $\Phi$ is a finite $\ZZ[1/p]$-linear combination of rigged virtual lattices of the form $\Psi_e^{(\pi^i,1)}$ where $e\in \{0,1\}$
and $i\in \ZZ\cap[0,\infty)$.
Now for every integer $i\geq 0$ and $e\in \{0,1\}$ we have
$$\Psi_e^{(\pi^i,1)}=
q^{-i}\sum_{t\in \pi^{i-1}\FF_q[\pi^{-1}]\cap \OO}
\Psi_{\pi^{i}e+t}$$
by scaling rule (\ref{equation:RiggedScaling})
and the definitions. Therefore we have
$$\Phi=\sum_{i=1}^N \alpha_i \Psi_{t_i}$$
for some numbers $\alpha_i\in \ZZ[1/p]$
and $t_i\in \OO$. After grouping terms
we may assume that the $t_i$ are distinct. 
Now in general we have
\begin{equation}\label{equation:LocalNeed}
(\Psi_t)_*(0,0)=\delta_{t0},\;\;\;
\GGG_0[\Psi_t]_*(0,0)=1,
\end{equation}
by the definitions in the former case
and by Lemma~\ref{Lemma:RationalRVFT} in the latter case.
Since $\Phi_*(0,0)=0$ by hypothesis, we may assume by (\ref{equation:LocalNeed}) that $t_i\in k^\times$
for all $i$.
Further, since $\GGG_0[\Phi]_*(0,0)=0$ by hypothesis, we have
$\sum_{i=1}^n \alpha_i=0$ by (\ref{equation:LocalNeed}) 
and hence
$$\Phi=\sum_{i=1}^n \alpha_i(\Psi_{t_i}-\Psi_1).$$
Thus $\Phi$ is interpolable of conductor $\leq \pi$ by Theorem~\ref{Theorem:InterpolationExamples}.
  \qed

\begin{Proposition}\label{Proposition:piRegularInterpolable}
Let $\Phi$ be a proper rigged virtual lattice.
If $\Phi$ is  $\pi$-regular of level $\leq M$, then $\Phi$ is interpolable of conductor $\leq \xi_M$.
\end{Proposition}
\proof By hypothesis, after some evident rearrangement, we have
$$\begin{array}{rcl}
\Phi&=&\displaystyle\sum_{i=1}^N \alpha_i(\Psi_{t_i}-\Psi_{e_i})^{(\pi^{\beta_i},\pi^{\gamma_i})}\\\\
&&\displaystyle+\sum_{i=1}^N \alpha_i(
\Psi_{e_i}^{}-
q^{-\gamma_i}\Psi_{e_i}^{(\pi^{-\gamma_i},\pi^{-\gamma_i})})^{(\pi^{\beta_i},\pi^{\gamma_i})}\\\\
&&\displaystyle+\sum_{i=1}^N \alpha_iq^{-\gamma_i}\Psi_{e_i}^{(\pi^{\beta_i-\gamma_i},1)}
\end{array}$$
where
$$\alpha_i\in \ZZ[1/p],\;\;\beta_i,\gamma_i\in \ZZ,\;\;t_i\in k,\;\;
e_i=\delta_{t_i,0},\;\;\;\norm{t_i}\leq \norm{\pi^{-M}}.$$  
 Call the sums  on the right $\Phi_1$, $\Phi_2$ and $\Phi_3$, respectively. The sum $\Phi_1$ is interpolable of conductor $\leq \xi_M$ by Theorem~\ref{Theorem:InterpolationExamples}. The sum $\Phi_2$ is interpolable of conductor $\leq \pi$ by Proposition~\ref{Proposition:CheapExamples}. Since $\Phi$, $\Phi_1$ and $\Phi_2$ are proper, so is $\Phi_3$.
Finally,  $\Phi_3$ is interpolable of conductor $\leq \pi$ 
by the preceding proposition. \qed

\begin{Proposition}\label{Proposition:InterpolationLimit}
Let $\{\Phi\}\cup\{\Phi_i\}_{i=1}^\infty$ be a family of soft rigged virtual lattices. 
Let $M\geq 0$ be an integer.
Assume that for every $a\in k^\times$ we have
$$\lim_{i\rightarrow\infty}
\left\Vert\left(\begin{array}{c}
a\\
\Phi-\Phi_i
\end{array}\right)-1\right\Vert=0.$$
Assume that $\Phi_i$ for every $i$ is strictly interpolable and yoked to an interpolating gadget with basepoint $\xi_M$. 
Then $\Phi$ is strictly interpolable and yoked to an interpolating gadget with basepoint $\xi_M$.
\end{Proposition}
\proof For each $i$ let $(F_i,\xi_M)$ be an interpolating gadget yoked to $\Phi_i$,
noting that $F_i\in \FF_q[[X,Y]]^\times$ for all $i$.
 Assume for the moment that $F=\lim F_i$ exists \linebreak $(X,Y)$-adically in $\FF_q[[X,Y]]$. Then necessarily $F\in \FF_q[[X,Y]]^\times$, and moreover
for every $a\in k^\times$ we have convergence
$$\left(\begin{array}{c}
a\\
F_i,\xi_M
\end{array}\right)_\sh=
F_i(\xi_M,(\rho(a)^{-1}\xi_M)^{\norm{a}})$$
$$\rightarrow_{i\rightarrow\infty}
F(\xi_M,(\rho(a)^{-1}\xi_M)^{\norm{a}}) =\left(\begin{array}{c}
a\\
F,\xi_M
\end{array}\right)_\sh$$
with respect to $\norm{\cdot}$.
By hypothesis we have the corresponding convergence of Catalan symbols. Therefore $(F,\xi_M)$
must be yoked to $\Phi$. Thus in order to prove the proposition it is enough just to show 
that $\lim F_i$ exists $(X,Y)$-adically in $\FF_q[[X,Y]]$. Fix a positive integer $N$ arbitrarily. For each $i$ use the Weierstrass Division Theorem to write
$$-1+F_{i+1}(X,Y)/F_i(X,Y)=R_i(X,Y)+Q_i(X,Y)\cdot\prod_{j=1}^N (Y-X^{q^j})$$
where $Q_i(X,Y)\in \FF_q[[X,Y]]$ and $R_i(X,Y)\in \FF_q[[X]][Y]$
is of degree $<N$ in $Y$. By hypothesis $$R_i(\xi_M,\xi_M^{q^j})=-1+\left(\begin{array}{c}
\pi^{-j}\\
F_{i+1}/F_i,\xi_M
\end{array}\right)_\sh=
-1+\left(\begin{array}{c}
\pi^{-j}\\
\Phi_{i+1}-\Phi_i
\end{array}\right)\rightarrow_{i\rightarrow\infty}0$$
 for $j=1,\dots,N$. By applying the Lagrange Interpolation Theorem 
 one deduces that $R_i(\xi_M,Y)\in k(\xi_M)[Y]$
 tends coefficient-by-coefficient to $0$ with respect to $\norm{\cdot}$ as $n\rightarrow\infty$, and hence that $R_i(X,Y)$ tends $X$-adically to $0$ as $i\rightarrow\infty$. We conclude that $F_{i+1}/F_i\equiv 1\bmod{(X,Y)^N}$ for all $i\gg 0$.
Since $N$ was arbitrarily chosen, the $(X,Y)$-adic convergence of the sequence $\{F_i\}_{i=1}^\infty$ is proved, and with it the proposition. \qed

\subsection{An approximation scheme}\label{subsection:Approximation}
We give the set up for the last proposition of the paper, which is the core of the proof of Theorem~\ref{Theorem:Interpolability}.
 Fix $\ell,w\in k$ and $r\in k^\times$.
Fix a cocompact discrete subgroup $L\subset k$.
 Let $\Phi$ be the unique primitive rigged virtual lattice such that 
 $$\Phi_*=\one_{\ell+L}\otimes \one_{w+r\OO}.$$
Note that $\Phi$ is the general example of a primitive rigged virtual lattice.
Fix an integer $n_0\geq 0$ such that 
$$p^{n_0}\mu(r\OO)/\mu(k/L)\in \ZZ.$$
Fix an integer $M\geq 0$ such that 
$$\norm{\pi^{-M}}\geq \norm{w/r}.$$
 For each $n\in \ZZ$ put
$$L_n=(\pi^{-n}\OO\cap L)\oplus 
\left(\bigoplus_{i=n+1}^\infty \FF_q\cdot \pi^{-i}\right)\subset k,
$$
and let $\Phi_n$ be the unique primitive rigged virtual lattice such that 
$$(\Phi_n)_*=
\one_{\ell+L_n}\otimes \one_{w+r\OO}.$$
By construction $\Phi_n$ is a $\pi$-regular primitive rigged virtual lattice of level $\leq M$ for every $n$. The intuition here is that $\Phi_n$
is a good approximation to $\Phi$ for $n\gg 0$.
The next result makes precise what we mean by good approximation.
\begin{Subproposition}\label{Subproposition:Approximation}
The rigged virtual lattice $p^{n_0}(\Phi-\Phi_n)$ is 
soft and strictly
interpolable  for $n\gg 0$.
\end{Subproposition}
\noindent The proof is structured as a (long) series of lemmas.
\begin{Lemma}\label{Lemma:LastSeparable}
$p^{n_0}(\Phi-\Phi_n)$ is separable for 
$n\gg 0$.
\end{Lemma}
\proof If $n$ is sufficiently large, then $\mu(k/L)=\mu(k/L_n)$, in which case both $p^{n_0}\Phi$ and $p^{n_0}\Phi_n$
are separable by Lemma~\ref{Lemma:RationalRVFT}.
\qed

\begin{Lemma}\label{Lemma:LastSoft}
$\Phi-\Phi_n$ is soft for $n\gg 0$.
\end{Lemma}
\proof
We have
$$(\Phi-\Phi_n)_*(0,t)=
(\one_{\ell+L_n}-\one_{\ell+L})(0)\one_{w+r\OO}(t)$$
for all $t\in k$, and hence $(\Phi-\Phi_n)_*(0,t)$ vanishes identically in $t$ for $n\gg 0$.
By Lemma~\ref{Lemma:RationalRVFT}
we have
$$\GGG_0[(\Phi-\Phi_n)]_*(0,t)=
\ee_0(wt)\mu(r\OO)
\left(\frac{1}{\mu(k/L)}-
\frac{1}{\mu(k/L_n)}\right)
\one_{r^{-1}\DDD^{-1}}(t)
$$
and hence $\GGG_0[(\Phi-\Phi_n)]_*(0,t)$
vanishes identically in $t$ for $n\gg 0$.

It remains only to prove that $\Theta(a,\Phi-\Phi_n)$ vanishes identically in $a$ for $n\gg 0$. Given $n\in \ZZ$,  $a\in k^\times$ and $x\in k$ put
$$F_{n,a}(x)=
\left(\one_{(\ell+L)\cap a(w+r\OO)}-
\one_{(\ell+L_n)\cap a(w+r\OO)}\right)(x),
$$
$$G_{n,a}(x)=
\ee_0(x(aw-\ell))
\left(\one_{L^\perp\cap a^{-1}r^{-1}\DDD^{-1}}-
\frac{\one_{L_n^\perp\cap a^{-1}r^{-1}\DDD^{-1}}}{\mu(k/L_n)/\mu(k/L)}\right)(x).$$
By scaling rule (\ref{equation:RiggedScaling}), we have
\begin{equation}\label{equation:ZerothSeriousFormula}
((\Phi-\Phi_n)^{(1,a)})_*(x,x)=F_{n,a}(x).
\end{equation}
By Lemma~\ref{Lemma:RationalRVFT} and scaling rule (\ref{equation:RationalFourierScalingRule}), we have
\begin{equation}\label{equation:FirstSeriousFormula}
\GGG_0[(\Phi-\Phi_n)^{(1,a)}]_*(x,x)=\frac{\norm{a}\mu(r\OO)}{\mu(k/L)}G_{n,a}(x).
\end{equation}
Therefore by the definition of the theta symbol, along with the scaling rule (\ref{equation:ThetaScaling}) and functional equation 
(\ref{equation:NoughtThetaFE})
satisfied by the theta symbol, we have
$$\Theta(a,\Phi-\Phi_n)=
\Theta(1,(\Phi-\Phi_n)^{(1,a)})=\sum_{x\in k}F_{n,a}(x),
$$
$$\Theta(a,\Phi-\Phi_n)=
\Theta(1,\GGG_0[(\Phi-\Phi_n)^{(1,a)}])=
\frac{\norm{a}\mu(r\OO)}{\mu(k/L)}\sum_{x\in k}G_{n,a}(x).
$$
Now look closely at the formulas for $F_{n,a}$ and $G_{n,a}$.
We can find some $a_1\in k^\times$ and $n_1$ such that
$G_{n,a}$ vanishes identically for $\norm{a}\geq \norm{a_1}$
and $n\geq n_1$.
 We can find some $n_2$
such that $F_{n,a}$ vanishes identically
for $\norm{a}<\norm{a_1}$ and $n\geq n_2$.
It follows that $\Theta(a,\Phi-\Phi_n)$ vanishes identically in $a$
for $n\gg 0$. 
\qed

\begin{Lemma}
For $i\geq n\gg 0$ the rigged virtual lattice $p^{n_0}(\Phi_i-\Phi_n)$
is strictly interpolable and yoked to an interpolating gadget with basepoint $\xi_M$.
\end{Lemma}
\proof It is clear that $p^{n_0}(\Phi_i-\Phi_n)$ is in all cases $\pi$-regular of level $\leq M$.
By the preceding two lemmas,
 $p^{n_0}(\Phi_i-\Phi_n)$ is separable and soft for $i\geq n\gg 0$. Since softness implies properness,
  $p^{n_0}(\Phi_i-\Phi_n)$ is interpolable of
 conductor $\leq \xi_M$ for $i\geq n\gg 0$
by Proposition~\ref{Proposition:piRegularInterpolable}.
Since softness implies effectiveness, $p^{n_0}(\Phi_i-\Phi_n)$
is strictly interpolable and yoked to an interpolating gadget with basepoint $\xi_M$ for $i\geq n\gg 0$ by Proposition~\ref{Proposition:StrictnessCriterion}.   \qed

\begin{Lemma}
Let $U\subset k$ be an open compact subset of $k$.
We have 
$$\card (L^\perp \cap U)=\card (L_n^\perp \cap U)$$
for $n\gg 0$.
\end{Lemma}
\proof
Put $\varphi=\FFF^{-1}[\one_U]$
and select $N$ such that $\varphi$ is supported
in the set $\pi^{-N}\OO$. 
By the Poisson summation formula (\ref{equation:BasicPoisson}) we have
$$ \card ( U\cap L_n^\perp)=
\mu(k/L_n)\sum_{x\in \pi^{-N}\OO\cap L_n}\varphi(x).$$
For $n\gg 0$ the right side does not change if we replace
$L_n$ by $L$.
\qed

\begin{Lemma}
For each $a\in k^\times$, we have
$$\lim_{n\rightarrow\infty}
\left\Vert\left(\begin{array}{c}
a\\
\Phi-\Phi_n
\end{array}\right)-1\right\Vert=0.
$$
\end{Lemma}
\noindent This is the heart of the proof of the proposition.
\proof
By (\ref{equation:ZerothSeriousFormula}), (\ref{equation:FirstSeriousFormula}), the definition of the Catalan symbol,
and the scaling rule (\ref{equation:CatalanScaling}) for the Catalan symbol, we have
$$
\left(\begin{array}{c}
a\\
\Phi-\Phi_n
\end{array}\right)=\left(\begin{array}{c}
1\\
(\Phi-\Phi_n)^{(1,a)}
\end{array}\right)=
\prod_{x\in k^\times} x^{F_{n,a}(x)}\cdot
\left(\prod_{x\in k^\times} x^{G_{n,a}(x)}\right)^{\frac{\norm{a}\mu(r\OO)}{\mu(k/L)}}.
$$
Now the function $F_{n,a}$ vanishes identically for $n\gg 0$ and moreover we can find an annulus
$$A=\{x\in k^\times\mid 0<\norm{a_0}
\leq \norm{x}\leq \norm{a_1}<\infty\}$$ such that $G_{n,a}$
is supported in $A$ for $n\gg 0$. (It is important here that $a$ is fixed; no claim of uniformity of convergence is being made in the lemma.)  Then we have simply
\begin{equation}\label{equation:SecondSeriousFormula}
\left(\begin{array}{c}
a\\
\Phi-\Phi_n
\end{array}\right)^{\frac{\mu(k/L)}{\norm{a}\mu(r\OO)}}=
\prod_{x\in A}x^{G_{n,a}(x)}
\end{equation}
for $n\gg 0$. Note that $G_{n,a}$ is $\ZZ$-valued for $n\gg 0$. Now select any open compact subgroup $U\subset k^\times$ small enough so that the function 
$x\mapsto \ee_0(x(aw-\ell))$ restricted to the compact set $A\subset k^\times$ is constant
on cosets of $U$. Then by the preceding lemma
the summation of $G_{n,a}(x)$ over any coset of $U$
contained in $A$ vanishes for $n\gg 0$, and hence
$\prod_{x\in A}x^{G_{n,a}(x)}\in U$ for $n\gg 0$. 
Since $U$ is arbitrarily small, convergence is proved.
\qed

\proof[Proof of the proposition]
By the preceding lemmas, the hypotheses of \linebreak Proposition~\ref{Proposition:InterpolationLimit} 
 are fulfilled by the family
 $$\{p^{n_0}(\Phi-\Phi_{n})\}\cup\{p^{n_0}(\Phi_{i+n}-\Phi_{n})\}_{i=1}^\infty$$ provided that $n\gg 0$.
\qed
\subsection{End of the proof of Theorem~\ref{Theorem:Interpolability}}
By Proposition~\ref{Subproposition:Approximation}, every proper rational rigged virtual lattice $\Phi$ can be decomposed as a $\ZZ[1/p]$-linear combination of soft strictly interpolable rigged virtual lattices
plus a proper $\pi$-regular rigged virtual lattice which is interpolable by Proposition~\ref{Proposition:piRegularInterpolable}. \qed

 \end{document}